\documentclass{article}

\usepackage[bottom=2.0cm,top=2.0cm,left=2.0cm,right=2.0cm]{geometry}

\usepackage{titlesec}
\usepackage{color}

\usepackage[utf8]{inputenc}
\usepackage[english]{babel}
\usepackage[T1]{fontenc} 

\usepackage{graphicx}
\graphicspath{{Images/}}
\usepackage{eso-pic} 
\usepackage{subfig} 
\usepackage{caption} 
\usepackage{transparent}
\usepackage{float}
\usepackage{tikz}
\usetikzlibrary{fit,positioning,calc,backgrounds, 3d,arrows.meta, perspective}
\usepackage[export]{adjustbox} 

\usepackage{amsmath}
\usepackage{amsthm}
\usepackage{bm}
\usepackage[overload]{empheq}  

\usepackage{tabularx}
\usepackage{longtable} 
\usepackage{colortbl}

\usepackage{algorithm}
\usepackage{algorithmic}

\usepackage{listings}

\usepackage[colorlinks=true,linkcolor=black,anchorcolor=black,citecolor=black,filecolor=black,menucolor=black,runcolor=black,urlcolor=black]{hyperref} 

\usepackage{cleveref}
\usepackage[square, numbers, sort&compress]{natbib} 

\usepackage{appendix}

\usepackage{enumitem}

\usepackage{amsthm,thmtools,xcolor} 
\usepackage{comment} 
\usepackage{fancyhdr} 
\usepackage{lipsum} 
\usepackage{tcolorbox} 
\usepackage{amssymb, dsfont}
\usepackage{siunitx}

\usepackage{authblk}

\newcommand{\R}{\mathbb{R}}
\newcommand{\Rplus}{{\R^{+}}}

\newcommand{\bea}{\begin{eqnarray}} 
\newcommand{\eea}{\end{eqnarray}}

\newcommand{\ol}[1]{\overline{#1}}
\newcommand{\mean}[1]{\{\mskip-6mu\{#1\}\mskip-6mu\}}
\newcommand{\jump}[1]{[\mskip-3mu[#1]\mskip-3mu]}

\newcommand{\DotGrid}[6]{%
  \foreach \x in {0,...,#1}
    \foreach \y in {0,...,#2}
      \fill ({#5 + #3*\x},{#6 + #3*\y}) circle (#4);
}

\title{A hybrid reduced-order and high-fidelity discontinuous Galerkin Spectral Element framework for large-scale PMUT array simulations}
\author[$\S$]{Paola F. Antonietti}
\author[$*$]{Omer M.O. Abdalla}
\author[$\S$]{Michelangelo G. Garroni}
\author[$\S$]{Ilario Mazzieri}
\author[$\S$]{Nicola Parolini}

\affil[$\S$]{\small {MOX-Laboratory for Modeling and Scientific Computing, Department of Mathematics, Politecnico di Milano, Piazza Leondardo da Vinci 32, 20133 Milan, Italy.}}
\affil[$*$]{\small {Department of Civil and
Environmental Engineering, Politecnico di Milano, Piazza Leondardo da Vinci 32, 20133 Milan, Italy.}}

\begin{document}

\everymath{\displaystyle}
\maketitle

\begin{abstract}
Piezoelectric Micromachined Ultrasonic Transducers (PMUTs) are essential for next-generation ultrasonic sensing and imaging due to their bidirectional electromechanical behavior, compact design, and compatibility with low-voltage electronics. As PMUT arrays grow in size and complexity, efficiently modeling their coupled electromechanical-acoustic behavior becomes increasingly challenging. This work presents a novel computational framework that combines model order reduction with a Discontinuous Galerkin Spectral Element Method (DG-SEM) paradigm to simulate large PMUT arrays. Each PMUT’s mechanical behavior is represented using a reduced set of vibration modes, which are coupled to an acoustic domain model to describe the full array. To further improve efficiency, a secondary acoustic domain is connected via DG interfaces, enabling non-conforming mesh refinement, with variable approximation order, and accurate wave propagation. The framework is implemented in the SPectral Elements in Elastodynamics with Discontinuous Galerkin (SPEED) software, an open-source, parallelized platform leveraging domain decomposition, high-order polynomials, METIS graph partitioning, and MPI for scalable performance. The proposed methodology addresses key challenges in meshing, supporting high-fidelity simulations for both PMUT transmission and reception phases. Numerical results demonstrate the framework’s accuracy, scalability, and efficiency for large PMUT array simulations.
\end{abstract}

\section{Introduction}\label{sec:introduction}

Piezoelectricity is the property of certain materials to generate an electric potential when subjected to mechanical stress and, conversely, to deform elastically when an electrical voltage is applied. These phenomena are referred to as the direct and inverse piezoelectric effects, respectively. In Piezoelectric Micromachined Ultrasonic Transducers (PMUTs), this coupled mechanical behaviour enables both ultrasonic wave emission and reception: an applied alternating voltage induces diaphragm vibrations that radiate acoustic waves, while incident pressure waves deform the membrane and generate electrical signals. Thin piezoelectric films such as Lead Zirconate Titanate (PZT), Aluminum Nitride (AlN), and Scandium-doped Aluminum Nitride (ScAlN) are commonly integrated into multilayer diaphragms to achieve efficient electromechanical transduction \cite{MASSIMINO}.

\begin{figure}
    \centering

    \begin{tikzpicture}
    
    \definecolor{colorUSG}{RGB}{230, 230, 250}
    \definecolor{colorAlO}{RGB}{222, 184, 105}
    \definecolor{colorTopElec}{RGB}{255, 165, 0}
    \definecolor{colorPiezo}{RGB}{189, 183, 107}
    \definecolor{colorBotElec}{RGB}{255, 255, 50}
    \definecolor{colorPolySi}{RGB}{240, 230, 170}
    
    \definecolor{colorSubstrate}{RGB}{180, 180, 180}

    \def\fullW{10}      
    \def\activeW{6}    
    \def\cavityW{8}    
    
    \pgfmathsetmacro{\fullStart}{-0.5 * \fullW}
    \pgfmathsetmacro{\cavityStart}{-0.5 * \cavityW}
    \pgfmathsetmacro{\cavityEnd}{0.5 * \cavityW}

    \draw[fill=colorSubstrate] (\fullStart, 0) rectangle (\cavityStart, 2.0);
    \draw[fill=colorSubstrate] (\cavityEnd, 0) rectangle (-\fullStart, 2.0);
    \node[right, black] at (\fullStart+1.1*\cavityEnd,1.0){$\text{Cavity}$};
    
    \draw[->, thick] (5.15,1.0) -- (5.7,1.0);
    \node[right, black] at (5.8,1.0)
    {$\text{Substrate}$};

    \draw[fill=colorPolySi] (\fullStart, 2.0) rectangle (-\fullStart, 3.5);
    \node[right, black] at (\fullStart+1.18*\cavityEnd,2.7){$\text{Si}$};
    
    \draw[fill=colorAlO] (\fullStart, 3.5) rectangle (-\fullStart, 3.9);
    \node[right, black] at (\fullStart+1.14*\cavityEnd,3.68){$\text{AlN}$};
    
    \draw[fill=colorBotElec] (\fullStart+0.4*\cavityEnd, 3.9) rectangle (-\fullStart-0.4*\cavityEnd, 4.1);
    \draw[->, thick] (3.55,5.1) -- (5.0,4.65);
    \draw[->, thick] (3.55,4.0) -- (5.0,4.45);
    \node[right, black] at (5.1,4.525){$\text{Electrodes}$}; 

    \draw[fill=colorPiezo] (\fullStart+0.4*\cavityEnd, 4.1) rectangle (-\fullStart-0.4*\cavityEnd, 5.0);
    \node[right, black] at (\fullStart+1.14*\cavityEnd,4.525){$\text{PZT}$};
    
    \draw[fill=colorTopElec] (\fullStart+0.4*\cavityEnd, 5.0) rectangle (-\fullStart-0.4*\cavityEnd, 5.2);
    
    \draw[fill=colorUSG] (\fullStart+0.4*\cavityEnd, 5.2) rectangle (-\fullStart-0.4*\cavityEnd, 5.57);
    \node[right, black] at (\fullStart+1.14*\cavityEnd,5.38){$\text{SiO}_{2}$};   
    
    \end{tikzpicture}

    \caption{Cross-section of a representative PMUT stack, illustrating the multilayer structure in which each layer has a defined thickness and distinct material properties.}
    \label{fig:pmut_stack}
\end{figure}
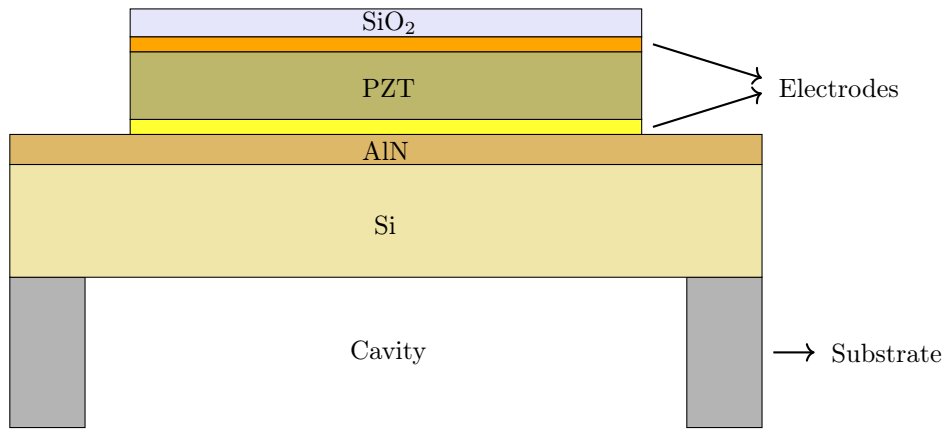

With technological progress, systems increasingly demand higher miniaturization, complexity, and automation.
Ultrasonic devices have followed this trend, particularly through the development of Micromachined Ultrasonic Transducers (MUTs) using Micro-Electro-Mechanical Systems (MEMS) technology, which allow compact, high-performance systems \cite{Birjis2022, He2022}. 
Among them, PMUTs have gained significant attention due to their advantageous electrical properties and compatibility with low-voltage electronics \cite{MASSIMINO}. Their Complementary Metal-Oxide
Semiconductor compatibility enables integration of sensing and control electronics on a single substrate, supporting the development of miniaturized and efficient ultrasonic systems \cite{Birjis2022}. Figure \ref{fig:pmut_stack} provides a representation of a multilayered PMUT stack cross-section.

Compared to bulk piezoelectric transducers, PMUTs benefit from batch microfabrication, resulting in compact arrays with improved yield and reduced cost. In contrast to capacitive MUTs, they do not require high bias voltages and exhibit lower impedance, simplifying electronic interfacing and improving energy efficiency \cite{Birjis2022, He2022}. 

These advantages have enabled PMUT arrays to find applications in medical imaging, fingerprint sensing, sonar, non-destructive testing, flow measurement, and underwater communication \cite{Birjis2022, MASSIMINO, He2022}. 

As PMUT arrays grow in geometric complexity and face increasingly demanding operational constraints, the
numerical models required to accurately simulate their behaviour become significantly more computationally inten
sive. Consequently, there is a pressing need to develop innovative modeling strategies that deliver fast, accurate,
and reliable tools capable of keeping up with the rapid evolution of PMUT technology \cite{ProceedingUtrecht}.

The multiphysics nature of PMUT operation—combining structural deformation, piezoelectric coupling, electrical excitation, and acoustic radiation—makes modeling particularly challenging. Early analytical formulations based on classical plate and membrane theories provided closed-form estimates of resonance frequencies and modal behavior, but these simplified approaches cannot capture multilayer stacks, non-uniform electrode layouts, residual stresses, or fluid loading effects encountered in practical devices \cite{Birjis2022, roy2023thin}.

To overcome these limitations, high-fidelity Finite Element Method (FEM) simulations have become the reference tool for detailed PMUT analysis. Three-dimensional multiphysics FEM models simultaneously solve structural mechanics and piezoelectric equations, enabling accurate prediction of mode shapes, displacement amplitudes, and electromechanical coupling coefficients across geometric and material parameters \cite{he2022sensitivity, sensors2024boltedjoints}. Symmetry-based reductions and parametric sweeps are frequently employed to optimize diaphragm geometry and piezoelectric layer thickness \cite{Sensors2019_4450, PMC5821422}. For realistic performance evaluation, structural models are often coupled with acoustic domains through pressure acoustics or fluid–structure interaction formulations. These coupled simulations allow estimation of radiated pressure fields, radiation impedance, directivity patterns, and bandwidth in air- or liquid-loaded conditions \cite{roy2023thin, nature2025highframe, PMC10673182, Micromachines2022_1979}. Such analyses are crucial for air-coupled sensing and imaging arrays, where acoustic loading significantly affects resonance frequency and damping characteristics. Beyond linear frequency-domain analysis, nonlinear and transient multiphysics simulations have been explored to investigate complex dynamic phenomena. Time-domain coupling of piezoelectric, structural, and fluid equations has been used to model acoustic streaming and micropump actuation in ScAlN-based devices, revealing nonlinear mode interactions and even chaotic behavior under specific excitation regimes \cite{Smyth2015744, ChaosPMUT}. These studies highlight the limitations of purely linear modal approximations when operating at large amplitudes. 

While full 3D multiphysics FEM provides high predictive accuracy, it becomes prohibitively computationally expensive for large parameter sweeps and extended PMUT arrays. To address this issue, reduced-order modeling strategies have been developed. Modal superposition techniques project the electromechanical response onto a truncated basis of dominant eigenmodes extracted from high-fidelity simulations, significantly reducing computational cost while preserving accuracy \cite{Sensors2019_4450, MASSIMINO}. Equivalent circuit models based on Mason or Butterworth–Van Dyke formalisms \cite{CHAUDHARY20222556} offer further simplification by translating electromechanical behavior into electrical analogs suitable for system-level analysis \cite{ijeti2020crosstalk, PMC5821422}. At the array level, modeling complexity increases due to mutual acoustic coupling, mechanical crosstalk, and phase synchronization effects. Hybrid analytical–numerical formulations and reduced-dimensional acoustic models have been proposed to approximate inter-element interactions without explicitly resolving all structural degrees of freedom \cite{ijeti2020crosstalk, SingleChannel3D, arXiv2101_04443, ZHAO2024115677}. Large-domain acoustic simulations typically incorporate absorbing boundary treatments such as perfectly matched layers (PMLs) to prevent artificial reflections \cite{grote2010PML, KALTENBACHERPML, ZHANG2024105765}.

Despite substantial progress, key challenges remain: balancing computational efficiency and predictive accuracy for large arrays, accounting for fabrication variability and residual stresses, capturing nonlinear behavior at high excitation levels, and ensuring stable acoustic–structure coupling in time-domain simulations. These limitations motivate the development of scalable numerical frameworks that combine model reduction with high-performance wave propagation solvers.\\\\
The main objective of this study is to explore computational strategies that can efficiently handle large array
simulations involving an extensive number of degrees of freedom (DOFs). The main idea is to consider an efficient
model order reduction approach derived from \cite{MASSIMINO} to reduce the computational complexity associated with simulating
large PMUT arrays. The proposed novel framework is based on first solving the three-dimensional eigenvalue piezoelectric problem on a
reference PMUT, from which a reduced set of significant vibration modes is extracted. The mechanical response
of each transducer is then expressed as a linear combination of these modal shapes. The modal displacements are
interpolated within a purely acoustic domain that models the PMUT arrays as bi-dimensional membranes. The proposed framework is implemented and validated in SPEED \cite{SPEED}, an open-source software developed for large-scale wave propagation phenomena in complex three-dimensional media. Our framework integrates the Discontinuous Galerkin Spectral Element Method (DG-SEM) with domain decomposition and high-order polynomial approximation to ensure geometric flexibility and high-order accuracy. Mesh partitioning and parallel communication rely on METIS \cite{karypis2013metis} and MPI libraries, while explicit time integration ensures excellent scalability and performance on distributed-memory architectures.

The primary challenges addressed include an optimized mesh generation and partitioning strategy, as well as the optimization of algorithmic operations to ensure both high-performance computing and accuracy in the simulations. By improving these aspects, the proposed approach lead to scalable, high-fidelity modeling of complex PMUT arrays both in transmission (TX) and reception (RX) phases.
\\\\
The remainder of this article is organized as follows. Section \ref{sec:Acoustic problem} formulates the acoustic problem coupled with the modal equations of the PMUT, addressing both the transmission and reception scenarios. This section also includes a discussion on Perfectly Matched Layer (PML) \cite{grote2010PML}, highlighting its role in eliminating spurious reflections at the
boundaries of the computational domain.

Section \ref{sec:numerical_discretization} presents the numerical discretization of the coupled problem in
both space (using the Discontinuous Galerkin Spectral Element Method) and time with Newmark scheme. Section \ref{sec:domain_decomposition} presents the domain decomposition strategy employed for efficient mesh partitioning, which is critical for handling large-scale PMUT arrays. Section \ref{sec:dg_optimization} focuses on the implementation optimization designed to
accelerate computations involving non-conforming meshes. Section \ref{sec:numerical_results}  presents representative numerical results that validate the methodology and demonstrate its performance for large PMUT array simulations, and includes the evaluation of the scalability of the proposed approach on parallel computing architectures. Finally, in Section \ref{sec:conclusion} we draw some conclusions and discuss future perspectives.

\section{Multi-physics problem formulation for large PMUT array simulations}\label{sec:Acoustic problem}
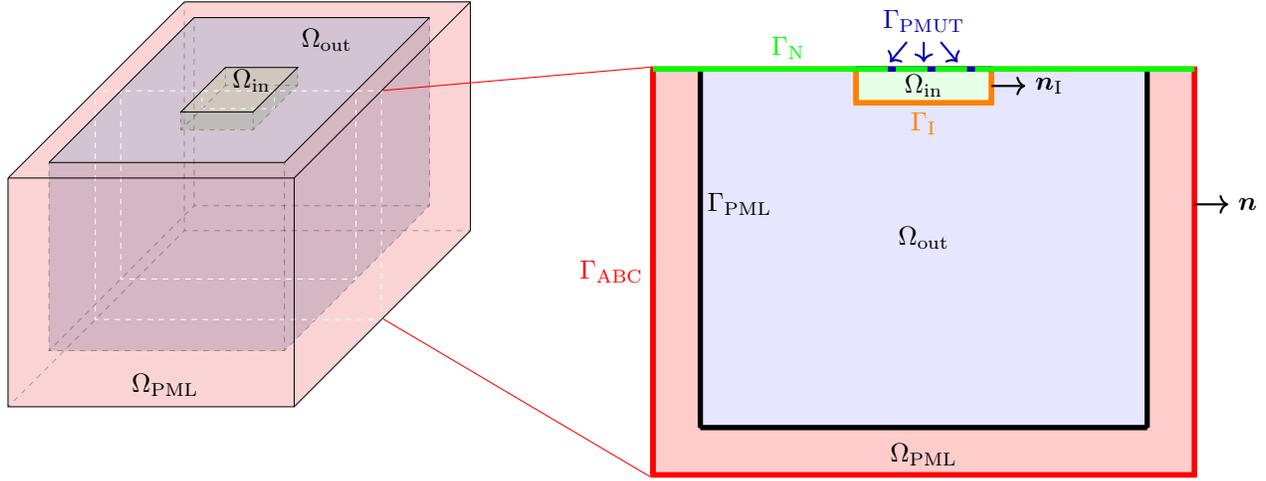
\begin{figure}[H]
\centering
\begin{minipage}{0.35\linewidth}  
    \centering
    \begin{tikzpicture}
        \node[anchor=south west,inner sep=0] (image) at (0,0) 
        {\includegraphics[width=1\linewidth]{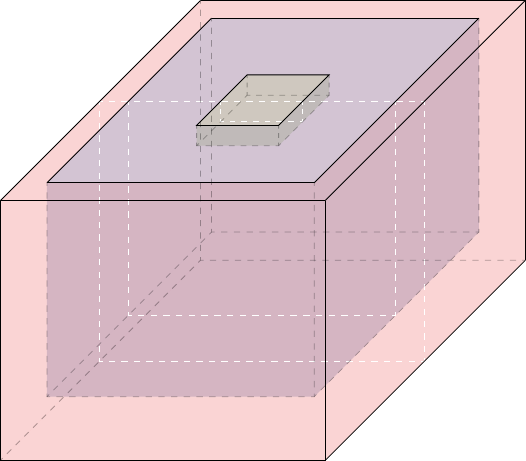}};
        \node at (3.25,4.35) {$\Omega_\text{in}$};
        \node at (4.25,4.9) {$\Omega_\text{out}$};
        \node at (2.1,0.3) {$\Omega_\text{PML}$};

        \draw[-, red] (5,4.22) -- (8.55,4.53);
        \draw[-, red] (5,1.18) -- (8.55,-0.92);
        
    \end{tikzpicture}
    
\end{minipage}
\hspace{0.05\linewidth}
\begin{minipage}{0.55\linewidth}  
    \centering
    \begin{tikzpicture}[scale=0.9, thick]

    \fill[red!20!white] (0,0) rectangle (8,6);
    \node at (4,0.3) {$\Omega_\text{PML}$};

    \fill[blue!10!white] (0.7,0.7) rectangle (7.3,6);
    \node at (4,3.5) {$\Omega_\text{out}$};

    \fill[green!10!white] (3,5.5) rectangle (5,6);
    \node at (4,5.75) {$\Omega_\text{in}$};

    \draw[thick,red,line width=2pt] (0,0) rectangle (8,6);
    \draw[thick] (0.7,0.7) rectangle (7.3,6);
    \draw[thick,orange,line width=2pt] (3,5.5) rectangle (5,6);

    \node[left, red] at (0,3) {$\Gamma_\text{ABC}$};
    \draw[line width=2pt, red] (0,0) -- (8,0);
    \draw[line width=2pt, red] (0,0) -- (0,6);
    \draw[line width=2pt, red] (8,0) -- (8,6);

    \node[left, black] at (1.9,4) {$\Gamma_\text{PML}$};
    \draw[line width=2pt, black] (0.7,0.7) -- (7.3,0.7);
    \draw[line width=2pt, black] (0.7,0.7) -- (0.7,6);
    \draw[line width=2pt, black] (7.3,0.7) -- (7.3,6);

    \node[below, orange] at (4,5.5){$\Gamma_\text{I}$};
    \draw[line width=2pt, orange] (3,5.5) -- (5,5.5);
    \draw[line width=2pt, orange] (3,5.5) -- (3,6);
    \draw[line width=2pt, orange] (5,5.5) -- (5,6);

    \node[above, blue!70!black] at (4,6.4) {$\Gamma_\text{PMUT}$};
    \draw[line width=2pt, blue!70!black] (3,6) -- (5,6);
    \draw[line width=2pt, green, dash pattern=on 12pt off 3pt] (3,6) -- (5,6);

    \node[above, green] at (2,6) {$\Gamma_\text{N}$};
    \draw[line width=2pt, green] (0,6) -- (3,6);
    \draw[line width=2pt, green] (5,6) -- (8,6);

    \draw[->, thick] (5,5.75) -- (5.5,5.75);
    \node[right, black] at (5.5,5.75){$\bm{n}_\text{I}$};

    \draw[->, thick] (8,4) -- (8.5,4);
    \node[right, black] at (8.5,4){$\bm{n}$};

    \draw[->, thick, blue!70!black] (3.77,6.4) -- (3.52,6.1);
    \draw[->, thick, blue!70!black] (4.1,6.4) -- (4.1,6.1);
    \draw[->, thick, blue!70!black] (4.42,6.4) -- (4.67,6.1);

    \end{tikzpicture}
\end{minipage}%

    \caption{A three-dimensional view of the domain decomposition (left) and the corresponding two-dimensional slice with boundary labels (right) for (\ref{eq:full_problem}).}
    \label{fig:domain_label}
\end{figure}

In this section, we introduce the general formulation of the coupled problem.
Let $\Omega=\Omega_\text{in}\cup\Omega_\text{out}\cup\Omega_\text{PML}$ be a three-dimensional open Lipschitz domain, and $T\in\Rplus$ be the final time. The boundary of the acoustic domain $\Omega$ is decomposed as $\partial\Omega=\Gamma_\text{PMUT}\cup\Gamma_\text{N}\cup\Gamma_\text{ABC}$, consisting of disjoint $\Gamma_\text{PMUT}, \Gamma_\text{N}$ and $\Gamma_\text{ABC}$ parts. $\Gamma_\text{PMUT}$ is the portion of the boundary where the action of PMUT membranes is accounted for, $\Gamma_\text{N}$ is where homogeneous Neumann conditions are imposed, and $\Gamma_\text{ABC}$ is the remaining part of the boundary where absorbing condition are considered. 
Additionally, 
$\Gamma_\text{I}$ represents an internal interface between the sub-domain $\Omega_\text{in}$ and $\Omega_\text{out}$, cf. Figure \ref{fig:domain_label}. 
The splitting of $\Omega_{\text{in}}$ and $\Omega_{\text{out}}$ enables the treatment of the two subdomains as potentially non-conforming regions while simultaneously allowing the assignment of distinct material properties within each subdomain. Moreover, the inclusion of the subdomain $\Omega_{\text{PML}}$ ensures the proper implementation of a Perfectly Matched Layer, thereby enabling the effective absorption of outgoing waves and the mitigation of artificial reflections at the truncated computational boundary. 

\noindent We consider the following problem

\begin{equation}
\begin{dcases}
        p_{\text{tt}} + \alpha p_{\text{t}} + \beta p = \nabla \cdot (c^2 \nabla p) + \nabla \cdot \bm{\Phi} - \gamma \psi, &\text{in } \Omega\times\left(0,T\right], \\
        \bm{\Phi}_{\text{t}} = Z_1 \bm{\Phi} + c^2 Z_2 \nabla p + c^2 Z_3 \nabla \psi, &\text{in } \Omega\times\left(0,T\right], \\
        \psi_{\text{t}} = p, &\text{in } \Omega\times\left(0,T\right], \\
        p_{\text{t}} = p=\psi=0, &\text{in } \Omega\times\left\{0\right\},\\
        \bm{\Phi} = \bm{0}, &\text{in } \Omega\times\left\{0\right\}, \\
        \nabla p \cdot \bm{n} = \frac{1}{c}p_{\text{t}}, &\text{on } \Gamma_\text{ABC}\times\left[0,T\right],\\
        \nabla p \cdot \bm{n} = 0, &\text{on } \Gamma_\text{N}\times\left[0,T\right],\\
        \nabla p \cdot \bm{n} = -\kappa\sum_{m}\ddot{q}_{m}\bm{U}_{m}\cdot\bm{n}, &\text{on } \Gamma_\text{PMUT}\times\left[0,T\right],\\
        p_{\lvert_{\Omega_\text{in}}} = p_{\lvert_{\Omega_\text{out}}}, &\text{on } \Gamma_\text{I}\times\left[0,T\right], \\       \left(\nabla p_{\lvert_{\Omega_\text{in}}} - \nabla p_{\lvert_{\Omega_\text{out}}}\right)\cdot\bm{n_\text{I}} =0, &\text{on } \Gamma_{\text{I}}\times\left[0,T\right],
\end{dcases}
\label{eq:full_problem}
\end{equation}
where $p(\bm{x},t)$ denotes the acoustic pressure, $\bm{\Phi}(\bm{x},t)$ represents the (vector) variable that introduces additional diffusion to the problem, and $\psi(\bm{x},t)$ is an auxiliary variable enabling the formulation of the PML damping in a purely local differential form in time. The term $\ddot{q}_{m}$ is computed from the following ordinary differential equation
\begin{equation}
\begin{dcases}
        \ddot{q}_{m} = - \omega_{m}^{2} q_{m} +\int_{\Gamma_\text{PMUT}} p \bm{U}_{m}\cdot\bm{n} \,\mathrm{d\sigma} + \eta_{m} \phi, &\text{in } \left(0,T\right],\\
        \phi = \phi_0(t)\mathds{1}_{\{0 \le t \le \bar{t}\}}(t)+
        \frac{1}{C}\sum_{m}\eta_{m} q_{m}
        \mathds{1}_{\{\bar{t} < t \le T\}}(t),&\text{in } \left[0,T\right],\\
        q_{m} = \dot{q}_{m} = 0, &\text{at } t=0,
\end{dcases}
\label{eq:ode}
\end{equation}
where $\bm{U}_{m}$ encodes the mechanical displacement of the PMUTs with the modal coordinates $q_{m}$, while $\phi(t)$ is a low-voltage activation function. 

The acoustic problem is thus coupled with the piezo-electro-mechanical problem describing the behavior of PMUTs. Due to the large number of membranes typically employed in realistic scenarios, solving the full physics of the PMUTs may be computationally unaffordable. Thus, it is preferable to consider model order reduction techniques, where the behavior of bi-dimensional PMUT membranes is described by ordinary differential equations. Here, we consider the model discussed in \cite{MASSIMINO}, whose aim is to reduce the constitutive law of piezoelectricity by solving an eigenvalue problem to represent the displacement nodal values on a reduced basis of modal vectors. 
Notice that the equations describing the PMUTs' contribution in (\ref{eq:full_problem})-(\ref{eq:ode}) can be expanded as follows:
\begin{equation}
\begin{dcases}
    \ddot{q}_{k,m} = - \omega_{k,m}^{2} q_{k,m} +\int_{\gamma_k} p \bm{U}_{k,m}\cdot\bm{n} \,\mathrm{d\sigma} + \eta_{k,m} \phi_k, &\text{in } \left(0,T\right],\\
    \phi_k = \phi_{k,0}(t)\mathds{1}_{\{0 \le t \le \bar{t}\}}(t)+
    \frac{1}{C}\sum_{m}\eta_{k,m} q_{k,m}
    \mathds{1}_{\{\bar{t} < t \le T \}}(t),&\text{in } \left[0,T\right],\\
    \nabla p \cdot \bm{n} = -\kappa\sum_{k,m}\ddot{q}_{k,m}\bm{U}_{k,m}\cdot\bm{n}, &\text{on } \gamma_k\times\left[0,T\right],\\
    q_{k,m} = \dot{q}_{k,m} = 0, &\text{at } t=0,
\end{dcases}
\label{eq:pmut_problem}
\end{equation}
where $k\in \{1,\dots,N_p\}$, and $m\in \{1,\dots,N_m\}$ represent the $k$-th PMUT and the $m$-th mode, respectively, while $\bigcup\limits_{1\le k\le N_p}\gamma_k=\Gamma_\text{PMUT}$.
Precisely, if $\bm{U}_{k}(\mathbf{x},t)$ denotes the mechanical displacement of the $k$-th PMUT in the reduced-order approach, one has
\begin{equation*}
\bm{U}_{k}(\mathbf{x},t)=
\sum_{m=1}^{N_m} q_{k,m}(t)\,\bm{U}_{k,m}(\mathbf{x}),
\label{eq:modal_projection}
\end{equation*}
where $\bm{U}_{k,m}$ is the $m$-th mode shape, and $q_{k,m}(t)$ is the corresponding generalized modal coordinate. Projecting the mechanical balance equations onto each retained mode, one obtains the ordinary differential equation for $q_{k,m}$ described in (\ref{eq:pmut_problem}). The coefficient $\omega_{k,m}$ is the natural frequency of mode $m$ for the $k$-th PMUT, $p$ is the acoustic pressure on the PMUT’s surface $\gamma_{k}$, and $\eta_{k,m}$ is the electro-mechanical coupling coefficient relating the applied potential $\phi_{k}$ to mode $m$. Due to their nature, PMUTs are able to capture both the transmission (TX) and reception (RX) phases. Setting $\bar{t} \in (0,T)$, for $t \leq \bar{t}$ the voltage $\phi_{k,0}$ is applied to produce an acoustic pressure, while for $t> \bar{t}$ it is directly produced by the pressure on the surface $\gamma_k$. 
The constant $C$ in (\ref{eq:pmut_problem}) models the piezoelectric capacitance. 
For any $k \in \{1,...,N_p \}$, the system \eqref{eq:pmut_problem} represents
the coupling condition at the boundary $\gamma_k$, i.e., at $\Gamma_\text{PMUT}$.
In \eqref{eq:pmut_problem},
 $\kappa = \rho c^{2}$ is the bulk modulus and $\rho=\rho(\bm{x})$ is the acoustic density, $c$ being the sound speed of the medium.

The first three equations of (\ref{eq:full_problem}) describe the acoustic wave propagation and incorporate a perfectly matched layer, i.e. $\Omega_\text{PML}$, which
effectively attenuates nonphysical reflections at the boundary.
The derivation of the PML conditions can be found in \cite{KALTENBACHERPML, grote2010PML}. The scalars $\alpha = (\zeta_1 + \zeta_2 + \zeta_3)$, $\beta = (\zeta_1 \zeta_2 + \zeta_2 \zeta_3 + \zeta_3 \zeta_1)$, $\gamma = \zeta_1 \zeta_2 \zeta_3$, and the matrices
\begin{align}
    Z_1 &= \text{diag}(-\zeta_1, -\zeta_2, -\zeta_3), \nonumber\\
    Z_2 &= \text{diag}(\alpha - 2\zeta_1, \alpha - 2\zeta_2, \alpha - 2\zeta_3), \nonumber\\
    Z_3 &= \text{diag}(\zeta_2 \zeta_3, \zeta_3 \zeta_1, \zeta_1 \zeta_2), \nonumber
\end{align}
are related to the specific damping profile $\zeta_i(x_i), i\in\{1,2,3\}$ selected in the three spatial directions. 
The idea of the PML equations is to provide additional damping and diffusivity to the acoustic equation to shrink the waves' amplitude before they touch the external boundary of the PML, i.e., $\Gamma_\text{ABC}$. In this paper, in agreement with \cite{grote2010PML}, we choose the following damping profile 

\begin{equation*}
    \zeta_i(x_i) = \tilde{\zeta}_i \left(\frac{(|x_i| - h_i)}{\ell_i} - \frac{\sin \left( \frac{2 \pi (|x_i| - h_i)}{\ell_i} \right)}{2\pi} \right)\mathds{1}_{\{h_i \leq |x_i| \leq h_i + \ell_i\}}(x_i),
\end{equation*}
where, for each $i-th$ direction, $i\in{\{1,2,3\}}$, $h_i$ denotes half of the thickness of the propagation domain and $\ell_i$ the thickness of the PML.
The function $\zeta_i(x_i)$ vanishes with continuity at the interface $\Gamma_\text{PML}$ with the purely acoustic domain, and the constant $\tilde{\zeta_i}$ is related to the sound speed $c$ in the acoustic medium, the thickness of the PML, and the relative reflection coefficient $R$ of the acoustic wave at the PML interface, as follows

\begin{equation*}
    \tilde{\zeta_i} = \frac{c}{\ell_i} \log \left( \frac{1}{R} \right), \quad i=1,2,3.
\end{equation*}
The inverse of the reflection coefficient $R$ has to be selected large enough to avoid the wave being reflected at the interface, but it has to be kept small enough to have a less steep function $\zeta_i(x_i)$ in the transition region between the interface $\Gamma_\text{PML}$ and the end of the PML. In the following, $R\in[10^{-5},10^{-4}]$. A first-order absorbing boundary condition is imposed on $\Gamma_\text{ABC}$, while a hard wall boundary condition is set on $\Gamma_\text{N}$. Continuity of the pressure and of its normal derivative is imposed on $\Gamma_\text{I}$.
The terms $\bm{n}$, and $\bm{n}_\text{I}$ represent the unit outer normal to the domain $\Omega$ and to the domain $\Omega_\text{in}$, respectively. The full set of unknowns and parameters defined in equations (\ref{eq:full_problem})-(\ref{eq:ode}) is reported in Table \ref{tab:unk_parameters}.

\begin{table}[h!]
\centering

\begin{minipage}[t]{0.25\textwidth}
\centering

\vspace{0.3em}

\begin{tabular}{|c|c|}
\hline
\textbf{Unknown} & \textbf{Unit} \\
\hline
$p$ & [\si{\pascal}] \\
$\psi$ & [\si{\pascal\second}] \\
$\bm{\Phi}$ & [\si{\pascal\meter\per\second\squared}] \\
$q_m$ & [\si{\meter}] \\
$\phi$ & [\si{\volt}] \\
\hline
\end{tabular}
\end{minipage}
\hfill
\begin{minipage}[t]{0.32\textwidth}
\centering
\vspace{0.3em}

\begin{tabular}{|c|c|}
\hline
\textbf{Acoustic parameter} & \textbf{Unit} \\
\hline
$\alpha, Z_1, Z_2$ & [\si{\per\second}] \\
$\beta, Z_3$ & [\si{\per\second\squared}] \\
$\gamma$ & [\si{\per\second\cubed}] \\
$c$ & [\si{\meter\per\second}] \\
\hline
\end{tabular}
\end{minipage}
\hfill
\begin{minipage}[t]{0.32\textwidth}
\centering

\vspace{0.3em}

\begin{tabular}{|c|c|}
\hline
\textbf{PMUTs parameter} & \textbf{Unit} \\
\hline
$\omega_m$ & [\si{\radian\per\second}] \\
$\eta_m$ & [\si{\coulomb\per\meter}] \\
$\kappa$ & [\si{\pascal}] \\
$C$ & [\si{\farad}] \\
$\phi_0$ & [\si{\volt}] \\
\hline
\end{tabular}
\end{minipage}

\caption{The three tables report the unknowns and the parameters of equations (\ref{eq:full_problem})--(\ref{eq:ode}) with the corresponding units of measurement.}
\label{tab:unk_parameters}
\end{table}

\section{Space discretization: Discontinous Galerkin
Spectral Element  Method}\label{sec:numerical_discretization}
The numerical discretization of (\ref{eq:full_problem}) is based on a Discontinuous Galerkin Spectral Element Method, which allows both for high-order polynomials of SEM and grid flexibility of the DG method \cite{SPEED}. 
%
%

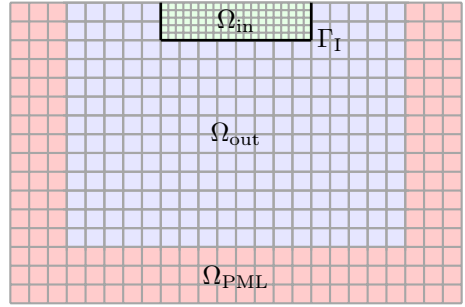
\begin{figure}[H]
    \centering
    
    \begin{tikzpicture}[thick]
    
    \fill[red!20!white] (0,0) rectangle (6,4);
    
    \fill[blue!10!white] (0.75,0.75) rectangle (5.25,4);
    
    \draw[step=0.25, gray!70] (0.7,0.7) grid (5.3,4);
    
    \draw[step=0.25, gray!70] (0,0) grid (6,4);
    
    \fill[green!10!white] (2,3.5) rectangle (4,4);
    
    \draw[step=0.1, gray!70] (2,3.5) grid (4,4);

    \draw[line width=1pt, black] (2,3.5) -- (4,3.5);
    \draw[line width=1pt, black] (4,3.5) -- (4,4);
    \draw[line width=1pt, black] (2,3.5) -- (2,4);
    
    \node at (3,0.35) {$\Omega_\text{PML}$};
    \node at (3,2.25) {$\Omega_\text{out}$};
    \node at (3,3.75) {$\Omega_\text{in}$};
    \node at (4.25,3.5) {$\Gamma_\text{I}$};
    
    \end{tikzpicture}
    \caption{Example of mesh employed in the DG-SE approach. Different mesh sizes and polynomial degrees are employed for $\Omega_\text{in}$ and $\Omega_\text{out}$, resulting in a non-conforming mesh. 
    }\label{fig:dg_mesh}
\end{figure}

    
    
    
    
    
    
    

    
    
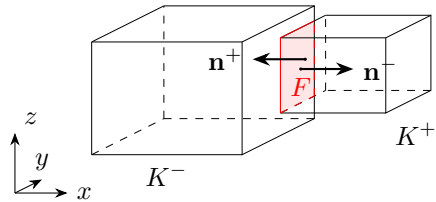
\begin{figure}[H]
\centering
\pgfdeclarelayer{bg}
\pgfsetlayers{bg,main}
\begin{tikzpicture}[
  x={(1cm,0cm)},
  y={(0.6cm,0.3cm)},
  z={(0cm,1cm)},
  >=Stealth,
  box/.style={draw=black,very thin},
  hidden/.style={draw=black,very thin,dashed},
  face/.style={fill=#1!30,opacity=0.7}
]

\newcommand{\drawBox}[6]{%
  \pgfmathsetmacro{\X}{#2}\pgfmathsetmacro{\Y}{#3}\pgfmathsetmacro{\Z}{#4}%
  \coordinate (O) at #1;
  \coordinate (v000) at (O);
  \coordinate (v100) at ($(O)+(\X,0,0)$);
  \coordinate (v010) at ($(O)+(0,\Y,0)$);
  \coordinate (v001) at ($(O)+(0,0,\Z)$);
  \coordinate (v110) at ($(O)+(\X,\Y,0)$);
  \coordinate (v101) at ($(O)+(\X,0,\Z)$);
  \coordinate (v011) at ($(O)+(0,\Y,\Z)$);
  \coordinate (v111) at ($(O)+(\X,\Y,\Z)$);

  
  \draw[hidden] (v110) -- (v010);
  \draw[hidden] (v010) -- (v011);
  \draw[hidden] (v000) -- (v010);
  \draw[#6] (v000) -- (v001);
  \draw[#6] (v000) -- (v100);

  \draw[#6] (v100) -- (v101) -- (v111) -- (v110) -- cycle;
  \draw[#6] (v001) -- (v101);
  \draw[#6] (v011) -- (v111);
  \draw[#6] (v001) -- (v011);
}

\def\Ldx{2.0}   
\def\Ldy{1.6}   
\def\Ldz{1.5}   
\drawBox{(0,0,0)}{\Ldx}{\Ldy}{\Ldz}{white}{box};

\def\Rdx{1.4}
\def\Rdy{1.0}
\def\Rdz{1.0}
\def\RoriginY{0.85} 
\def\RoriginZ{0.3}  
\drawBox{(\Ldx,\RoriginY,\RoriginZ)}{\Rdx}{\Rdy}{\Rdz}{white}{box};

\coordinate (c000) at (\Ldx,\RoriginY,\RoriginZ);                    
\coordinate (c010) at (\Ldx,\RoriginY+\Rdy,\RoriginZ);             
\coordinate (c001) at (\Ldx,\RoriginY,\RoriginZ+\Rdz);             
\coordinate (c011) at (\Ldx,\RoriginY+\Rdy,\RoriginZ+\Rdz);       

\def\ymin{\RoriginY}   
\def\ymax{\RoriginY + 0.75}   

\coordinate (pc000) at (\Ldx, \ymin, \RoriginZ);
\coordinate (pc010) at (\Ldx, \ymax, \RoriginZ);
\coordinate (pc001) at (\Ldx, \ymin, \RoriginZ+\Rdz);
\coordinate (pc011) at (\Ldx, \ymax, \RoriginZ+\Rdz);

\begin{pgfonlayer}{bg}
  \fill[red!10!white]
    (pc000) -- (pc010) -- (pc011) -- (pc001) -- cycle;
\end{pgfonlayer}

\draw[red,dashed] (pc010) -- (pc011);
  \draw[red,dashed] (pc000) -- (pc010);
  \draw[red]        (pc000) -- (pc001);
  \draw[red]        (pc001) -- (pc011);

\coordinate (nP1) at ($(c000)!.45!(c011)$); 
\coordinate (nP2) at ($(c000)!.55!(c011)$); 
\coordinate (fP) at ($(c000)!.45!(c011)$); 
\def\Nlen{0.7}
\draw[->,thick] (nP1) -- ++(\Nlen,0,0) node[right] {$\mathbf{n}^{-}$};
\draw[->,thick] (nP2) -- ++(-\Nlen,0,0) node[left] {$\mathbf{n}^{+}$};
\node[below] at (fP) {$\color{red}{F}$};

\fill[black] (nP1) circle (0.6pt);
\fill[black] (nP2) circle (0.6pt);

\node[below left] at (\Rdx,0,0) {$K^{-}$};
\node[below right] at (\Ldx+\Rdx, \RoriginY, \RoriginZ) {$K^{+}$};

\draw[->] (0, -1.7, 0) -- ++(0.7,0,0) node[right]{$x$};
\draw[->] (0, -1.7, 0) -- ++(0,0.6,0) node[above]{$y$};
\draw[->] (0, -1.7, 0) -- ++(0,0,0.8) node[above right]{$z$};
\end{tikzpicture}
\caption{Representation of two adjacent hexahedra $K^{+}$ and $K^{-}$ with the shared face $F$ and the corresponding normal vectors $\bm{n}^{+}$ and $\bm{n}^{-}$.}
\label{fig:DG_hex}
\end{figure}

We next discuss the derivation of the DG-SEM formulation by assuming, without loss of generality,  that $\Omega_\text{PML}=\emptyset$.
Let $\Omega = \Omega_\text{in}\cup\Omega_\text{out}$ be a polygonal domain and consider its partition $\mathcal{T}_h = \mathcal{T}_\text{in}\cup\mathcal{T}_\text{out}$, s.t. $\mathcal{T}_\text{in}\cap\mathcal{T}_\text{out}=\emptyset$ with $\mathcal{T}_\text{in}$ and $\mathcal{T}_\text{out}$ containing hexahedral elements $K_\text{in}$ and $K_\text{out}$ of different granularity $h_\text{in}$ and $h_\text{out}$, respectively. By defining  $N_\text{in}$, and $N_\text{out}$, 
the number of elements in each partitions, one has $\mathcal{T}_{\text{in}}=\left\{K_\text{in}^j \right\}_{j=1}^{N_\text{in}}$ and $\mathcal{T}_{\text{out}} = \left\{ K_\text{out}^j \right\}_{j=1}^{N_\text{out}}$, respectively.
%
Let $\hat{K}=[-1, 1]^3$ be the reference element, and let $\bm{X}_{K_\text{a}}:\hat{K}\rightarrow K_\text{a}$ be a trilinear invertible 
map with positive Jacobian $\bm{J}_{K_\text{a}}$ such that $ \forall\,K_\text{a}\in\mathcal{T}_\text{a}$,  $K_\text{a} = \bm{X}_{K_\text{a}}(\hat{K})$, for $\text{a} = \{ \text{in}, \text{out} \}$, with $\bm{X}_{K_\text{a}}^{-1}:K_\text{a}\rightarrow \hat{K}$ . Next, we define the internal faces as
$F = \ol{K}^{j}_\text{in}\cap\ol{K}^{k}_\text{out}$ for certain $j\in\{1,...,N_\text{in}\}$ and $k \in \{1,...,N_\text{out}\}$ and collect them in $\mathcal{F}^I$.
Similarly, we define the boundary faces 
$F = \ol{K}^{j_\text{a}}_\text{a}\cap\ol{\partial\Omega}$ for certain $j_\text{a}\in\{1,...,N_\text{a}\}$ and collect them in $\mathcal{F}^B$.
Moreover, we split $\mathcal{F}^B$ as $\mathcal{F}^B=\mathcal{F}_\text{N}^B\cup\mathcal{F}_\text{ABC}^B\cup\mathcal{F}_\text{PMUT}^B$, where Neumann, absorbing, and coupling conditions (with PMUT membranes) are imposed, respectively. Note that boundary faces are compliant with this subdivision.
%
%
%
%
%
%
\\Consider now two adjacent elements $K^-\in\mathcal{T}_{\text{out}}$ and $K^+\in\mathcal{T}_{\text{in}}$, sharing a face $F\in \mathcal{F}^I$ and define $\bm{n}^+$ and $\bm{n}^-$  to be the outward unit normal vectors to the face $F$ from $K^+$ and $K^-$, respectively. Figure \ref{fig:DG_hex} provides a visual representation of two adjacent elements.
\noindent
For each face $F\in \mathcal{F}^I$ and for any sufficiently regular scalar function $u$ and vector $\bm{v}$ 
 we define the jump operator $\jump{\cdot}$ and the average operator $\mean{\cdot}$ by
\begin{alignat*}{2}
    & \mean{u} = \frac{1}{2}(u^{+} + u^{-}), \quad && \jump{u} = u^{+}\bm{n}^{+} + u^{-}\bm{n}^{-},
\end{alignat*}
\begin{alignat*}{2}
    & \mean{\bm{v}} = \frac{1}{2}(\bm{v}^{+} + \bm{v}^{-}), \quad && \jump{\bm{v}} = \bm{v}^{+}\cdot\bm{n}^{+} + \bm{u}^{-}\cdot\bm{n}^{-},
\end{alignat*}
where $\pm$ denotes the trace on $u$ and $\bm{v}$ taken within the interior of  $K^\pm$, respectively. Similarly, for each face $F\in \mathcal{F}^B$ we set  $\mean{u} = u$, $\jump{u} = u\bm{n}$ and $\mean{\bm{v}} = \bm{v}$, $\jump{\bm{v}} = \bm{v}\cdot\bm{n}$.

\noindent For $\text{a} = \{ \text{in}, \text{out} \}$  we consider the space
\begin{equation*}
    V_\text{a}^{r_\text{a}} = \left\{ v\in \mathcal{C}^0\left(\overline{\Omega}_\text{a}\right): v\big\rvert_{K_\text{a}}\circ\bm{X}^{-1}_{K_\text{a}}\in\mathbb{Q}_{r_\text{a}}(\hat{K}),\quad\forall K_\text{a}\in\mathcal{T}_\text{a} \right\},
\end{equation*}
where $\mathbb{Q}_{r_\text{a}}(\hat{K})$ is the space of polynomials of degree less than or equal
to $r_\text{a} \ge 1$ in each coordinate direction. Next, we introduce 
\begin{equation*}
    V^{\text{DG}} =  \left\{ v\in L^2\left(\Omega\right): v\big\rvert_{{\Omega}_{\text{a}}}\in V^{r_\text{a}}_{\text{a}} \right\},
\end{equation*}
and provide the semi-discrete formulation:
$\forall\, t\in \left(0,T\right]$ find $p_h = p_h(t) \in V^{\text{DG}}$ s.t. 
\begin{equation}
\label{eq:dg}
    m(\ddot{p_h}, v_h) + a(p_h, v_h) + d(\dot{p}_h, v_h)=  -\kappa\sum_{m=1}^{N_m}\int_{\Gamma_\text{PMUT}}\ddot{q}_{m}\bm{U}_{m}\cdot\bm{n}v_h\,\mathrm{d\sigma},\quad\forall\, v_h\in V^\text{DG},
\end{equation}
with $p_h(0) = \dot{p}_h(0) = 0$ and where, for any $1 \leq m  \leq N_m$, it holds 
\begin{equation}\label{eq:pmuts_eqs}
    \ddot{q}_{m}(t) + \omega_{m}^{2} q_{m}(t) =\int_{\Gamma_\text{PMUT}} p_h(t) \bm{U}_m\cdot\bm{n} \,\mathrm{d\sigma} + \eta_{m} \phi(t) \quad\forall\, t\in \left(0,T\right].
\end{equation}
For any $\omega,v\in V^\mathrm{DG}$, the bilinear forms in  (\ref{eq:dg}) are defined as:

\begin{align}
m(\omega, v)
  &= 
\sum_{K_{\text{s}}\in\mathcal{T}_{\text{in}} \cup \mathcal{T}_{\text{out}}}
     \!\int_{K_{\text{s}}} \omega v\, \mathrm{d}\Omega, \label{eq:Wf_1}
\\[6pt]
a(\omega, v)
  &= \sum_{K_{\text{s}}\in\mathcal{T}_{\text{in}} \cup \mathcal{T}_{\text{out}}}
\int_{K_{\text{s}}}c^{2}\nabla\omega\cdot\nabla v\,\mathrm{d}\Omega 
-\sum_{F\in\mathcal{F}^{I}}\!\int_{F} \mean{c^{2}\nabla\omega}\cdot\jump{v}\,\mathrm{d}\sigma \nonumber \\ & \qquad \qquad\qquad \qquad\qquad \qquad -\sum_{F\in\mathcal{F}^{I}}\!\int_{F} \mean{c^{2}\nabla v}\cdot\jump{\omega}\,\mathrm{d}\sigma   +\sum_{F\in\mathcal{F}^{I}}\!\int_{F} \gamma_{_F}\,\jump{\omega}\cdot\jump{v}\,\mathrm{d}\sigma,
\label{eq:Wf_2}
\\[6pt]
d(\omega, v)
  &= \sum_{F\in\mathcal{F}^{\mathrm{B}}_{\mathrm{ABC}}}
     \int_{F}c \omega v\,\mathrm{d}\sigma,\label{eq:Wf_4}
\end{align}
where the parameter $\gamma_{_F} = \alpha \bar{c}^2\frac{\max({r_\text{in}},{r_\text{out}})^2}{\min(h_\text{in},h_\text{out})}$, with $\bar{c}$ being the harmonic average of $c$ on $F$
and with $\alpha$ a sufficiently large positive parameter at our disposal, see, e.g., \cite{AntoniettiArtoniMazzieriParoliniRocchi2023}.

\subsection{Algebraic formulation and time integration}
Let $\{\psi_{i}\}_{i=1}^{N_h}$ be the Lagrangian basis functions for $V^\text{DG}$, cf. \cite{MazRapQA}, and set $v_h(\bm x)=\psi_i(\bm x)$, and $p_h(\bm{x},t)=\sum_{j=1}^{N_h}p_j(t)\psi_{j}(\bm{x})$ in 
(\ref{eq:Wf_1})--(\ref{eq:Wf_4}). Rearranging the terms in matrix form, the semi-discrete algebraic formulation of problem (\ref{eq:dg})-(\ref{eq:pmuts_eqs}) reads: for any $t\in (0,T]$ find $\bm p(t) = (p_1(t),...,p_{N_h}(t))^T \in \mathbb{R}^{N_h}$ and $\bm q(t) = (q_1(t),...,q_m(t))^T \in \mathbb{R}^{N_m}$,  s.t.
\begin{align}
M\bm{\ddot{p}}(t) + K\bm{p}(t) + C\bm{\dot{p}}(t) & = \bm{f}(\ddot{\bm q}(t)),\label{eq:semi_SEM}\\
\ddot{\bm q}(t) +  W \bm q(t)  & = \bm g (\bm p(t), \phi(t)), 
\label{eq:semi_PMUT}
\end{align}
with $\bm p = \dot{\bm p} = \bm 0$, $\bm q = \dot{\bm q} = \bm 0$. In \eqref{eq:semi_SEM},  $M\in\mathbb{R}^{N_h\times N_h}$ is the diagonal mass matrix associated with $m\left(\cdot, \cdot \right)$,
i.e., the matrix representation of the $L^2$ inner product, evaluated with reduced integration, $K\in\mathbb{R}^{N_h\times N_h}$ is the stiffness matrix related to $a\left(\cdot, \cdot \right)$, $C\in\mathbb{R}^{N_h\times N_h}$ is the matrix encoding the absorbing conditions $d\left(\cdot,\cdot\right)$, while $\bm{f}\in\mathbb{R}^{N_h}$ is the vector describing 
the forcing term acting on $\Gamma_\text{PMUT}$, i.e.,
\begin{equation}\label{eq:f_i}
\bm f_i(t) = -\kappa\sum_{m=1}^{N_m}\ddot{q}_{m}(t)\int_{\Gamma_\text{PMUT}}\bm{U}_{m}(t,\bm x)\cdot\bm{n}\psi_i(\bm x)\,\mathrm{d\sigma}, \quad 1\le i\le N_h.
\end{equation}
\\In \eqref{eq:semi_PMUT}, $W \in \mathbb{R}^{N_m \times N_m}$ is the matrix $W = {\rm diag}(\omega_{1}^{2},...,\omega_{N_m}^{2})$ and $\bm g \in \mathbb{R}^{N_m}$ is such that 
\begin{equation}\label{eq:g_m}
g_m(t) = \sum_{j=1}^{N_h}p_j(t) \int_{\Gamma_\text{PMUT}} \psi_{j}(\bm{x}) \bm{U}_{m}(\bm x)\cdot\bm{n} \,\mathrm{d\sigma} + \eta_{m} \phi(t) \quad 1 \le m \le N_m.
\end{equation}
\\Next,  we let $S\in\mathbb{N^+}$ be the number of sub-intervals of length $\Delta t$ in which the time interval $[0,T]$ is divided in such a way that $t^n = n\Delta t, n\in \{0,\dots, S\}$. The equations (\ref{eq:semi_SEM})-(\ref{eq:semi_PMUT}) can be discretized with the following second-order predictor-corrector  Newmark scheme, 
shown in Algorithm \ref{alg:newmark}. 
Note that a similar scheme is employed when a PML boundary is introduced in the computational model.

The Newmark time-integration method is employed with $\gamma=1/2$ and $\beta=0$. For this choice, the scheme is second-order accurate, and the corresponding update relations reduce to the explicit leapfrog formulation. Since the method is explicit, it is conditionally stable under the CFL condition $\Delta t\lesssim h_\text{min}$, being $h_\text{min}$ the minimum diameter among the mesh elements \cite{raviart1983introduction}.

\begin{algorithm}[H]
\caption{Newmark time integration scheme}
\label{alg:newmark}
\begin{algorithmic}[1]
\STATE Set $\bm q^0 = \dot{\bm q}^0 =  \ddot{\bm q}^0 = \bm 0$ and 
$\bm{p}^0 = \dot{\bm{p}}^0 =  \ddot{\bm{p}}^0 = \bm 0$
and compute the voltage
$\phi^0$
\FOR{$n=0,...,S$} 
    \vspace{1.5mm}\STATE Predictor (\textit{PMUT}): set $\bm q^{n+1} \leftarrow \bm q^{n} + \Delta t\, \dot{\bm q}^{n} + \tfrac{1}{2}\Delta t^2\, \ddot{\bm q}^{n}$
    \STATE Predictor (\textit{PMUT}): set $\dot{\bm q}^{n+1} \leftarrow \dot{\bm q}^n + \tfrac{1}{2}\Delta t\, \ddot{\bm q}^{n}$ 
    \vspace{1.5mm}\STATE Solve (\textit{PMUT}): \quad \quad \, \quad  $\ddot{\bm q}^{n+1} \leftarrow - W \bm q^{n+1} + \bm g(\bm{p}^n, \phi^n)$
    \vspace{1.5mm}\STATE Corrector (\textit{PMUT}): \quad \; $\dot{\bm q}^{n+1} \leftarrow \dot{\bm q}^{n+1} + \tfrac{1}{2}\Delta t \, \ddot{\bm q}^{n+1}$
    %
%
    \vspace{1.5mm}\STATE Predictor (\textit{acoustic}): set  $\bm{p}^{n+1} \leftarrow \bm{p}^n + \Delta t\, \dot{\bm{p}}^n + \tfrac{1}{2}\Delta t^2\, \ddot{\bm{p}}^n,$
    \STATE Predictor (\textit{acoustic}): set $\dot{\bm{p}}^{n+1} \leftarrow \dot{\bm{p}}^n + \tfrac{1}{2}\Delta t\, \ddot{\bm{p}}^n$ 
    \vspace{1.5mm}\STATE Solve (\textit{acoustic}): \quad \quad \quad \;$\ddot{\bm{p}}^{n+1} \leftarrow M^{-1}(\bm{f}(\ddot{\bm q}^{n+1})- K\bm{p}^{n+1} - C\dot{\bm{p}}^{n+1})$
    \vspace{1.5mm}\STATE Corrector (\textit{acoustic}): \;  \, $\dot{\bm{p}}^{n+1}  \leftarrow \dot{\bm{p}}^{n+1} + \tfrac{1}{2}\Delta t\, \ddot{\bm{p}}^{n+1}$
    \vspace{1.5mm}\STATE Update: 
    $\bm q^n \leftarrow \bm q^{n+1}, \quad  \dot{\bm q}^n \leftarrow \dot{\bm q}^{n+1}, \quad  \ddot{\bm q}^n \leftarrow \ddot{\bm q}^{n+1}$
    \STATE Update: $\bm{p}^n \leftarrow \bm{p}^{n+1}, \quad   \dot{\bm{p}}^n \leftarrow \dot{\bm{p}}^{n+1}, \quad  \ddot{\bm{p}}^n \leftarrow \ddot{\bm{p}}^{n+1} $
\vspace{1.5mm}
\ENDFOR
\end{algorithmic}
\end{algorithm}

\section{Parallel implementation strategies for PMUT simulations}\label{sec:domain_decomposition}
In the context of High-Performance Computing, achieving competitive runtimes for large-scale simulations is of paramount importance. This requires fully leveraging the computational capabilities of parallel codes on multi-core clusters. A key challenge in this context is the efficient distribution of computational workload across cores, which often involves adopting domain-partitioning strategies tailored to the specific requirements of the simulation. In this work, we address a particularly challenging load-balancing problem arising from the need to account for multiple constraints simultaneously. Specifically, our goal is to balance the computational cost associated with the
\begin{enumerate}[label=\roman*)]
    \item boundary conditions for the PMUTs-array;
    \item acoustic propagation in $\Omega_\text{out}\cup\Omega_\text{PML}$;
    \item  connectivity setup for non matching grids on  $\Gamma_\text{I}$.
\end{enumerate}
Additionally, we aim to ensure that the membranes representing the PMUTs, namely the boundary $\Gamma_\text{PMUT}$, are not split across multiple cores, as such fragmentation can degrade both performance and efficiency due to increased inter-process communication and loss of structural coherence. To this end, we propose a tailored partitioning strategy that respects the structural integrity of the PMUTs while distributing the simulation load evenly across processing units. This approach enhances the scalability of the simulation framework and ensures efficient utilization of HPC resources without compromising the accuracy of the underlying discretization.

\subsection{Setup of the computational domain structure}
%
The computational domain considered in the present case study consists of three main regions. The first one, which is referred to as \textit{inner domain}, $\Omega_\text{in}$, contains the PMUT array and constitutes the region of primary geometrical and physical interest. The second one, which is referred to as \textit{outer domain}, $\Omega_\text{out}$, represents the acoustic wave propagation region. The third one, $\Omega_\text{PML}$, describes the external diffusive layer where Perfectly Matched Layer conditions are applied in order to mimic an unbounded medium.

For geometric fidelity and accuracy requirements associated with the vibrating membranes, the grid in $\Omega_\text{in}$ must be highly refined. Imposing a mesh of comparable resolution throughout $\Omega_\text{out}$ would lead to a prohibitive number of elements—potentially on the order of billions—thus rendering the computation intractable with standard computational resources. To reduce the overall number of degrees of freedom while preserving the required accuracy near the membranes, a coarser discretization is therefore adopted in $\Omega_\text{out}$.

Accordingly, the inner and the outer subdomains are meshed independently, allowing for local refinement in $\Omega_\text{in}$ and a coarser mesh in the propagation region $\Omega_\text{out}$. Conformity is instead imposed between $\Omega_\text{out}$ and $\Omega_\text{PML}$ to ensure a consistent treatment of the absorbing layer. The use of distinct discretizations in $\Omega_\text{in}$ and $\Omega_\text{out}$ results in a globally nonconforming mesh. To properly couple the two regions, an internal interface $\Gamma_\text{I}$ is introduced between $\Omega_\text{in}$ and $\Omega_\text{out}$, enabling the matching of the finer and coarser meshes across the subdomain boundary.

\subsubsection{Load balancing strategy for  \texorpdfstring{$\Omega_\text{out} \cup \Omega_\text{PML}$}{Omegaout U OmegaPML}
} 

\begin{figure}[H]
\begin{minipage}{0.23\linewidth} 
\centering
\begin{tikzpicture}[scale=1]
    \node[anchor=south west, inner sep=0] (img) at (0,0)
    {\includegraphics[width=\linewidth, angle=270]
    {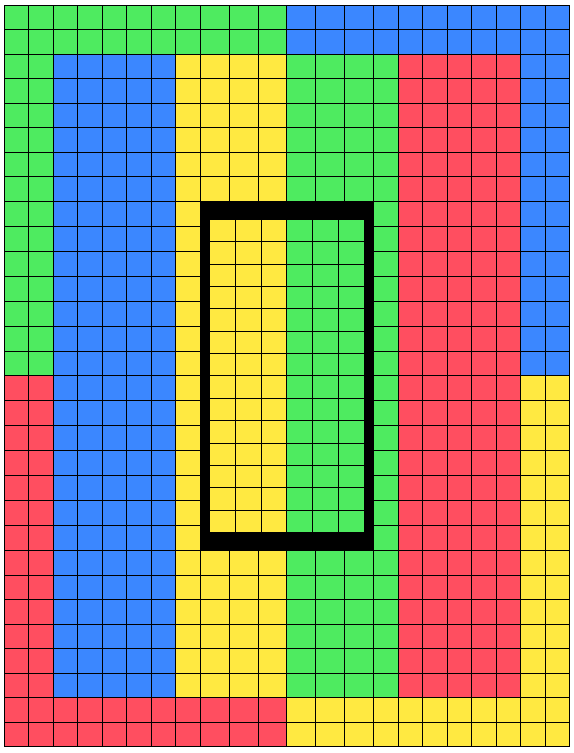}};              
          
\end{tikzpicture}
\end{minipage}%
\hspace{5.5em}
\begin{minipage}{0.23\linewidth}
\centering
\begin{tikzpicture}[scale=1]
    \node[anchor=south west, inner sep=0] (img) at (0,0)
    {\includegraphics[width=\linewidth, angle=90]
    {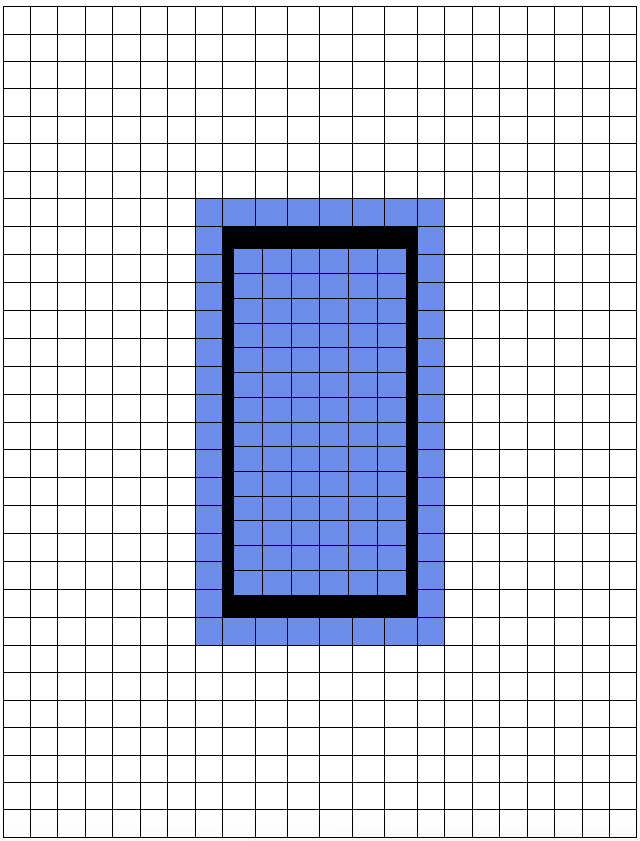}};
    \begin{scope}[x={(img.south east)}, y={(img.north west)}]
        \node[below left] at (0.55,0.55) {$\Omega_\text{DG}$};
        \node[below left] at (0.7,0.2)        {$\Omega_\text{out}\cup\Omega_\text{PML}\setminus\Omega_\text{DG}$};
    \end{scope}
\end{tikzpicture}
\end{minipage}%
\hspace{5.5em}
\begin{minipage}{0.23\linewidth}
\centering
\begin{tikzpicture}[scale=1]
    \node[anchor=south west, inner sep=0] (img) at (0,0)
    {\includegraphics[width=\linewidth, angle=90]
    {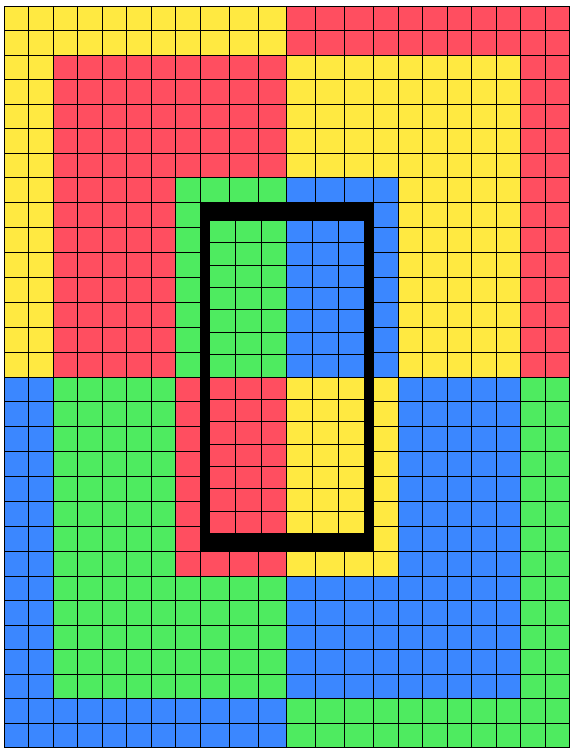}};

\end{tikzpicture}
\end{minipage}

\definecolor{c1}{rgb}{0.23137, 0.52941, 1.00000}
\definecolor{c2}{rgb}{1.00000, 0.91373, 0.25490}
\definecolor{c3}{rgb}{1.00000, 0.30588, 0.37647}
\definecolor{c4}{rgb}{0.30588, 0.92157, 0.37647}

\caption{Schematic illustration of the outer domain decomposition: 
non-optimized partitioning (left), DG layer $\Omega_\text{DG}$ (middle), optimal partitioning (right). 
The set of processors is 
$\{P_1, P_2, P_3, P_4\} = 
\{\textcolor{c1}{\raisebox{0.3ex}{$\bullet$}},\,
 \textcolor{c2}{\raisebox{0.3ex}{$\bullet$}},\,
 \textcolor{c3}{\raisebox{0.3ex}{$\bullet$}},\,
 \textcolor{c4}{\raisebox{0.3ex}{$\bullet$}}\}$.}

\label{fig:outer_domain_three_sided}
\end{figure}

The two domains $\Omega_\text{out}$ and $\Omega_\text{PML}$ are partitioned directly using Metis 5.1.0 \cite{karypis2013metis}.
The idea is to exploit multiple times a Metis routine, namely \texttt{METIS\_PartMeshDual}, that receives as input the set of hexahedra composing the two regions and returns the ownership of the mesh elements to the computational cores. 
The routine can be directly applied to the full set of mesh elements composing $\Omega_\text{out}$ and $\Omega_\text{PML}$ if the primary interest is to balance the number of hexahedra among the unit processes. However, in most cases, we will obtain an unbalanced workload 
for the processors, since not all processors will compute the integral terms on $\mathcal{F}^I$ in \eqref{eq:Wf_2}, cf. Figure \ref{fig:outer_domain_three_sided} (left).

To overcome this bottleneck, we first split $\Omega_{\rm out}$ into two subdomains: $\Omega_{DG}$ containing the hexahedral elements that share a face with $\Gamma_\text{I}$ and the remaining ones. 
Next, to equally distribute elements over processors, we call three times \texttt{METIS\_PartMeshDual}: one for $\Omega_{\rm PML}$, one for $\Omega_{\rm DG}$ and one for the remaining part $\Omega_{\rm out}\setminus \Omega_{\rm DG}$.
%
Figure \ref{fig:outer_domain_three_sided} shows the $\Omega_\text{DG}$ layer connected to the outer domain $\Omega_\text{out}$ (middle) and a schematic representation of the Metis partition (right).

\subsubsection{Mesh generation and load balancing strategy for \texorpdfstring{$\Omega_{\rm in}$}{Omegain}}
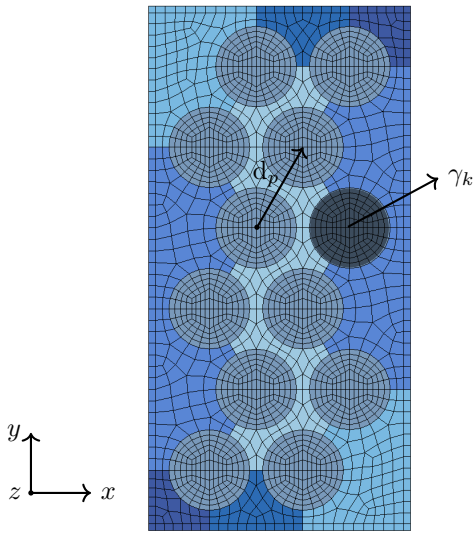
\begin{figure}[H]
\centering
\begin{tikzpicture}
  \node[inner sep=0pt] (img) at (4.5,0.8)
    {\includegraphics[angle=0, width=0.2\linewidth]{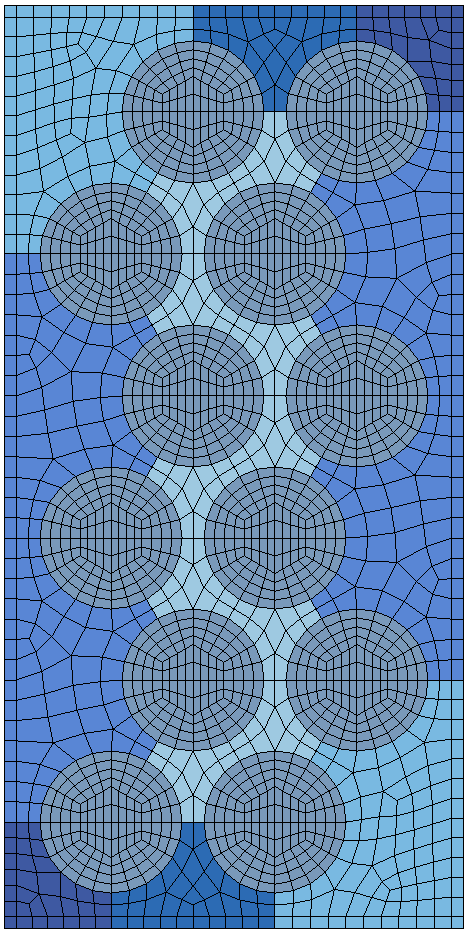}};

  \coordinate (node_z) at (1.2, -2.20);
  \fill[black] (node_z) circle (1pt);
  \draw[->, thick, black]  (node_z) -- +(0.8, 0) node[right] {$x$};
  \draw[->, thick, black](node_z) -- +(0, 0.8) node[left] {$y$};  
  \node[left] at (node_z) {\color{black}$z$};

  \draw[->, thick, black]  (4.215,1.35) -- +(0.605,1.05) node[] at (4.33,2.05) {$\text{d}_p$};
  \fill[black] (4.205,1.335) circle (1pt);
   \fill[black, opacity=0.5] (5.435,1.325) circle (14.9pt);
   \draw[->, thick, black]  (5.415,1.34) -- +(1.2,0.65) node[above, right] {$\gamma_k$};
\end{tikzpicture}
\caption{Top view of a minimal example of inner domain $\Omega_\text{in}$. The elementary structures originating $\Omega_\text{in}$ are depicted with different colors. The circular regions represent the membranes $\gamma_k$ belonging to $\Gamma_\text{PMUT}$, while the remaining part is $\Gamma_\text{N}$.}
\label{fig:minimal_inner_domain}
\end{figure}

The concept behind the mesh generation for the inner domain is to define a set of fundamental unit blocks that can be combined to form the complex membrane pattern of the PMUT array. In general, the membranes are not arranged in a Cartesian $x$-$y$ grid but are only aligned along the $x$-direction, as shown in Figure \ref{fig:minimal_inner_domain}. If we define \textit{row} the set of membranes aligned along $x$ and \textit{column} the set of PMUTs that winds along $y$, then each row in the array is offset relative to the previous one by a distance determined by the pitch of the PMUT, $\text{d}_p$, that is, the distance of the centers of two neighboring membranes, see Figure \ref{fig:minimal_inner_domain}.

\begin{figure}[H]
    \centering
        \begin{minipage}{0.4\linewidth}
        \centering
        \begin{tikzpicture}[scale=1]
            \node[opacity=0.6, anchor=south west, inner sep=0] (img) at (0,0) {\includegraphics[width=1\linewidth, angle=0]{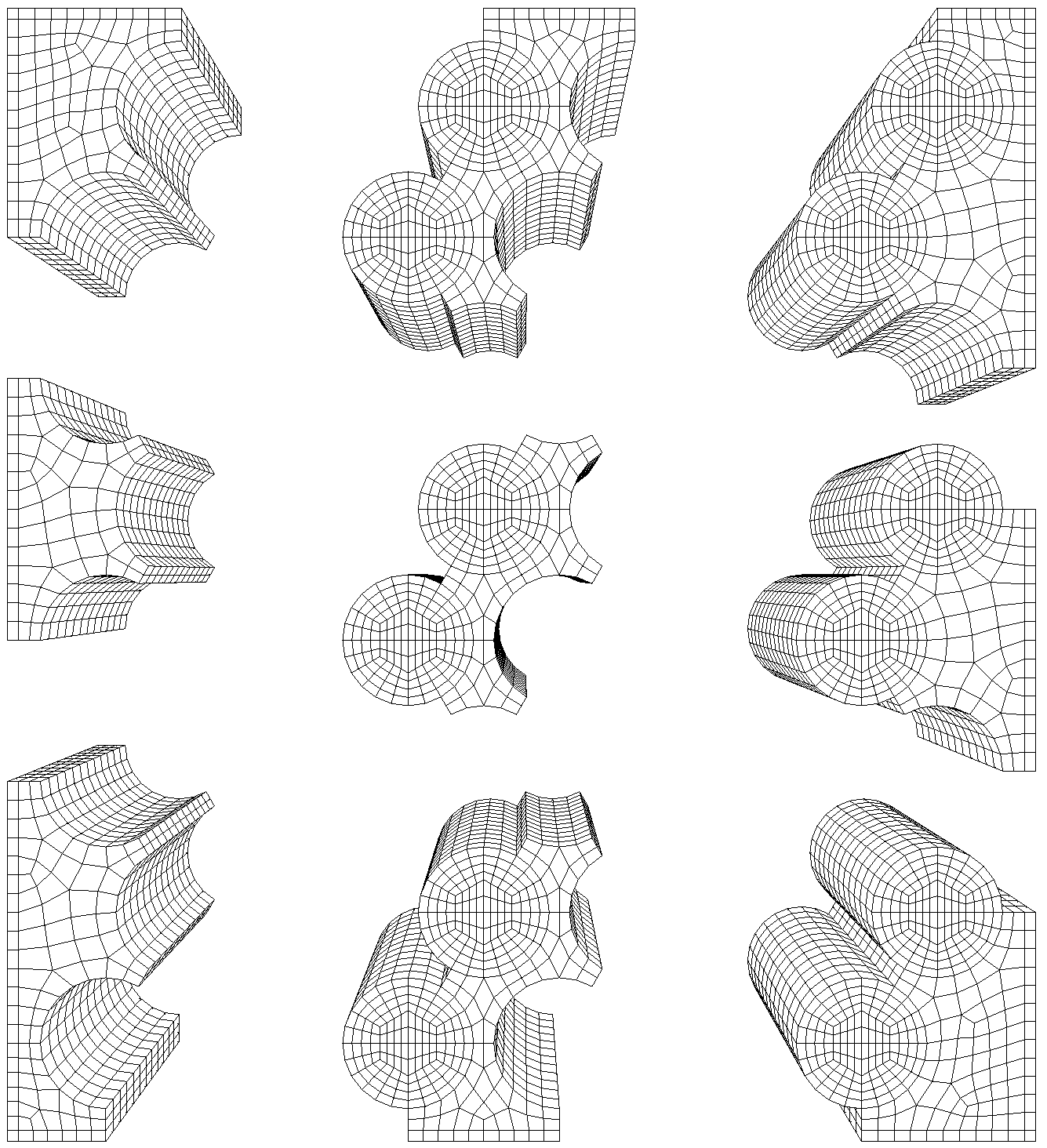}};
          \DotGrid{2}{0}{0.125}{0.75pt}{1.7}{6.4}
          \DotGrid{2}{0}{0.125}{0.75pt}{1.7}{4.1}
          \DotGrid{2}{0}{0.125}{0.75pt}{1.7}{1.7}
          \DotGrid{2}{0}{0.125}{0.75pt}{4.4}{6.4}
          \DotGrid{2}{0}{0.125}{0.75pt}{4.4}{4.1}
          \DotGrid{2}{0}{0.125}{0.75pt}{4.4}{1.7}
    
          \DotGrid{0}{2}{0.125}{0.75pt}{0.5}{5.25}
          \DotGrid{0}{2}{0.125}{0.75pt}{0.5}{2.95}
          \DotGrid{0}{2}{0.125}{0.75pt}{3.125}{4.85}
          \DotGrid{0}{2}{0.125}{0.75pt}{3.125}{2.55}
          \DotGrid{0}{2}{0.125}{0.75pt}{5.75}{4.85}
          \DotGrid{0}{2}{0.125}{0.75pt}{5.75}{2.55}

          \node[] at (0.5,1.5) {$\bm{B_{1,1}}$};
          \node[] at (0.5,4.3) {$\bm{B_{k,1}}$};
          \node[] at (0.5,7.0) {$\bm{B_{n,1}}$};

          \node[] at (3.3,1.5) {$\bm{B_{1,j}}$};
          \node[] at (3.3,4.3) {$\bm{B_{k,j}}$};
          \node[] at (3.3,7.0) {$\bm{B_{n,j}}$};

          \node[] at (6.3,1.5) {$\bm{B_{1,m}}$};
          \node[] at (6.3,4.3) {$\bm{B_{k,m}}$};
          \node[] at (6.3,7.0) {$\bm{B_{n,m}}$};
        \end{tikzpicture}
    \end{minipage}
    \hspace{0.1\linewidth}
\begin{minipage}{0.4\linewidth}
        \centering
        \begin{tikzpicture}[scale=1]
            \node[opacity=0.6, anchor=south west, inner sep=0] (img) at (0,0) {\includegraphics[width=1\linewidth, angle=0]{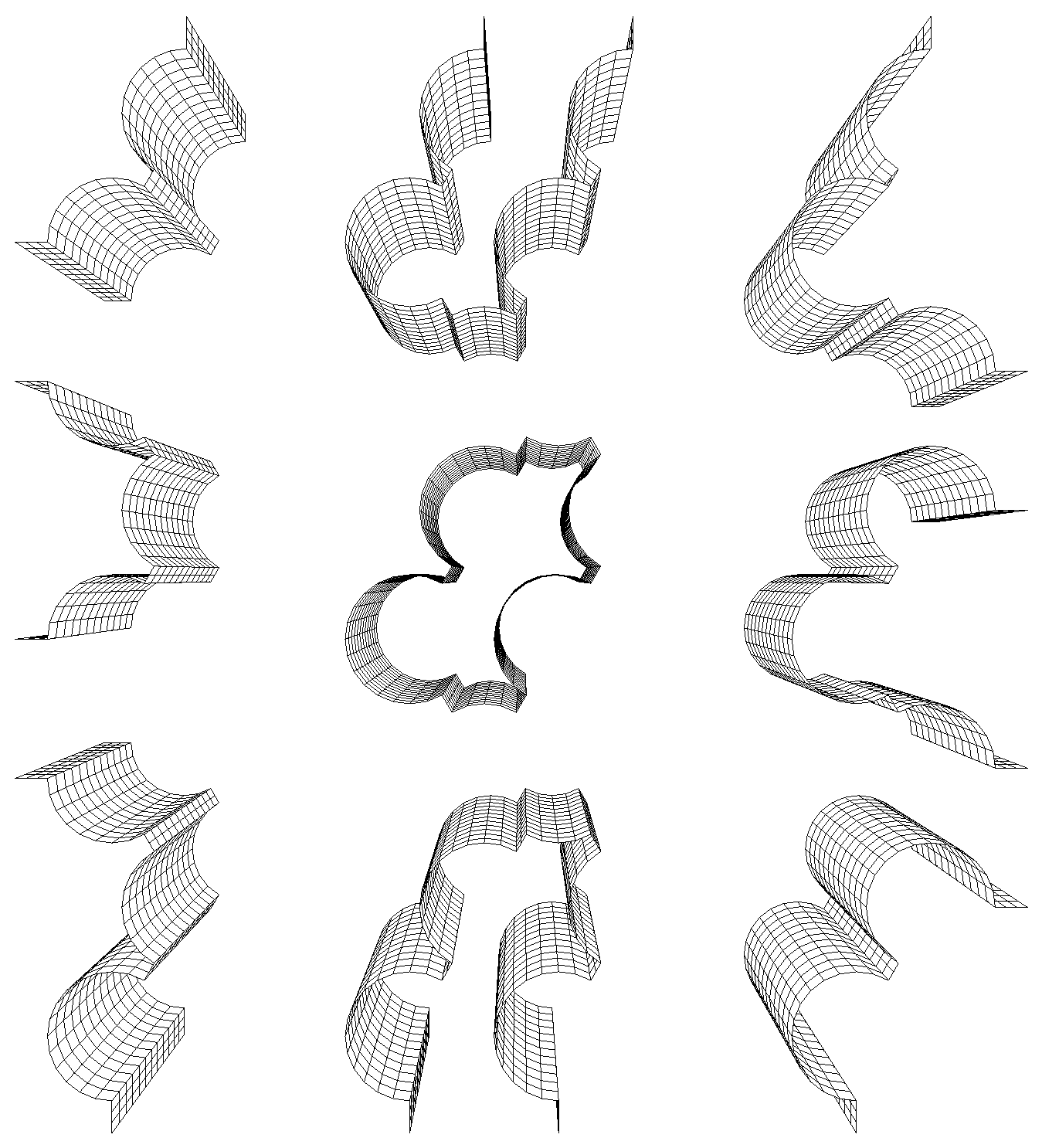}};
          \DotGrid{2}{0}{0.125}{0.75pt}{1.7}{6.4}
          \DotGrid{2}{0}{0.125}{0.75pt}{1.7}{4.1}
          \DotGrid{2}{0}{0.125}{0.75pt}{1.7}{1.7}
          \DotGrid{2}{0}{0.125}{0.75pt}{4.4}{6.4}
          \DotGrid{2}{0}{0.125}{0.75pt}{4.4}{4.1}
          \DotGrid{2}{0}{0.125}{0.75pt}{4.4}{1.7}
    
          \DotGrid{0}{2}{0.125}{0.75pt}{0.5}{5.25}
          \DotGrid{0}{2}{0.125}{0.75pt}{0.5}{2.95}
          \DotGrid{0}{2}{0.125}{0.75pt}{3.125}{4.85}
          \DotGrid{0}{2}{0.125}{0.75pt}{3.125}{2.55}
          \DotGrid{0}{2}{0.125}{0.75pt}{5.75}{4.85}
          \DotGrid{0}{2}{0.125}{0.75pt}{5.75}{2.55}

          \node[] at (0.5,1.5) {$\bm{\Gamma_{1,1}^{s}}$};
          \node[] at (0.5,4.3) {$\bm{\Gamma_{k,1}^{s}}$};
          \node[] at (0.5,7.0) {$\bm{\Gamma_{n,1}^{s}}$};
          
          \node[] at (3.3,1.5) {$\bm{\Gamma_{1,j}^{s}}$};
          \node[] at (3.3,4.3) {$\bm{\Gamma_{k,j}^{s}}$};
          \node[] at (3.3,7.0) {$\bm{\Gamma_{n,j}^{s}}$};
          
          \node[] at (6.3,1.5) {$\bm{\Gamma_{1,m}^{s}}$};
          \node[] at (6.3,4.3) {$\bm{\Gamma_{k,m}^{s}}$};
          \node[] at (6.3,7.0) {$\bm{\Gamma_{n,m}^{s}}$};
          
        \end{tikzpicture}
    \end{minipage}
    \caption{Schematic representation of the mesh partitioning into blocks $B_{n,m}$ (left) and the corresponding interfaces $\Gamma_{n,m}^{s}$, with $n=1,\ldots,N_{\text{block}}$ and $m=1,\ldots,M$, used for block merging (right). The blocks illustrate the content of the elementary files $f_s^{\text{el}}$ owned by the computational processes. Repetition of these blocks along the $x$- and $y$-directions allows the computational domain to be extended to an arbitrary size.}

    \label{fig:inner_domain_partition}
\end{figure}
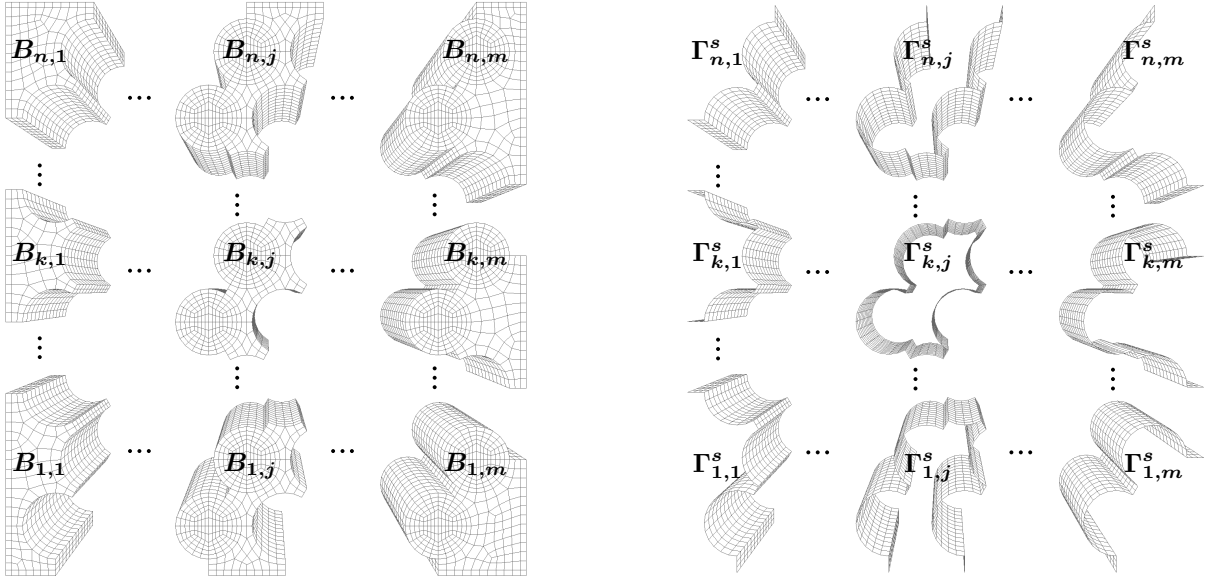

While element partitioning for $\Omega_{\rm out} \cup \Omega_{\rm PML}$ is performed by employing the Metis library,  
for the inner domain $\Omega_{\rm in}$, we consider a different strategy. 
Indeed, we need a sufficiently refined mesh for $\Omega_{\rm in}$ to accurately capture the geometry of the PMUT membranes and to properly resolve the applied voltage input signal on them. On the other hand we need to balance among processors: (i) the hexahedral elements that share a face with $\Gamma_\text{I}$, (ii) the number of PMUT membranes, and (iii) the remaining elements. 
%
%
%
Moreover, an additional issue to face in simulations involving a large computational domain is generating and extracting the mesh (to be partitioned) in a single step due to RAM limitations. 
To account for all these constraints, the main approach proposed in this study is to cut the mesh into smaller blocks and provide a tool to merge them into a single mesh directly within the code, since the of membranes involved is extremely large, reaching more than ten thousands PMUTs. Consequently, generating the complete mesh in a single step is often impractical due to memory limitations and time resources. 
\\\\The proposed strategy, which simultaneously considers the constraints imposed by PMUT arrays and the requirements of load balancing, is to generate a set of elementary files, $\{f^\text{el}_s\},$ with $s\in\{1,\dots,N_\text{sub}^\text{el}\}$, each containing the mesh of a subdomain of the inner region $\Omega_{\rm in}$. These subdomains are then assigned to individual computational processes, which collectively reconstruct the complete array structure. Because PMUT arrays are designed according to a predefined pattern, multiple processes can adopt the same elementary file, applying an appropriate coordinate shift to give rise to the global arrangement. Figure \ref{fig:inner_domain_partition} (left) provides an example of how those blocks appear. 
\\\\ 
In this way, four main constraints are directly satisfied: the PMUTs are distributed among different cores, keeping them entirely on the same process unit, the load balance of the inner domain is guaranteed, the DG faces on the finer side are well balanced, and the RAM limitation is overcame. Figure \ref{fig:minimal_inner_domain} illustrates a minimal example of the inner domain, along with the elementary structures used in its construction.

Next, we detailed the procedure developed to address the mesh generation and partitioning in $\Omega_{\rm in}$, when the number of PMUTs $N_p$ is sufficiently large. Let $N_{p,\text{row}}$ denote the number of PMUTs in a row, and let $N_\text{block} = N_{p,\text{row}} + 1$ represent the number of elementary blocks required to form a complete row in the inner domain.
The procedure begins by verifying that the number of computational cores $N_\text{proc}$ is an integer multiple of $N_\text{block}$, i.e. $N_\text{proc} = M N_\text{block}$ with $M \in \mathbb{N}^+$. Upon failure, the system raises an exception in order to avoid unfavorable workload imbalances, particularly within the finer regions of the mesh.
\\\\Subsequently, a connectivity structure $A_i$, for $i=0,...,N_\text{proc}-1$, is constructed for each process $i$ to define adjacency relationships. 
This ensures that each process $i$ has precise knowledge of its neighboring domains labeled in $A_i$, and thereby enabling efficient reconstruction of the global mesh connectivity. As an example, we can retrieve the mesh depicted in Figure \ref{fig:minimal_inner_domain} by setting $N_\text{block} = 3$, and $M = 3$ in Figure \ref{fig:inner_domain_partition}. The corresponding blocks $B_{n,m}$ and the core-ownership $L_{n,m}$ of each block, for $n=1,...,N_\text{block}$, $m=1,...,M$, can be grouped in the following matrices, respectively:
\[
B =
\begin{bmatrix}
B_{3,1} & B_{3,2} & B_{3,3} \\
B_{2,1} & B_{2,2} & B_{2,3} \\
B_{1,1} & B_{1,2} & B_{1,3}
\end{bmatrix}
,\quad
L =
\begin{bmatrix}
2 & 5 & 8 \\
1 & 4 & 7 \\
0 & 3 & 6
\end{bmatrix},
\]
from which the adjacency structure $A_i$, $\forall i =0,\dots,8$, is derived
\[
\begin{array}{lll}
A_0 = \begin{bmatrix} 1 & 3 & 4 \end{bmatrix}, &
A_1 = \begin{bmatrix} 0 & 2 & 3 & 4 & 5 \end{bmatrix}, &
A_2 = \begin{bmatrix} 1 & 4 & 5 \end{bmatrix}, \\[6pt]
A_3 = \begin{bmatrix} 0 & 1 & 4 & 6 & 7 \end{bmatrix}, &
A_4 = \begin{bmatrix} 0 & 1 & 2 & 3 & 5 & 6 & 7 & 8 \end{bmatrix}, &
A_5 = \begin{bmatrix} 1 & 2 & 4 & 7 & 8 \end{bmatrix}, \\[6pt]
A_6 = \begin{bmatrix} 3 & 4 & 7 \end{bmatrix}, &
A_7 = \begin{bmatrix} 3 & 4 & 5 & 6 & 8 \end{bmatrix}, &
A_8 = \begin{bmatrix} 4 & 5 & 7 \end{bmatrix}.
\end{array}
\]
So far, all the elementary files $f^\text{el}_i$ are independent of each other, due to the fact that nodes, faces, and hexahedra have a local numeration starting from 1. Then, we need an algorthm to merge all elementary blocks and define a global numbering for the mesh elements. 
We start by writing all nodes lying on shared interfaces $\Gamma^\text{s}_{n,m}$, for $n=1,...,N_\text{block}$, $m=1,...,M$, cf. Figure \ref{fig:inner_domain_partition} (right), in specific files $S_i$, $i=0,...,N_\text{proc}-1$. A quicksort algorithm orders the nodes according to the coordinates $x$, $y$, and $z$; in this way, the procedure for matching the nodes' IDs of different interfaces $\Gamma_{n,m}^\text{s}$ that represent the same physical surface is fast and affordable.
%
%
Notice that this step must be performed sequentially, and the connectivity structures $A_i$ provide a way to speed up the procedure. 
Rank zero keeps its original identifiers, while all other processes receive the identifiers of shared nodes from neighboring ranks with smaller IDs, ensuring that each common node is assigned a unique identifier and duplicates are removed. This procedure leads to non-continuous local numbering due to missing identifiers. A final global renumbering step is then performed across all mesh files to restore continuity and establish a consistent global numbering of nodes, faces, and elements.
The full procedure is synthetically shown in Algorithm \ref{algo:mesh-merging}.

\begin{algorithm}[H]
\caption{Inner domain partition and merging algorithm}\label{algo:mesh-merging}
\begin{algorithmic}[1]
\REQUIRE $N_\text{proc}$ multiple of $N_\text{block}$
\FOR{$i\in\{0,\dots,N_\text{proc}-1\}$}
\STATE Construct core-connectivity structure $A_i$
\STATE Assign to core $i$ its mesh file $f^\text{el}_i$
\STATE Renumber nodes, faces, and elements of $f^\text{el}_i$
\STATE Sort and store blocks' interface nodes in file $S_i$

\FOR{$j$ in $A_{i}$ s.t. $j<i$, $i\neq 0$}
    \STATE Match nodes in file $S_i$ and $S_j$
    \STATE Retain node IDs of file $S_j$ and discard IDs in $S_i$
\ENDFOR
\ENDFOR

\STATE Renumber nodes according to global IDs
\end{algorithmic}
\end{algorithm}

\section{Fast assembly of interface terms for non-matching grids}\label{sec:dg_optimization}
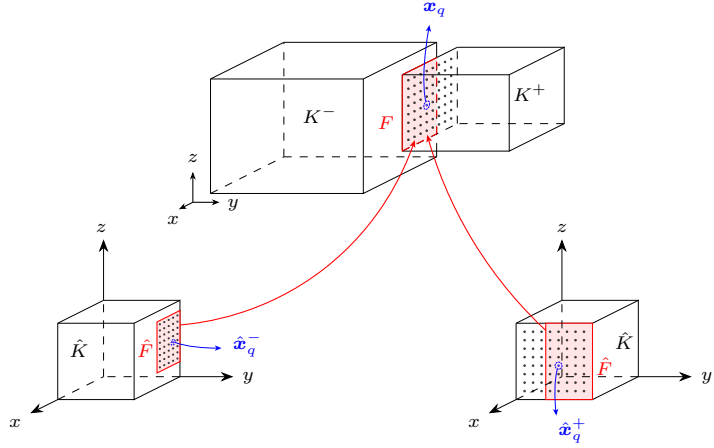
\begin{figure}[H]
    \centering
     \begin{tikzpicture}
        
        \node[inner sep=0pt] (img) at (0,0)
        {\includegraphics[angle=0, width=0.5\linewidth]{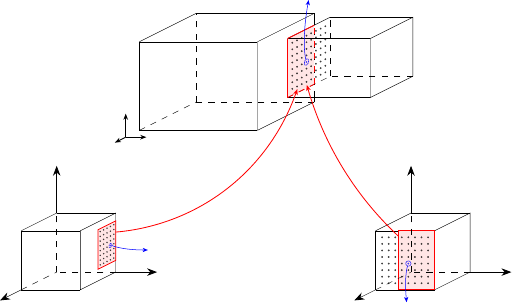}};

        \node[right] at (1.3,-2.7) {$\scriptstyle{x}$};
        \node[right] at (4.4,-2.1) {$\scriptstyle{y}$};
        \node[right] at (2.475,-0.1) {$\scriptstyle{z}$};
        \node[right] at (2.5,-2.8) {$\color{blue}\scriptstyle{\hat{\bm{x}}_q^+}$};
        \node[right] at (3.0,-1.94){$\color{red}\scriptstyle{\hat{F}}$};
        \node[right] at (3.25,-1.6){$\scriptstyle{\hat{K}}$};

        \node[right] at (-4.8,-2.7) {$\scriptstyle{x}$};
        \node[right] at (-1.7,-2.1) {$\scriptstyle{y}$};
        \node[right] at (-3.65,-0.1) {$\scriptstyle{z}$};
        \node[right] at (-1.825,-1.65){$\color{blue}\scriptstyle{\hat{\bm{x}}_q^-}$};
        \node[right] at (-3.1,-1.7){$\color{red}\scriptstyle{\hat{F}}$};
        \node[right] at (-4.0,-1.7){$\scriptstyle{\hat{K}}$};

        \node[right] at (-2.7,0) {$\scriptstyle{x}$};
        \node[right] at (-1.9,0.2) {$\scriptstyle{y}$};
        \node[right] at (-2.43,0.83) {$\scriptstyle{z}$};
        \node[right] at (0.7,2.8){$\color{blue}\scriptstyle{\bm{x}_q}$};
        \node[right] at (0.1,1.3){$\color{red}\scriptstyle{F}$};
        \node[right] at (-0.9,1.4){$\scriptstyle{K^-}$};
        \node[right] at (1.9,1.7){$\scriptstyle{K^+}$};
        
    \end{tikzpicture}
    \caption{The figure shows the physical hexahedra $K^-$ and $K^+$, sharing the common face $F$, together with their quadrature nodes $\bm{x}_q$, and the corresponding reference hexahedra $\hat{K}^\pm$ and their reference nodes $\hat{\bm{x}}_q^\pm$.}
    \label{fig:DG_hex_nodes_and_reference}
\end{figure}

The nonconformity between the inner and outer discretizations entails the construction of an interface mesh on $\Gamma_\text{I} = \partial \Omega_\text{in} \cap \partial \Omega_\text{out}$, defined as the geometric intersection of two hexahedral meshes with different resolutions. The number of resulting DG interfaces increases with the number of membranes under consideration. Efficiently identifying and matching interface elements becomes a critical task before the numerical solution of the governing equations can be carried out. While this strategy provides clear advantages in terms of geometric flexibility and reduction of the number of degrees of freedom, it also introduces significant challenges. In particular, a fast neighbor-matching algorithm is required to determine pairs of elements sharing DG interfaces, with the dual objective of reducing preprocessing overhead and improving performance during the time-stepping loop. In this work, we address these issues and propose an optimized methodology for interface construction and matching on $\Gamma_\text{I}$, significantly reducing the computational burden associated with the intersection process.
%
%
%
Note that the intersection faces $F \in \mathcal{F}^I$ between the coarse 
and the fine boundaries 
are defined starting from the quadrature nodes located on the surface with the finer resolution. These quadrature nodes therefore determine the connectivity between elements (without explicitly computing the geometric intersection of the faces), which is a crucial aspect for the evaluation of the integrals in \eqref{eq:Wf_2}, cf. \cite{mazzieri2012nonconforming}. In particular, we focus our attention to the terms in \eqref{eq:Wf_2}  of the form 
\begin{equation}\label{eq:integral_DG}
\int_{F}g^{-}(\bm x)f^+ (\bm x)\,\mathrm{d}\sigma, \end{equation}
being $F=\partial K^+ \cap \partial K^-$ and where superscript $\pm$ denotes the restriction of polynomial functions $g$ and $f$ to $K^+$ or $K^-$, respectively. 
To compute the integral in \eqref{eq:integral_DG}, we recall that every element $K$ of the mesh is the image through a trilinear map $\bm{X}_{K}$ of the reference element $\hat{K} = (-1,1)^3$. Therefore, using the change of variables $\bm x = \bm{X}_{K^\pm}(\hat{\bm x})$ it holds
\begin{align*}
\int_{F}g^{-}(\bm x)f^+ (\bm x)\,\mathrm{d}\sigma & =  \int_{\hat{F}}g^{-}(\bm{X}_{K^-}(\hat{\bm x}^-))f^+ (\bm{X}_{K^+}(\hat{\bm x}^+)) |J_{\bm X_{K^+}}(\hat{\bm x}^+)|\,\mathrm{d}\hat{\sigma},
\end{align*}
being $\hat{F}$ (a portion of) a face of $\partial \hat{K}$ mapped through $\bm X_{K^-}$ or $\bm X_{K^+}$ to $F$, see Figure \ref{fig:DG_hex_nodes_and_reference}.
The above integral is thus evaluated using
a suitable quadrature rule with $N_q$ nodes and $\{\omega_q\}_{q=1}^{N_q}$ weights as follows 

\begin{multline*}
\int_{\hat{F}}g^{-}(\bm{X}_{K^-}(\hat{\bm x}^-))f^+ (\bm{X}_{K^+}(\hat{\bm x}^+)) |J_{\bm X_{K^+}}(\hat{\bm x}^+)|\,\mathrm{d}\hat{\sigma} \\ \approx \sum_{q=1}^{N_q} g^{-}(\bm{X}_{K^-}(\hat{\bm x}^-_q))f^+ (\bm{X}_{K^+}(\hat{\bm x}^+_q)) |J_{\bm X_{K^+}}(\hat{\bm x}^+_q)| w_q = \sum_{q=1}^{N_q} \hat{g}^{-}(\hat{\bm x}^-_q)\hat{f}^+ (\hat{\bm x}^+_q) |J_{\bm X_{K^+}}(\hat{\bm x}^+_q)| w_q, 
\end{multline*}
where $w_q$, for $q=1,..,N_q$ are suitable weights and since $g^-(\bm x_q) = \hat{g}^{-}(\hat{\bm x}_q^-)$ and $f^+(\bm x_q) = \hat{f}^{+}(\hat{\bm x}_q^+)$.

In our process, we first define $\bm x_q \in F$ as the image of quadrature nodes $\hat{\bm x}_q \in \hat{F}$ through $\bm X_{K^+}$, i.e., $\bm x_q = X_{K^+}(\hat{\bm x}_q^+)$, 
and then, we map them back by using $\bm X_{K^-}$ 
to find $\hat{\bm x}_q^-\in \hat{F}$, i.e., $\hat{\bm x}_q^- =  {\bm X_{K^-}^{-1}}(\bm x_q)$, cf. Figure \ref{fig:DG_hex_nodes_and_reference}.
So, in practice, we 
implicitly define $F$ and $\hat{F}$ through quadrature points.

To accelerate the aforementioned computation, information regarding element faces $\partial K^+$ and $\partial K^-$ in $\mathcal{F}^I$ is stored in two separate files.
A third file contains the coordinates of a single vertex of each face $\partial K^-$ belonging to $\mathcal{T}_{\rm out}$. 
The algorithm for face-matching proceeds as follows: the finer face $\partial K^+$ loop over their quadrature nodes $\{\bm x_q\}_{q=1}^{N_q}$, and look for the coarser face containing the current node $\bm x_j$. \\\\In this scenario, where hundreds of thousands of faces appears, to narrow the number of faces among the preselected ones, we first consider all the coarse faces $\partial K^-$
whose normal $\bm{n}^-$ is opposite, up to a tolerance, to that of the current fine face $\bm{n}^+$. Hence, the search for the correct neighboring element is reduced to a smaller set of faces. Then, among the latter, we consider the ones whose vertex $V_{K^-}$ lies within a circular neighborhood $C_r(\bm x_j)$ of the quadrature node of the finer face, with radius $r=\max\{\text{diam}(\partial K^-),\text{diam}(\partial K^+)\} + \varepsilon$, $\varepsilon >0$. If those conditions are satisfied, the face is selected among the possible ones intersecting $\partial
K^+$.

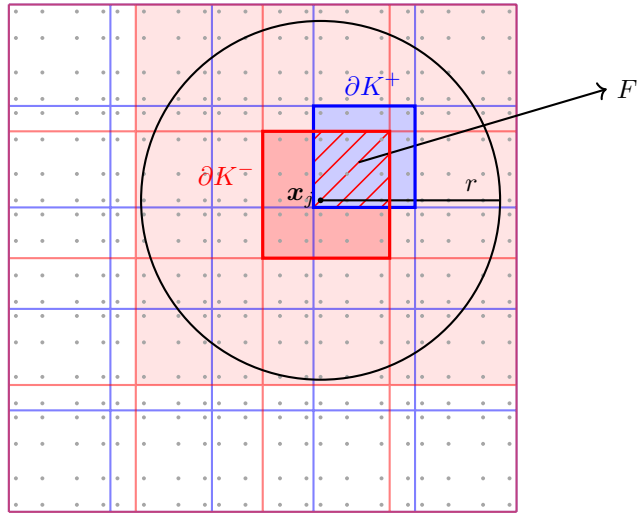
\begin{figure}[H]
    \centering
    \begin{tikzpicture}[scale=1.5]

        \begin{scope}
            \begin{scope}[on background layer]
                \fill[blue!70, opacity=1] (2.7,2.7) rectangle (3.375,3.375);
              \end{scope}
          \draw[step=0.9cm, blue, thick, opacity=0.5] (0,0) grid (4.5,4.5);
          \foreach \i in {0,...,4} {
            \foreach \j in {0,...,4} {
              \foreach \u in {0.06945,0.33005,0.66995,0.93055} {
                \foreach \v in {0.06945,0.33005,0.66995,0.93055} {
                  \fill[gray!70] (\i*0.9 + \u*0.9, \j*0.9 + \v*0.9) circle (0.5pt);
                }
              }
            }
          }
        \end{scope}
        
         \begin{scope}[on background layer]
            \fill[red!15, opacity=0.7] (1.125,1.125) rectangle (4.5,4.5);
            \fill[red!30, line width=1.2pt] (2.25,2.25) rectangle (3.375,3.375);
            \fill[blue!20, line width=1.2pt] (2.7,2.7) rectangle (3.6,3.6);

            \clip (2.7,2.7) rectangle (3.375,3.375);
        
            \foreach \i in {-2,-1.8,...,6} {
                \draw[red,  line width=0.6pt] (\i,0) -- ++(5,5);
            }
            
          \end{scope}
        \begin{scope}
          \draw[step=1.125cm, red, thick, opacity=0.5] (0,0) grid (4.5,4.5);
        \end{scope}

        \coordinate (C) at (2.762505,2.762505);
        \def\Radius{1.59099026}

        \coordinate (center) at (2.8,2.95);
        \node[below left] at (center) {$\bm{x}_j$};

        \draw[black, thick] (C) circle (\Radius);
        
        \fill[black] (C) circle (0.72pt);

        \draw[black, thick] (C) -- ++(\Radius,0);
        \node at (4.1,2.9) {$r$};

        \draw[->, thick, black]  (3.1,3.1) -- +(2.2,0.65) node[above, right] {$F$};

        \draw[blue, line width=1.2pt] (2.7,2.7) rectangle (3.6,3.6);
        \node at (3.23, 3.8){\color{blue}{$\partial K^+$}};

        \draw[red, line width=1.2pt] (2.25,2.25) rectangle (3.375,3.375);
        \node at (1.93, 3.0){\color{red}{$\partial K^-$}};

    \end{tikzpicture}
    \caption{Example of the matching algorithm of DG faces. The red and blue grids refer to the coarse and fine faces, respectively, and the gray dots are the LG quadrature nodes of the finer faces. The blue-filled face is a representative fine face, and the center $\bm{x}_j$ of the circle is its quadrature node considered; the red-filled faces are those selected by the algorithm. Notice that considering only one vertex of the coarse face instead of all four increases the number of residual faces, but it avoids excessive RAM consumption.}
    \label{fig:DG_matching}
\end{figure}

Finally, to properly find the face $\partial K^-$ sharing the node $\bm x_j$ with $\partial K^+$, a Newton-Raphson routine checks whether the quadrature point lies within the element $K^-$, computing its reference position in $[-1,1]^3$. This procedure is shown in Algorithm \ref{algo:DG-matching}, and Figure \ref{fig:DG_matching} provides a visual explanation.

\begin{algorithm}[H]
\caption{Face matching procedure}\label{algo:DG-matching}
\begin{algorithmic}[1]
\STATE Fix $\rm tol >0$
\FOR{$\partial K^+$  s.t.   $\partial K^+ \cap \Gamma_\text{I} \neq \emptyset$}
    \FOR{$\bm x_q$ with $q = 1, ..., N_q$}
        \FOR{$\partial K^-$  s.t.   $\partial K^- \cap \Gamma_\text{I} \neq \emptyset$} 
                \IF{$|\bm{n}_{+} +\bm{n}_{-}| < \rm tol$ \AND $V_{K^-} \in  C_{r}(\bm x_j)$}
                    \STATE Compute reference position $\hat{\bm x}_j$ with the Newton-Raphson method
                    \IF{$\hat{\bm x}_j \in [-1,1]^3$}
                    \STATE $K^-$ is a neighbour of $K^+$ 
                    \ENDIF
                \ENDIF
        \ENDFOR
    \ENDFOR
\ENDFOR

\end{algorithmic}
\end{algorithm}

\section{Numerical results}\label{sec:numerical_results}
This section presents the numerical results obtained from the PMUT array simulations. We first consider the case of a single PMUT configuration, which serves as a reference case for the proposed computational model. Once the model’s consistency is established, we investigate more complex scenarios involving multi-element PMUT arrays. The simulations are performed under two distinct operational conditions: a pure transmission (TX) framework, where the focus is on the emitted acoustic field, and a transmission–reception (TX–RX) framework, which additionally assesses the PMUTs sensing and reception capabilities.

\subsection{Single PMUT excitation in TX phase}\label{sec:TX_single_case}
\begin{figure}[H]
    \centering
    \begin{minipage}{0.3\linewidth}
        \centering
        \includegraphics[width=\linewidth]{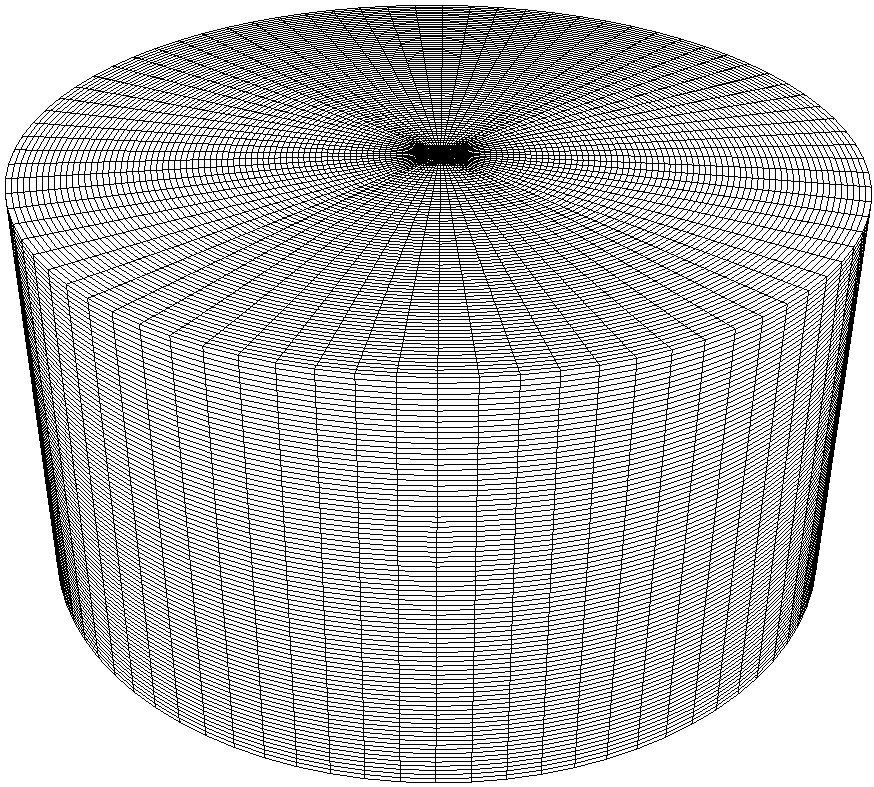}
        \caption*{(a) Lateral view.}
    \end{minipage}\hfill
    \begin{minipage}{0.3\linewidth}
        \centering
        \includegraphics[width=\linewidth]{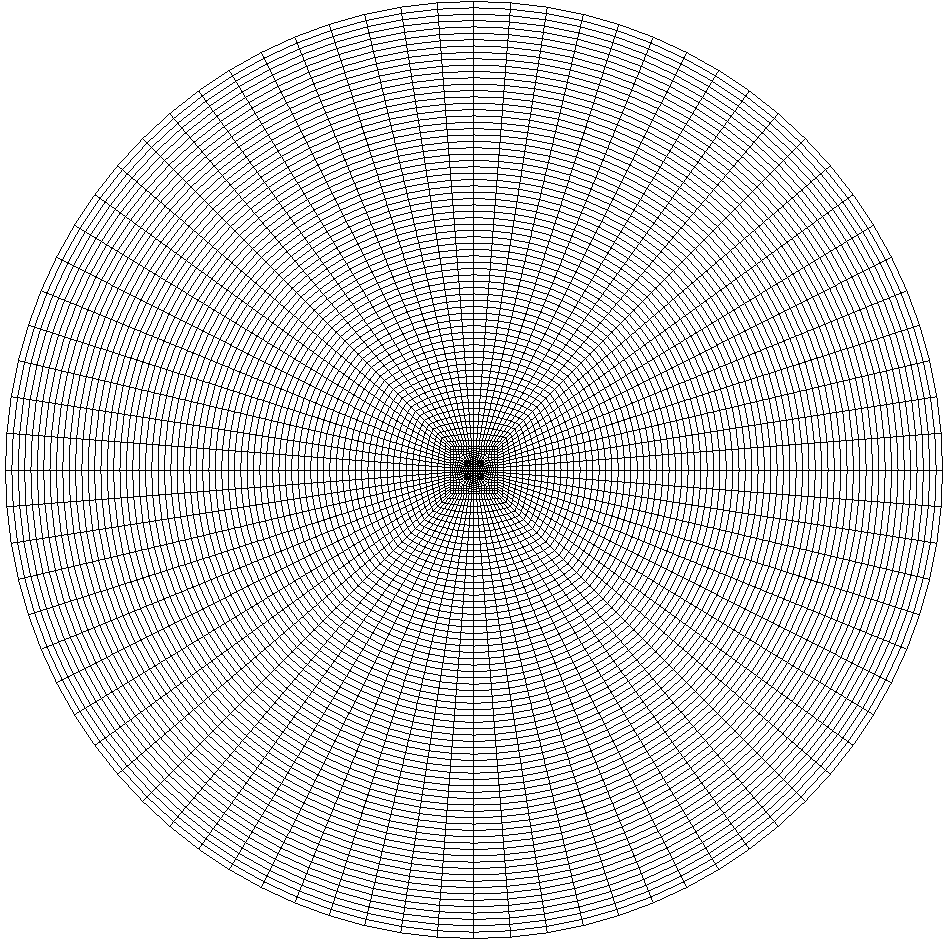}
        \caption*{(b) Top view.}
    \end{minipage}\hfill
\begin{minipage}{0.3\linewidth}
        \centering
        \begin{tikzpicture}
          \node[inner sep=0pt] (img) at (0,0)
            {\includegraphics[angle=90, width=1\linewidth]{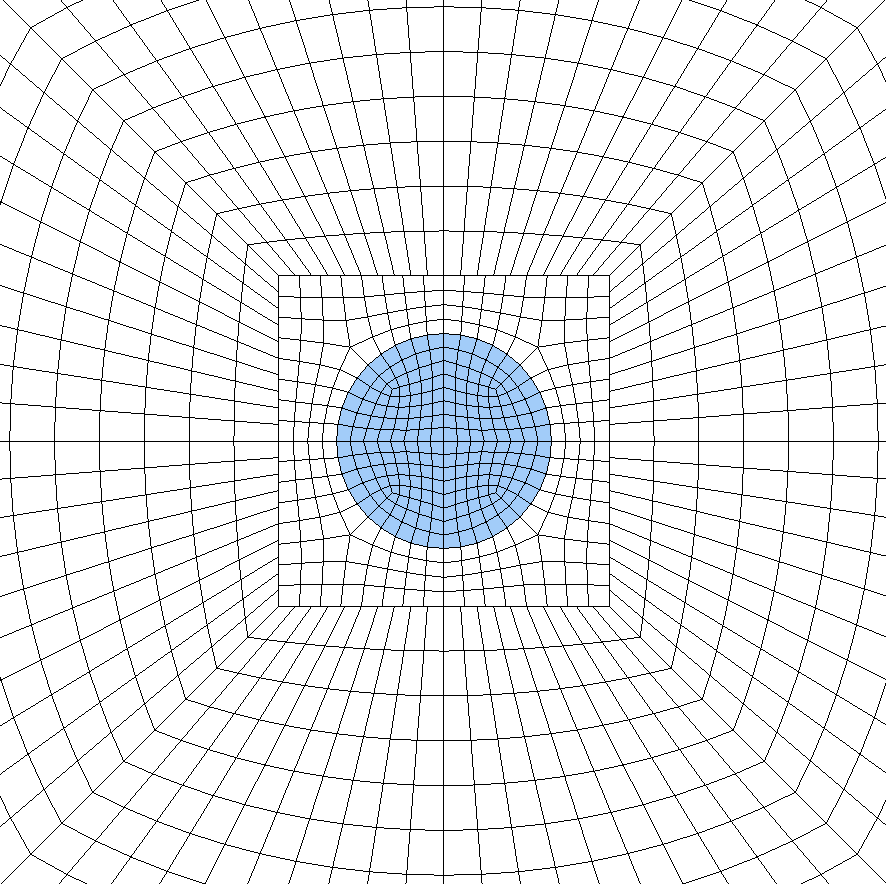}};

          \coordinate (node_z) at (-2.0, -2.0);
          \draw[->, thick, black]  (node_z) -- +(0.8, 0) node[right] {$x$};
          \draw[->, thick, black](node_z) -- +(0, 0.8) node[left] {$y$};
          \fill[black] (node_z) circle (1pt);
          \node[left] at (node_z) {$z$};
        \end{tikzpicture}
        \caption*{(c) Top view (zoom).}
    \end{minipage}
    \caption{Test case of Section \ref{sec:TX_single_case}. Mesh configuration used for the single PMUT simulation. Panels (a–c) show the lateral, top, and zoomed top views, respectively. The PMUT, highlighted in panel (c), is excited by an input voltage signal applied to its upper surface.}
    \label{fig:single_pmut_mesh}
\end{figure}

\begin{figure}[H]
    \begin{minipage}{0.3\linewidth}
        \centering
        \includegraphics[width=1.5\linewidth]{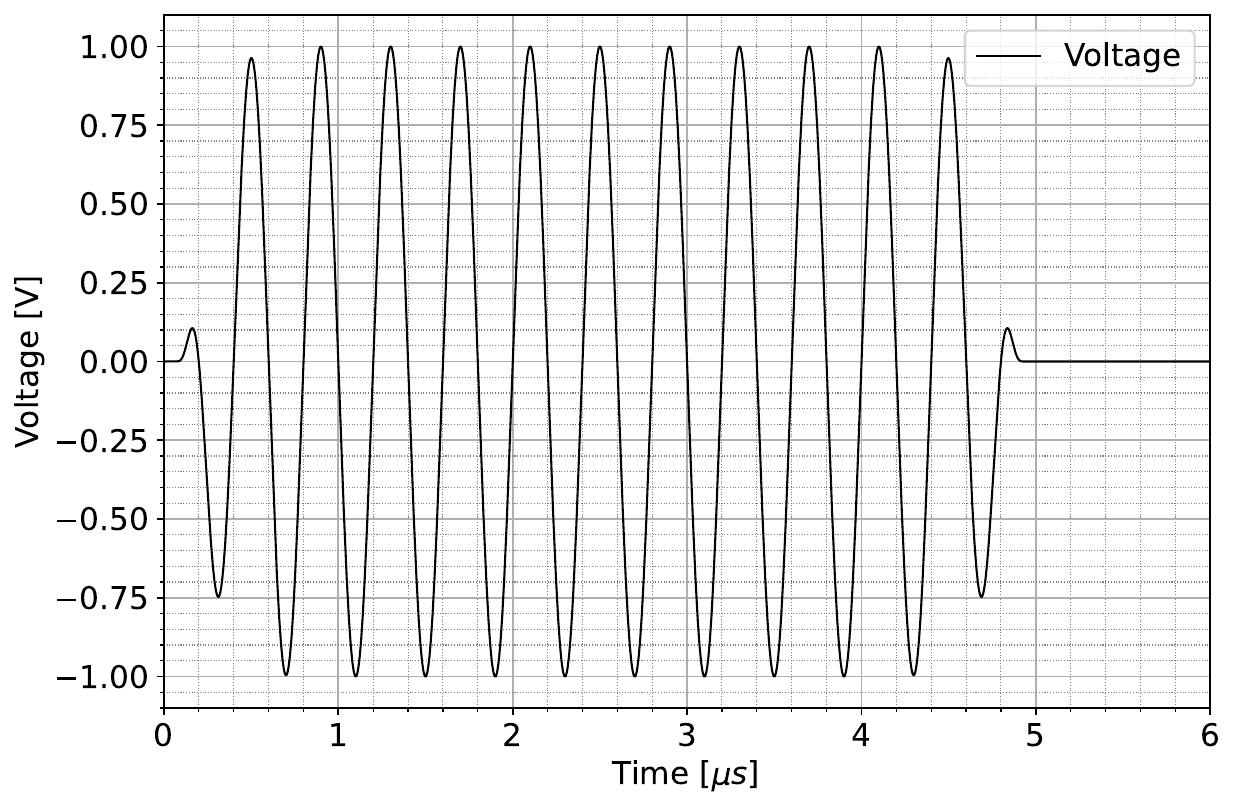}
    \end{minipage}\hspace{12em}
    \begin{minipage}{0.3\linewidth}
        \centering
        \includegraphics[width=1.5\linewidth]{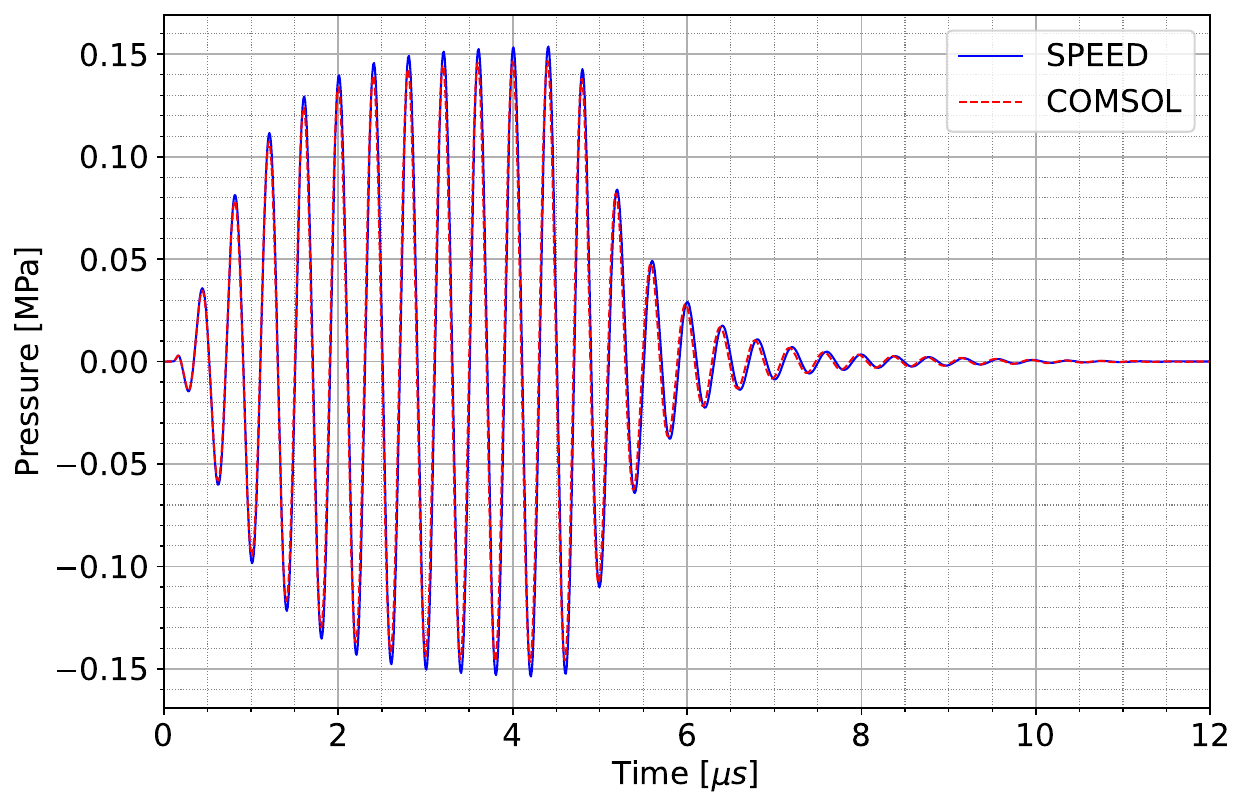}
    \end{minipage}
    \caption{Test case of Section \ref{sec:TX_single_case}. Time evolution of the input signal imposed on the PMUT surface (left), and output pressure computed at the PMUT center (right). SPEED and COMSOL signal are in blue and red, respectively. }
    \label{fig:voltage_plot_tx}
\end{figure}

The first analysis focuses on a single PMUT, i.e. $N_p=1$ in  \eqref{eq:pmut_problem}, operating in transmission mode. An alternating low-voltage signal of the form $\phi(t)=A\sin(2\pi ft)g(t)$ is applied to the surface of the PMUT, which is free to vibrate over the upper boundary of the computational domain. The modulation function $g(t)$ is introduced to ensure a smooth onset and offset of the excitation signal, thereby minimizing spurious frequency components and numerical artifacts associated with abrupt signal imposition. The function is defined as follows:
\begin{equation}\label{eq:voltage}
    g(t)=\exp\!\left(1 - \alpha(T_d-t)\right)\mathds{1}_{\{0 \le t < T_d/2\}}+ \exp\!\left(1 - \alpha(t)\right)\mathds{1}_{\{T_d/2 \le t \le T_d\}},
\end{equation}
where $\alpha(\tau) = \left(1 -\dfrac{2 \beta(\tau)^{1 + 2n}}{1 + \beta(\tau)^{2n}}\right)^{-2}$, with $\beta(\tau)=\left( 2\tau/T_d - 1 \right)$ and $n=10$. The excitation parameters are set to an amplitude of $A=\SI{1.0}{\volt}$, a driving frequency of $f=\SI{2.5}{\mega\hertz}$, and a total duration of $T_d=\SI{5}{\micro\second}$. The imposed voltage is plotted in Figure \ref{fig:voltage_plot_tx}. The computational domain comprises the PMUT surrounded by an acoustic medium for wave propagation, which extends for $\SI{2000}{\micro\meter}$  both in the radial and z directions. The radius of the PMUT is $\SI{65}{\micro\meter}$. To avoid spurious reflections, first-order absorbing boundary conditions are imposed along the lateral and bottom boundaries of the acoustic domain, while a hard-wall constraint is imposed at the top surface supporting the PMUT as in (\ref{eq:full_problem}), cf. Figure \ref{fig:single_pmut_mesh}. A time step $\Delta t = \SI{2e-4}{\micro\second}$ and a second-order polynomial degree in space are selected. The final time is set to $T=\SI{12}{\micro\second}$. Water is employed for both the inner and outer domains, hence the material density and the velocity of acoustic waves are $\rho = \SI{1e-3}{\kilo\gram\per\meter\cubed}$ and $c = \SI{1481}{\meter\per\second}$, respectively. The first three axially symmetric modes, i.e., $N_m=3$  in  (\ref{eq:pmut_problem}), are considered.

\begin{figure}[H]
    \begin{minipage}{0.3\linewidth}
        \centering
        \includegraphics[width=1.3\linewidth]{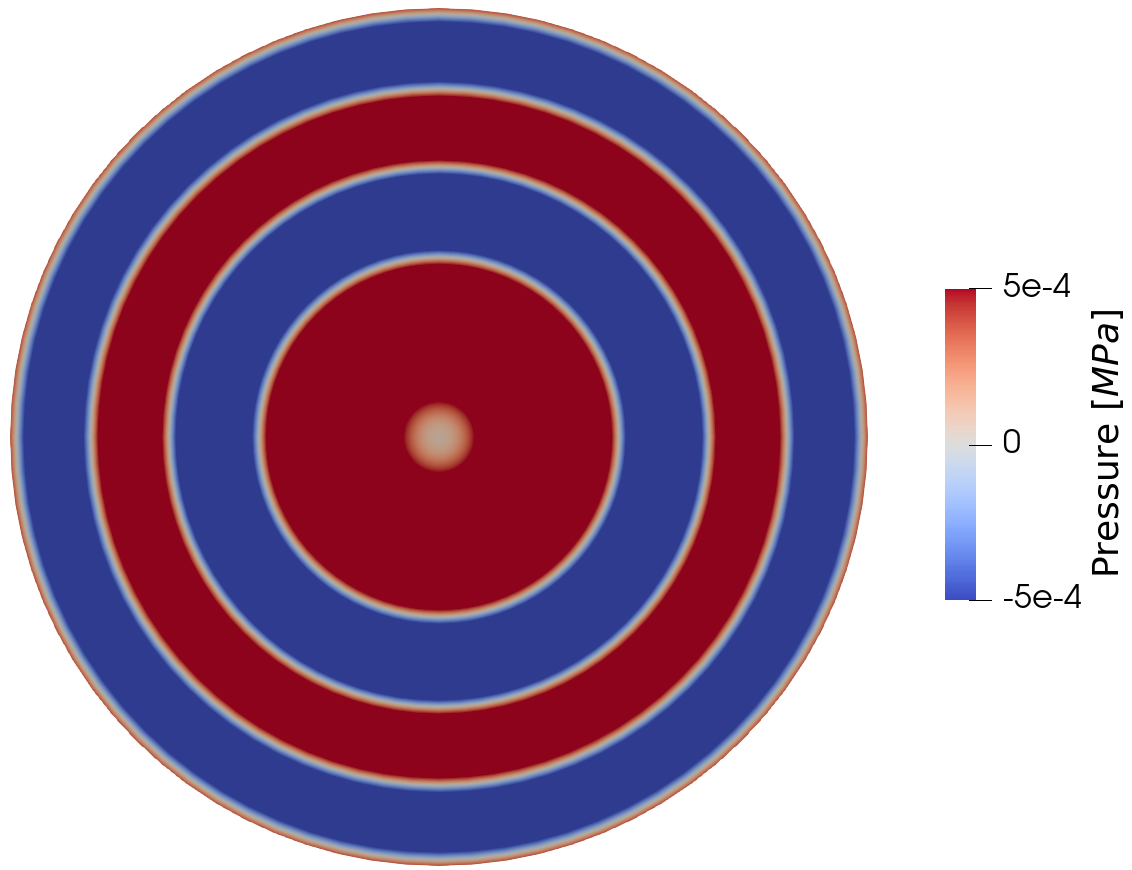}
    \end{minipage}\hspace{12em}
    \begin{minipage}{0.3\linewidth}
        \centering
        \includegraphics[width=1.5\linewidth]{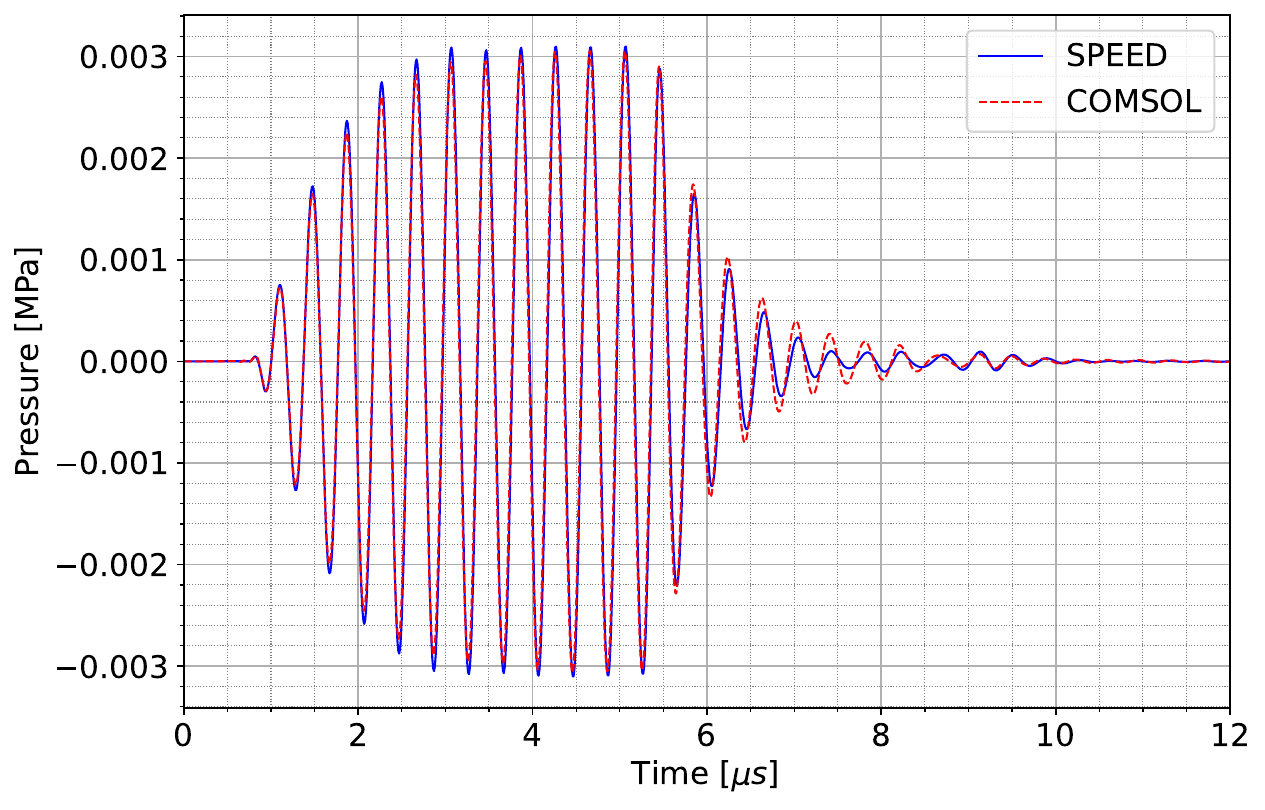}
    \end{minipage}
    \caption{Test case of Section \ref{sec:TX_single_case}. Snapshot at $t=\SI{4}{\micro\meter}$ of a slice in the $x$-$y$ plane (left), and output pressure computed at the $P=(\SI{0}{\micro\meter}, \SI{0}{\micro\meter}, \SI{-1000}{\micro\meter})$ for SPEED and COMSOL (right).}
    \label{fig:output_pressure}
\end{figure}

Figure \ref{fig:voltage_plot_tx} illustrates the temporal profile of the smooth input voltage applied to the PMUT surface, together with the resulting pressure excitation at the center of the membrane. The results are compared with a COMSOL Multiphysics V6.3 simulation of the PMUT array, where the device is modeled using a hybrid solid-shell approach: thin passive layers are represented using 2D first-order shear deformation theory shell elements, while the PZT layer and the silicon substrate are retained as 3D solid elements. Displacement continuity is enforced across shell--solid interfaces. The front encapsulation is modeled as an acoustic fluid with impedance-matched properties, and the surrounding acoustic domain solves the linearized acoustic wave equation  \cite{Abdalla2022PMUT}. This configuration provides an accurate approximation of a full 3D model while reducing computational cost, and the comparison demonstrates close correspondence between the two simulation frameworks. Figure \ref{fig:output_pressure}, instead, provides a visualization of the acoustic wavefield propagation on an $x$--$y$ plane extracted at a depth of $z=\SI{-1000}{\micro\meter}$, together with the temporal evolution of the pressure at $P=(\SI{0}{\micro\meter}, \SI{0}{\micro\meter}, \SI{-1000}{\micro\meter})$. The comparison with COMSOL further confirms the consistency and agreement between the two simulations.

\subsection{Single PMUT with an obstacle in TX-RX phase}\label{sec:single_tx-rx}
In general, PMUTs are designed to be used in both the transmission (TX) and the receiving (RX) phases, i.e., to trigger the propagation of the acoustic signal up to a specific target and to capture the response that they receive back from the target itself. To this aim, the following simulation deals with both TX and RX phases.

\begin{figure}[H]
    \centering
    \begin{minipage}{0.3\linewidth}
        \centering
        \begin{tikzpicture}
          \node[inner sep=0pt] (img) at (0,0)
            {\includegraphics[angle=0, width=0.7\linewidth]{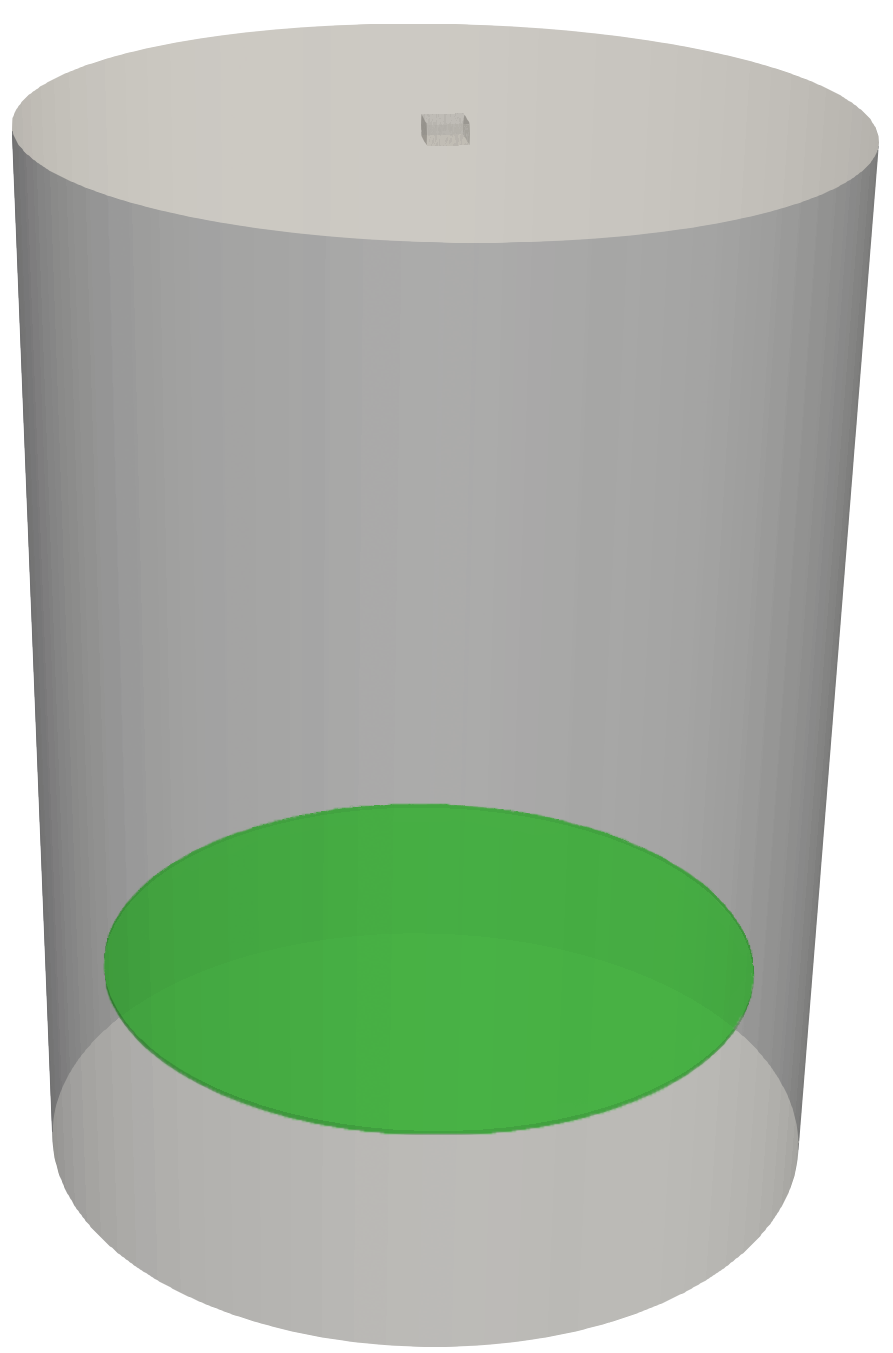}};

          \coordinate (node_z) at (-1.2, 0.5);
          \draw[->, thick, black]  (node_z) -- +(0.8, 0) node[right] {$x$};
          \draw[->, thick, black](node_z) -- +(0, 0.8) node[left] {$z$};
          \draw[->, thick, black]  (node_z) -- +(0.3, 0.3) node[right] {$y$};
        \end{tikzpicture}
        \caption*{(a) Full geometry.}
    \end{minipage}\hfill
    \begin{minipage}{0.3\linewidth}
        \centering
        \begin{tikzpicture}
          \node[inner sep=0pt] (img) at (0,0)
            {\includegraphics[angle=0, width=0.7\linewidth]{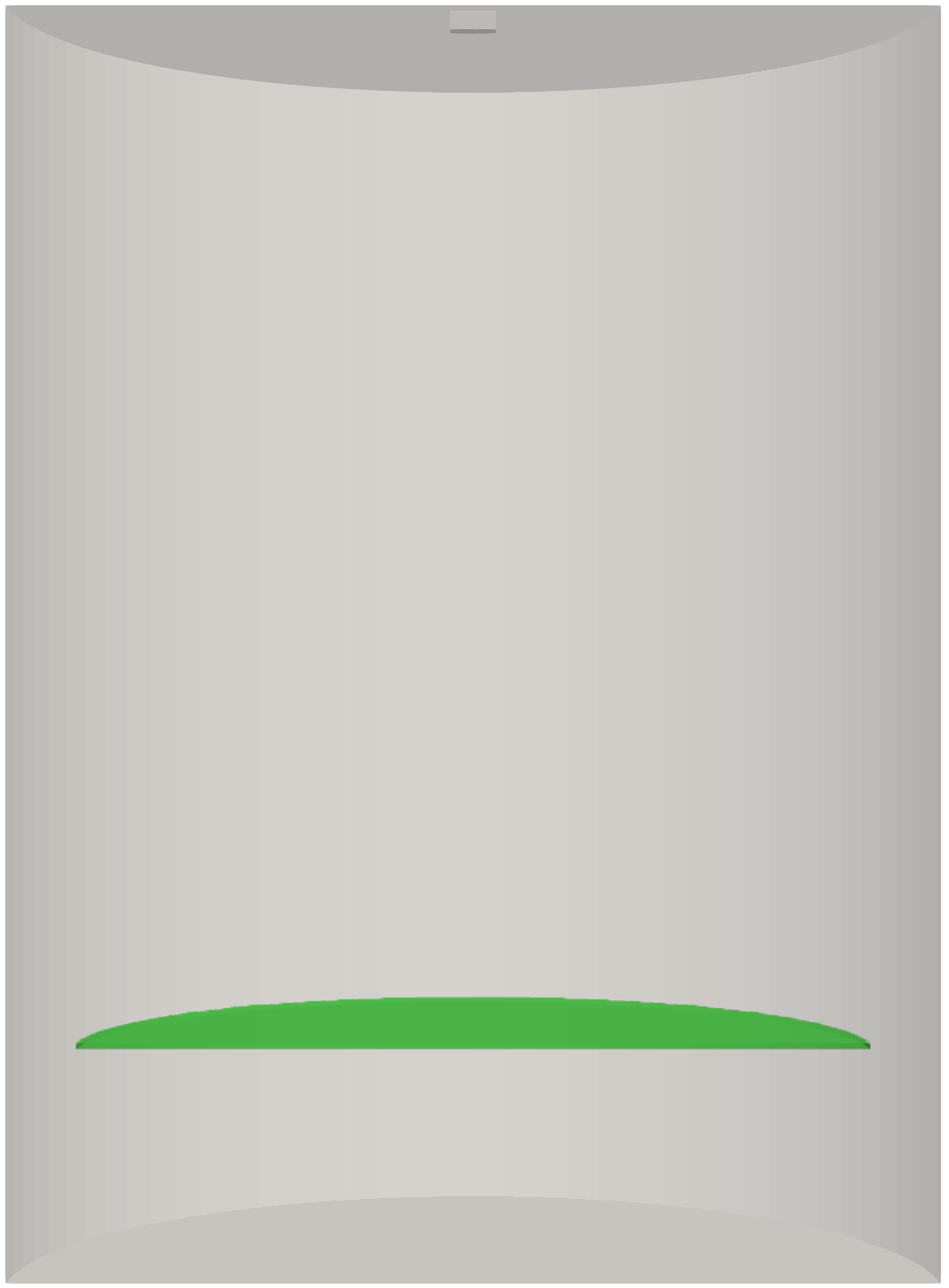}};

          \coordinate (node_h1) at (0.0, -1.56);
          \coordinate (node_h3) at (0.0, -2.485);
          \draw[-, thick, black](node_h1) -- +(0, 4.04) node[midway,left] {$h_1$};
          \draw[-, thick, black](node_h3) -- +(0, 0.9) node[midway,left] {$h_3$};
        \end{tikzpicture}
        \caption*{(b) Thickness of the layers.}
    \end{minipage}\hfill
    \begin{minipage}{0.3\linewidth}
        \centering
        \begin{tikzpicture}
          \node[inner sep=0pt] (img) at (0,0)
            {\includegraphics[angle=0, width=1.\linewidth]{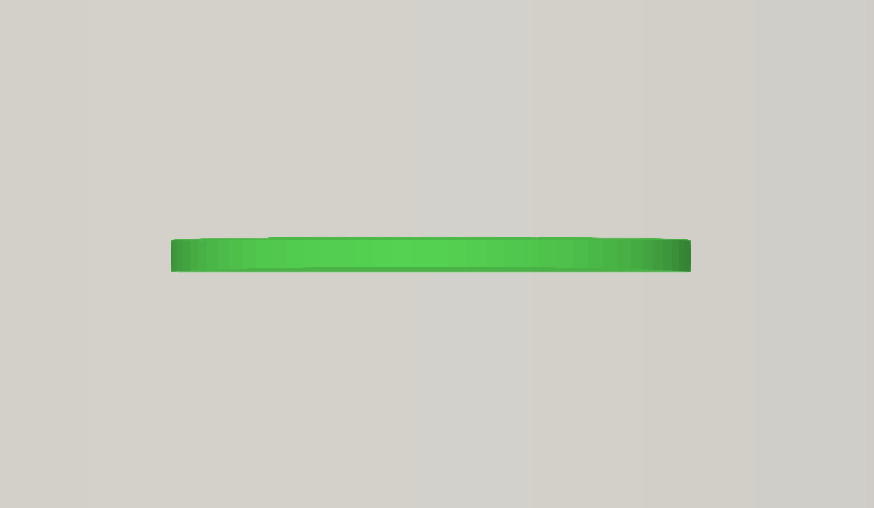}};

          \coordinate (node_h3) at (-1.67, -0.115);
          \coordinate (node_r) at (0.0, 0.15);
          \draw[-, thick, black](node_h3) -- +(0, 0.205) node[midway,left] {$h_2$};
          \draw[-, thick, black](node_r) -- +(1.52,0) node[midway, above] {$r_\text{obs}$};
        \end{tikzpicture}
        \caption*{(c) Obstacle zoom.}
    \end{minipage}

    \caption{Test case of Section \ref{sec:single_tx-rx}. Geometry used for the single PMUT simulation with an obstacle in TX-RX phases. Panels (a–c) show the geometry employed, the thickness of the cylinder layers, and a zoom of a representative obstacle for graphical reasons, respectively.}
    \label{fig:obstacle}
\end{figure}

\begin{figure}[H]
    \begin{minipage}{0.3\linewidth}
        \centering
        \includegraphics[width=1.5\linewidth]{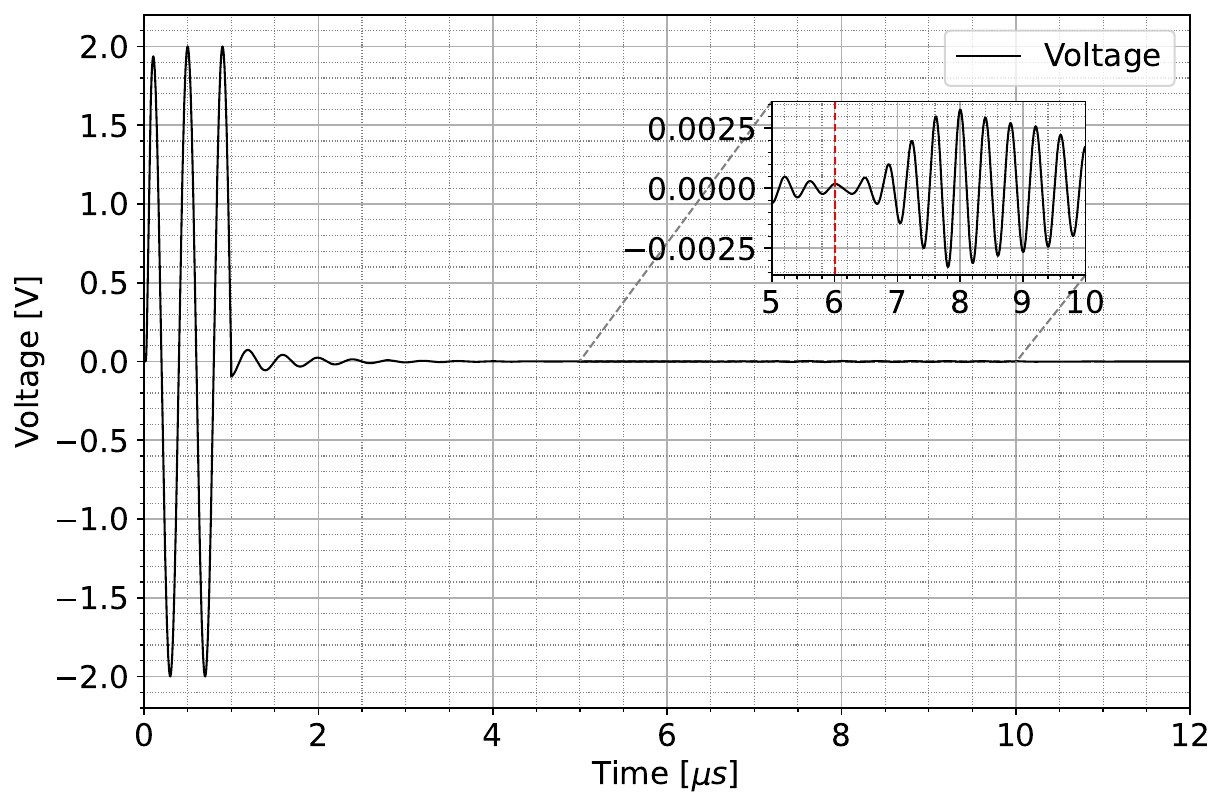}
    \end{minipage}\hspace{12em}
    \begin{minipage}{0.3\linewidth}
        \centering
        \includegraphics[width=1.5\linewidth]{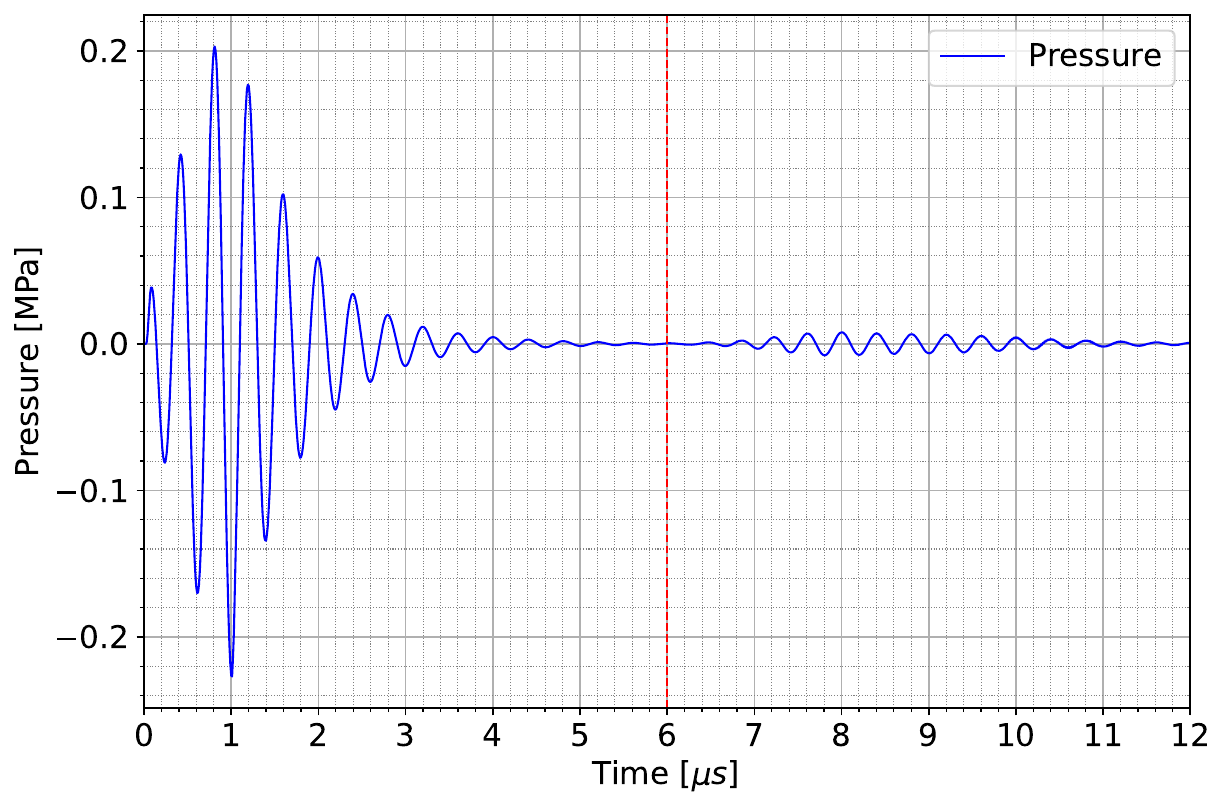}
    \end{minipage}
    \caption{Test case of Section \ref{sec:single_tx-rx}. Time evolution of the input signal imposed on the PMUT surface (left) with a zoom on the region for $t\in[5,10]\SI{}{\micro\second}$, and output pressure computed at the PMUT center (right). The time instant $t=6T_d$ is highlighted with a dashed line.}
    \label{fig:voltage_plot_tx_rx}
\end{figure}

An obstacle representing the target is put inside the acoustic propagation domain as shown in Figure \ref{fig:obstacle}. The setup is the following: an AC low-voltage $\phi(t)=A\sin(2\pi ft)h(t)$ is used with
\begin{equation*}
    h(t)=\exp\!\left(1 - \alpha(T_d-t)\right)\mathds{1}_{\{0 \le t \le T_d\}},
\end{equation*}
where $\alpha(\tau)$ is defined as in the test case of Section \ref{sec:TX_single_case}, with voltage parameters $A=\SI{2.0}{\volt}$,  $f=\SI{2.5}{\mega\hertz}$, and a total duration $T_d=\SI{1}{\micro\second}$. The acoustic water domain is a cylinder of radius $r=\SI{2000}{\micro\meter}$ composed of three layers of height $h_1 = 3c$, $h_2 = \SI{25}{\micro\meter}$, and $h_3 = \SI{1000}{\micro\meter}$, respectively, where $h_1$ is exactly the distance traveled by the acoustic wave in a time equivalent to $3T_d$ and $c = \SI{1481}{\meter\per\second}$ is the wave speed in water. In this way, the acoustic pressure wave fully reaches the obstacle with radius $r_\text{obs}=\SI{1700}{\micro\meter}$ placed at a distance $h_1$ from the PMUT center, and the pressure is received back from the membrane after a time equivalent to $6T_d$ (see Figure \ref{fig:voltage_plot_tx_rx}). From time $T_d$ up to the final time $T=\SI{12}{\micro\second}$, the PMUT is left in open-circuit condition with no externally applied voltage. The resulting mechanical vibrations and reflected pressure waves induce a voltage across the electrodes via the direct piezoelectric effect. Absorbing boundary conditions are imposed on the lateral and bottom sides of the outer cylinder, while a homogeneous Neumann condition is enforced on the top side and on the obstacle's surface.

\begin{figure}[H]
    \begin{minipage}{0.3\linewidth}
        \centering
        \includegraphics[width=1.5\linewidth]{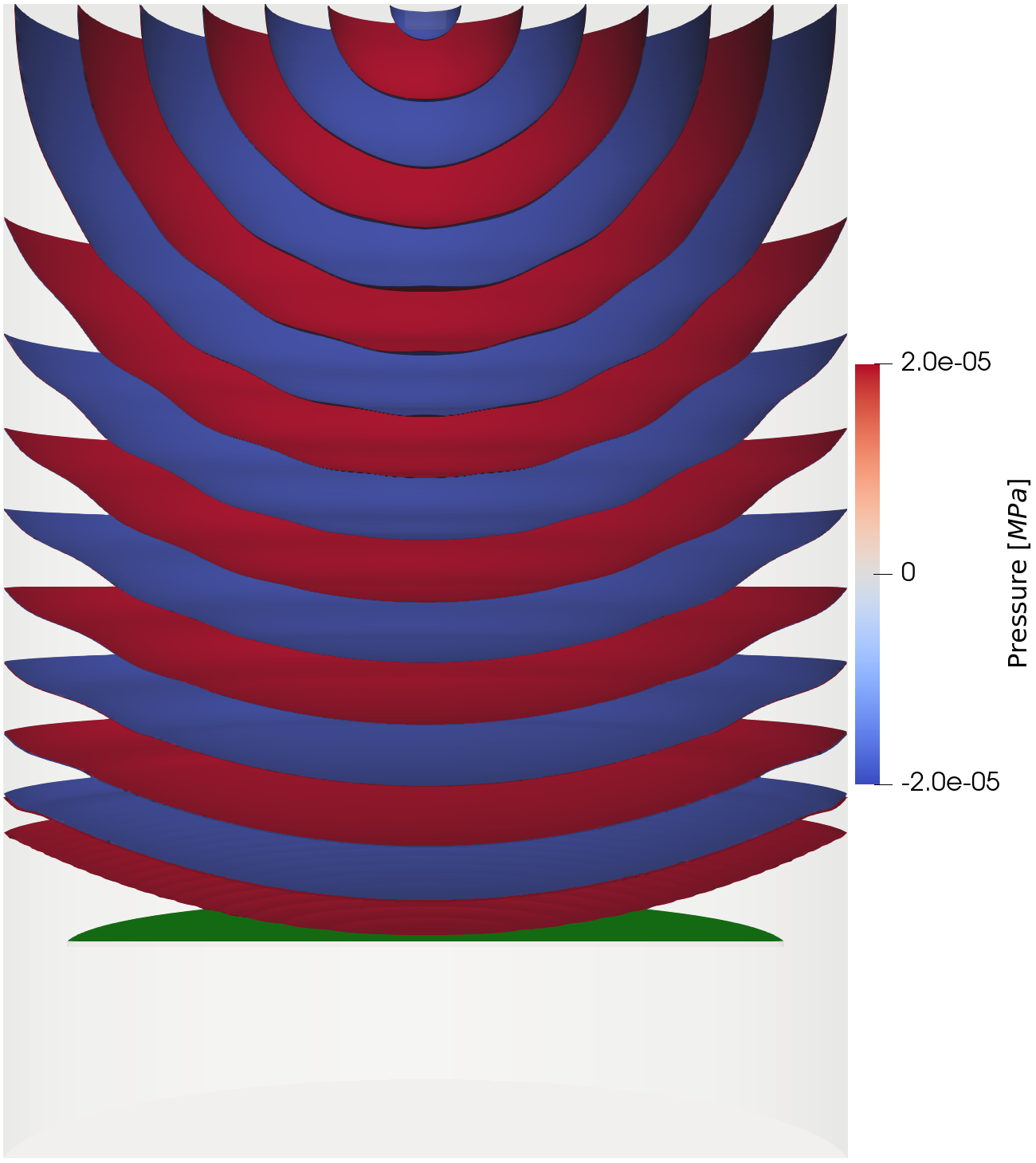}
    \end{minipage}\hspace{12em}
    \begin{minipage}{0.3\linewidth}
        \centering
        \includegraphics[width=1.5\linewidth]{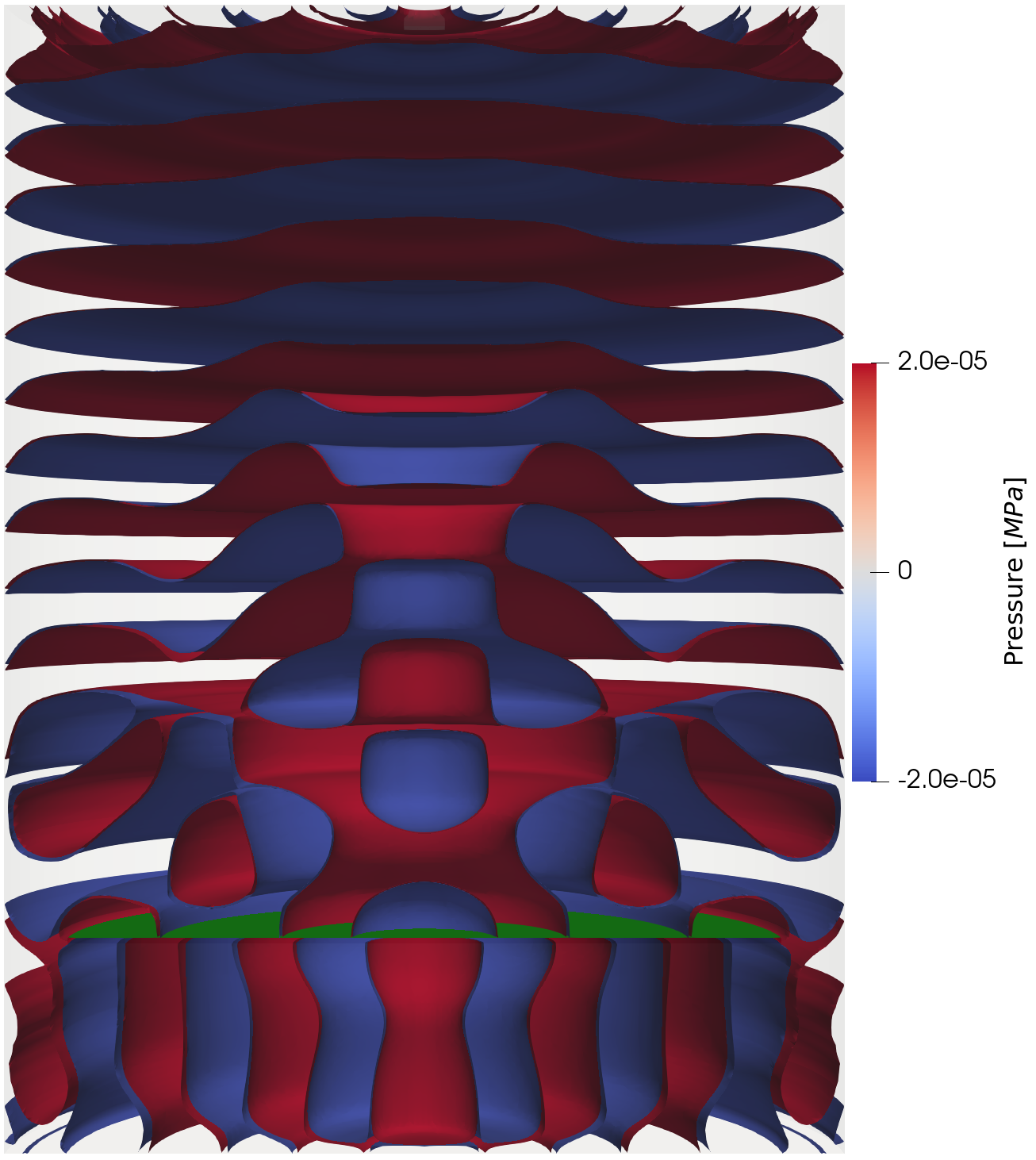}
    \end{minipage}
    \caption{Test case of Section \ref{sec:single_tx-rx}. Snapshots of the acoustic pressure waves contour plot: on the left, the one hitting the obstacle (in green), on the right, the one that travels back up to the PMUT.}
    \label{fig:contour_plot}
\end{figure}

\noindent Figure \ref{fig:voltage_plot_tx_rx} shows the time evolution of both the voltage and the corresponding pressure computed at the PMUT center. The behavior of the PMUT changes at time $T_d$: the voltage is no more imposed, hence the PMUT ends the transmission phase and enters in receiving phase. From time $T_d$ up to $t=6T_d$, the voltage amplitude decays since the pressure on the PMUT surface reduces, while after a time $t>6T_d$ the voltage slightly increases due to the pressure reflected from the obstacle. This trend is also visible in the pressure plot, where an increase in value is shown after the characteristic time $t=6T_d$, cf. Figure \ref{fig:voltage_plot_tx_rx} (right). Figure $\ref{fig:contour_plot}$ presents two snapshots of the iso-contours of the acoustic pressure; the one on the left is a snapshot at time $t=3T_d$, where the waves move towards the obstacle represented in green, while the picture on the right highlights the pressure waves that reach the PMUT surface at time $6T_d$ after the reflection.

\subsection{Multi-array PMUT simulation}\label{sec:multi_array}
\begin{figure}[H]
    \centering
    \begin{minipage}{0.6\linewidth}
        \centering
        \hspace{-1.93em}
        \begin{tikzpicture}
            \node[inner sep=0pt] (img) at (0,0)
                {\includegraphics[width=0.935\linewidth]{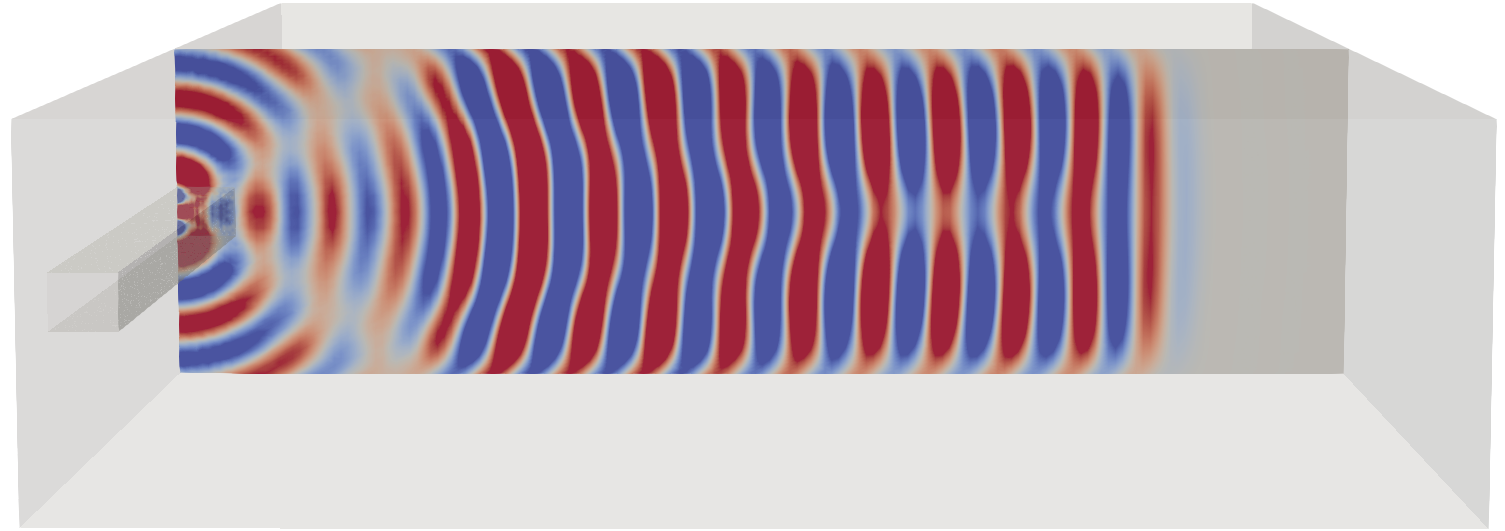}};
            \begin{scope}[]
                \node[below] at (3.5,-0.75) {$\Omega_\text{in}\cup\Omega_\text{out}$};
            \end{scope}
        \end{tikzpicture}
    \vspace{0.5em}
    \end{minipage}
    \begin{minipage}{0.6\linewidth}
        \centering
        \begin{tikzpicture}
            \node[inner sep=0pt] (img) at (0,0)
                {\includegraphics[width=1.002\linewidth]{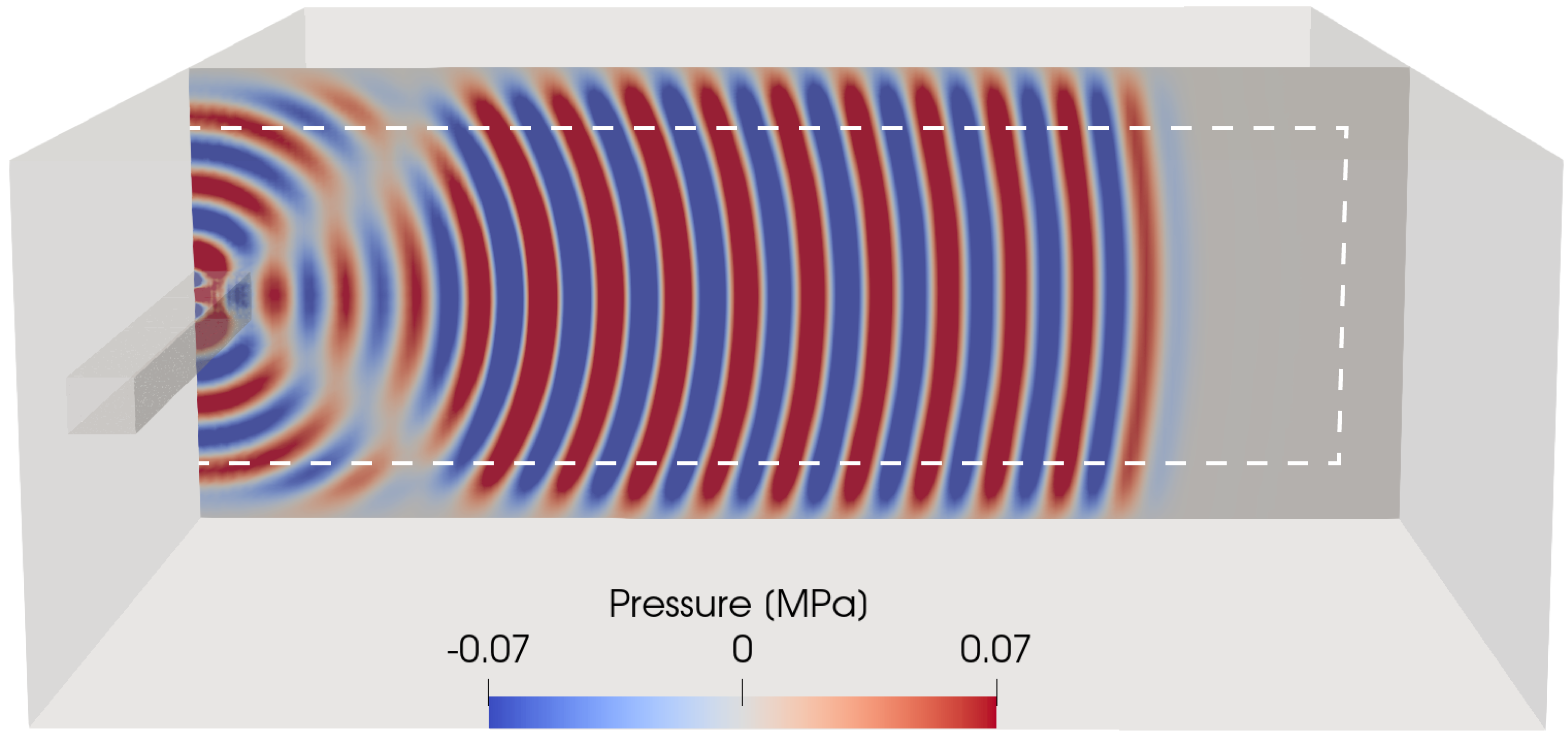}};
            \begin{scope}[]
                \node[below] at (3.5,-1.2) {$\Omega_\text{in}\cup\Omega_\text{out}\cup\Omega_\text{PML}$};
            \end{scope}
            \coordinate (node_z) at (6.0, -1.2);
            \draw[->, thick, black] (node_z) -- +(0.8, 0) node[right] {$z$};
            \draw[->, thick, black] (node_z) -- +(0, 0.8) node[left] {$x$};
            \fill[black] (node_z) circle (1pt);
            \node[left] at (node_z) {$y$};
        \end{tikzpicture}
    \end{minipage}
    \caption{Test case of Section \ref{sec:multi_array}. Schematic representation of a slice of the acoustic domain $\Omega_{\text{in}} \cup \Omega_{\text{out}}$ and the surrounding perfectly matched layer region $\Omega_{\text{PML}}$. The upper figure corresponds to the simulation without the PML, while the lower figure shows the solution obtained when the PML is included. The dashed line represents the interface $\Gamma_\text{PML}$.
}
    \label{fig:pml_examples}
\end{figure}

In general, absorbing boundary conditions are well-suited for simulations where wave propagation is predominantly orthogonal to the boundary. This approach is particularly effective when the computational domain is spherical or cylindrical, namely when the shape of the boundary matches the shape of the propagating waves. However, classical absorbing conditions can lead to undesirable internal reflections in domains with corners. 
Figure \ref{fig:pml_examples} illustrates the influence of the PML on a single-element array within a computational domain featuring corners. The two images correspond to the same simulation performed with and without $\Omega_{\rm PML}$. In the absence of this layer, wave reflections arising at the domain boundaries lead to noticeable distortions in the wavefront, thereby compromising the accuracy of the numerical solution. Conversely, when the PML is applied in combination with standard absorbing boundary conditions, spurious reflections are effectively suppressed, allowing for a smoother and more physically consistent propagation of acoustic waves within the domain.

\begin{figure}[H]
\centering












\begin{minipage}[t]{0.7\linewidth}
\raggedleft
\begin{tikzpicture}

\node[anchor=south west, inner sep=0] (img) at (0,0)
{\includegraphics[width=\linewidth]{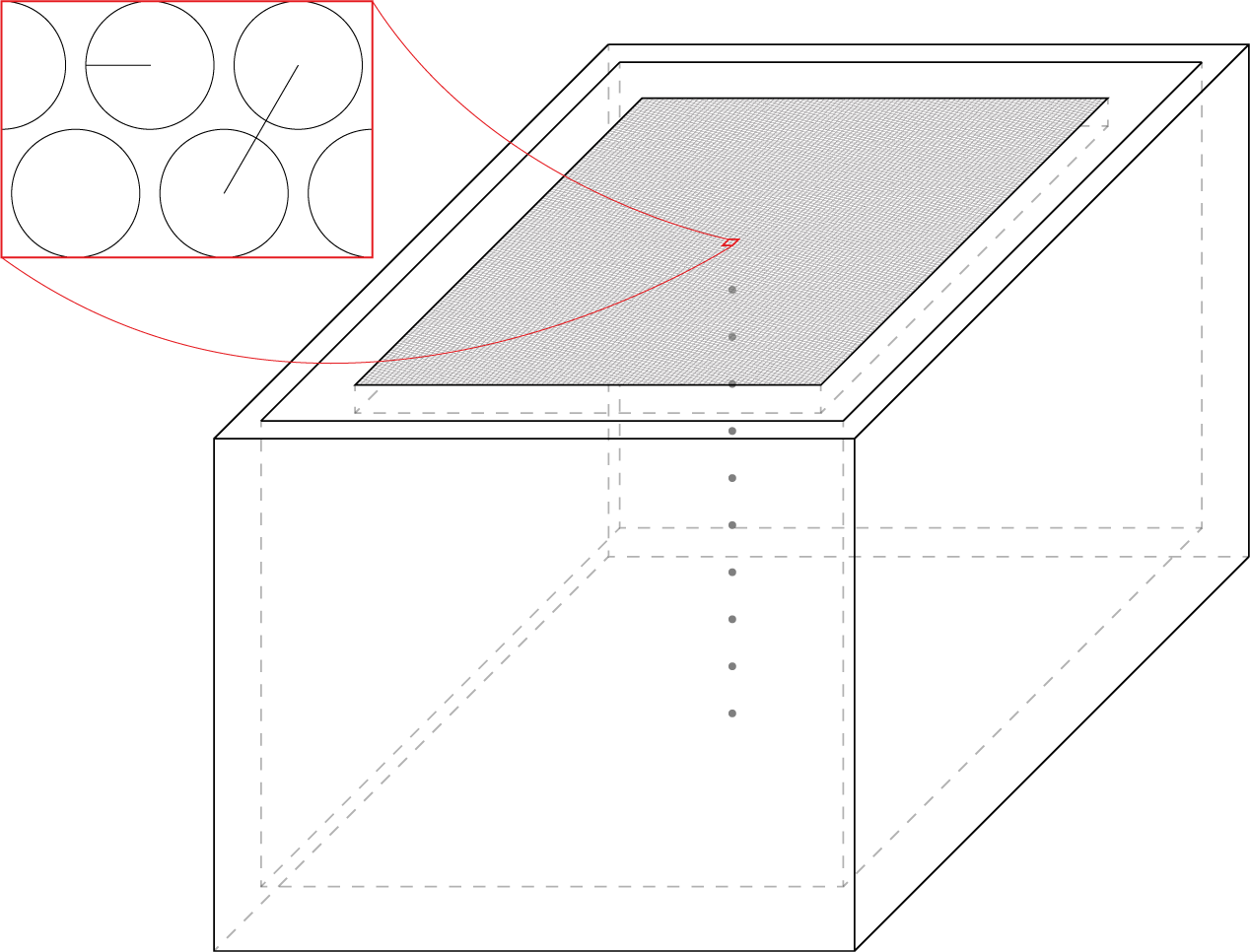}};

\begin{scope}[x={(img.south east)}, y={(img.north west)}]

\draw[{Stealth[length=0.7mm]}-{Stealth[length=0.7mm]}] (0.17,-0.025) -- (0.68,-0.025)
  node[midway, below] {$l_y$};

\draw[{Stealth[length=0.7mm]}-{Stealth[length=0.7mm]}] (0.815,0.8) -- (0.857,0.8);
\node at (0.85,0.82){\small{$2\lambda_\text{w}$}};

\draw[{Stealth[length=0.7mm]}-{Stealth[length=0.7mm]}] (0.45,0.9) -- (0.468,0.9);
\node at (0.43,0.91){\small{$\lambda_\text{w}$}};

\node at (0.10, 0.95) {$r_p$};
\node at (0.185, 0.89) {$d_p$};

\draw[{Stealth[length=0.7mm]}-{Stealth[length=0.7mm]}] (0.14,0.005) -- (0.14,0.54);
\node at (0.12,0.265) {$l_z$};

\draw[{Stealth[length=0.7mm]}-{Stealth[length=0.7mm]}] (0.52,0.568) -- (0.52,0.593);
\node at (0.55,0.58) {\small{$h_\text{in}$}};

\draw[{Stealth[length=0.7mm]}-{Stealth[length=0.7mm]}] (0.8,0.417) -- (0.8,0.445);
\node at (0.83,0.428) {\small{$\lambda_\text{w}$}};

\draw[{Stealth[length=0.7mm]}-{Stealth[length=0.7mm]}] (0.7,-0.025) -- (1.025,0.4);
\node at (0.905,0.2) {$l_x$};

\draw[{Stealth[length=0.7mm]}-{Stealth[length=0.7mm]}] (0.68,0.897) -- (0.71,0.933);
\node at (0.73,0.91) {\small{$2\lambda_\text{w}$}};

\draw[{Stealth[length=0.7mm]}-{Stealth[length=0.7mm]}] (0.564,0.937) -- (0.574,0.95);
\node at (0.58,0.97){\small{$\lambda_\text{w}$}};

\end{scope}
\end{tikzpicture}
\end{minipage}

\caption{Test case of Section \ref{sec:multi_array}. Representation of the three-dimensional view of the geometry employed for the 64-elements simulation. The dots represent the probes at which the solution is computed. Notice that the figure is not scaled.}
\label{fig:3d_view}

\end{figure}

\begin{figure}[H]
    \centering 
    \includegraphics[width=0.65\linewidth]{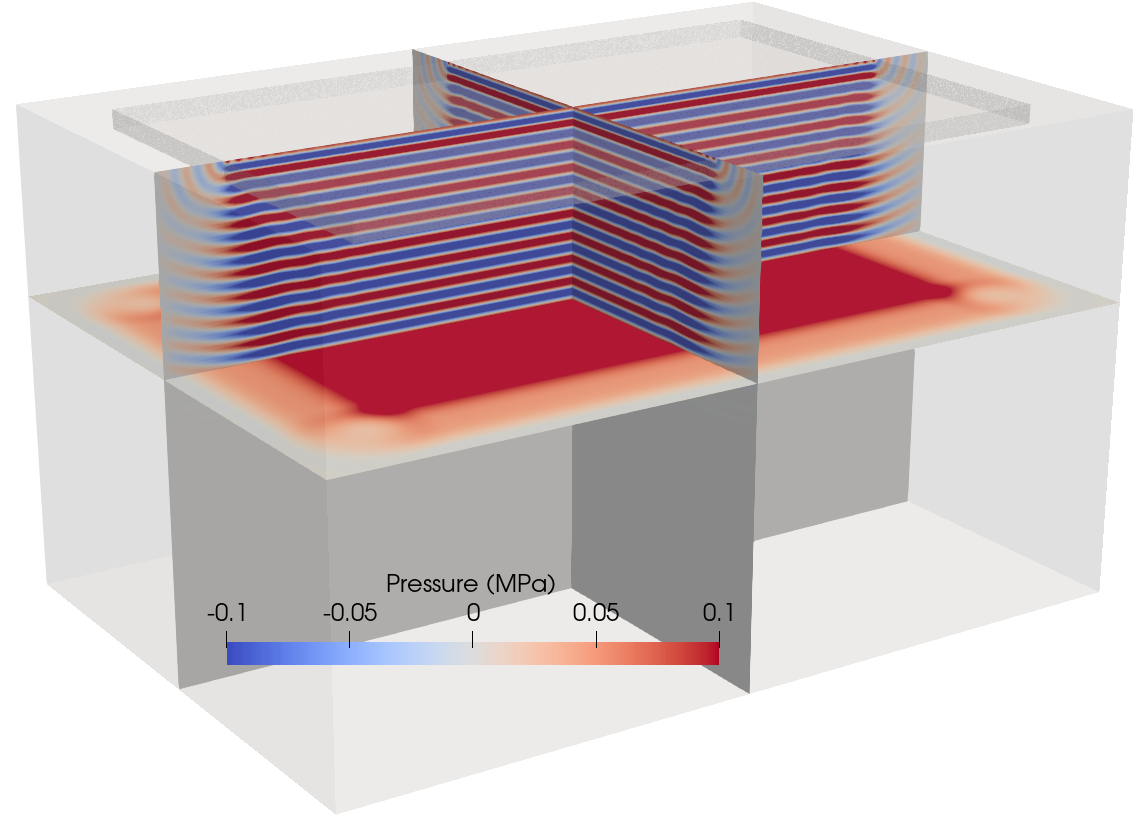}
    \caption*{}

    \caption{Test case of Section \ref{sec:multi_array}. The figure shows the snapshot of the surface containing the PMUTs at time $t=\SI{0.2}{\micro\second}$ with a zoom of the middle region (top), and the solution at time $t=\SI{4}{\micro\second}$ depicted over 3 slices (bottom).}
    \label{fig:result_64_elements}
\end{figure}

\begin{figure}[H]
    \centering 
    \includegraphics[width=0.6\linewidth]{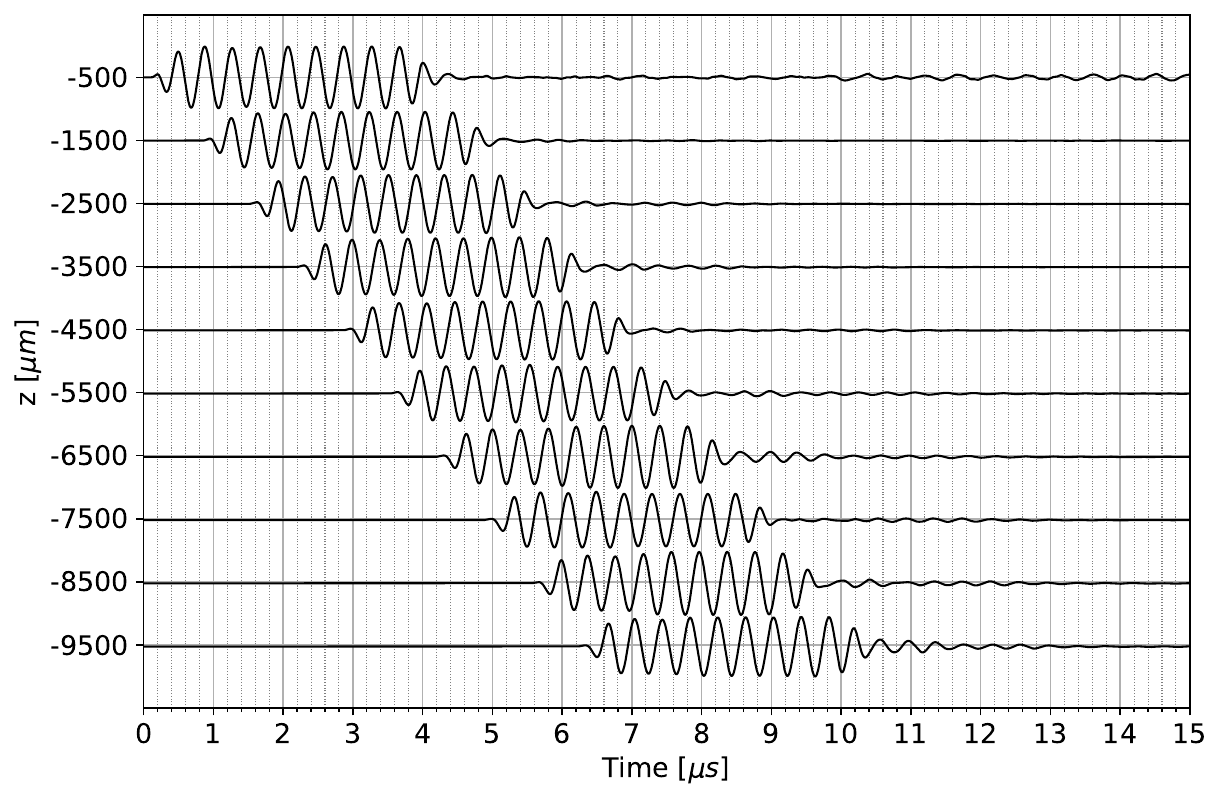}
    \caption{Test case of Section \ref{sec:multi_array}. The figure shows the acoustic pressure computed at ten equispaced points $P=(0,0,z)$, with $z\in\{\SI{-500}{\micro\meter},\,\dots,\,\SI{-9500}{\micro\meter}\}$.}
    \label{fig:result_64_elements_probe}
\end{figure}

Multi-array simulations involve a large number of membranes. This section focuses on the simulation of the 64-element array, cf. Figure \ref{fig:3d_view}. Each PMUT array consists of 184 membranes arranged in two adjacent columns of 92 each, distributed according to the configuration illustrated in Figure \ref{fig:minimal_inner_domain}. This configuration and the physical parameters used in the sequel have been taken from \cite{savoiaST21,savoiaST22}.
%
%
\\\\
The central idea behind this type of simulation, regardless of the number of arrays employed, is to define an acoustic propagation domain with a Cartesian geometry rather than a spherical or cylindrical one. Using Cartesian coordinates significantly reduces the total number of degrees of freedom required for the mesh compared to alternative geometries. Hard-wall conditions are imposed on the top surface $\Gamma_\text{N}$. A Perfectly Matched Layer encloses the computational domain, with absorbing boundary conditions applied along the outer boundary $\Gamma_\text{ABC}$, see equation \eqref{eq:full_problem}. Each array operates independently; PMUTs in different arrays may be driven by distinct voltage excitations. In the following simulation, the applied voltage signal follows the expression in equation (\ref{eq:voltage}), with parameters $f = \SI{2.5}{\mega\hertz}$ and $A = \SI{7.0}{\volt}$, over a total activation time of $T_d = \SI{4}{\micro\second}$. All arrays simultaneously begin and end the transmission phase. The domain consists of two regions described in Section \ref{sec:domain_decomposition}, i.e., an inner and an outer domain. In the inner region, the material properties are defined as $\rho_{\text{in}} = \SI{1140}{\kilo\gram\per\meter\cubed}$ and $c_{\text{in}} = \SI{1130}{\meter\per\second}$. The outer region corresponds to water, characterized by $\rho_{\text{out}} = \SI{1000}{\kilo\gram\per\meter\cubed}$ and $c_{\text{out}} = \SI{1481}{\meter\per\second}$. The inner domain has a thickness of $h_\text{in} = \SI{500}{\micro\meter}$. The perfectly matched layer extends over a distance equal to one wavelength in water, $\lambda_\text{w} = \SI{592.4}{\micro\meter}$, and the overall dimensions of the domain are $l_x = \SI{22874.4}{\micro\meter}$, $l_y = \SI{15570.64}{\micro\meter}$, and $l_z = \SI{12940.4}{\micro\meter}$ along the $x$, $y$, and $z$ directions, respectively. The radius of the PMUT is $r_p = \SI{65}{\micro\meter}$, and the distance between two adjacent PMUT centers, known as pitch, is $d_p = \SI{150}{\micro\meter}$.  Figure \ref{fig:3d_view} contains the geometrical information of the simulation. The time step of the simulation is $\Delta t = \SI{2e-4}{\micro\second}$, and a second-order polynomial degree in space is considered. The final time is set to $T=\SI{15}{\micro\second}$. The number of elements in the mesh is 33.8M, and the degrees of freedom are 262.6M.
Figure~\ref{fig:result_64_elements} illustrates the membrane excitation after $\SI{0.2}{\micro\second}$ and the resulting acoustic wave propagation within the computational domain along the three spatial directions. The pressure field is accurately captured within the computational domain, as shown by the alternating pattern of positive and negative planar pressure waves. Ten probes placed at $P=(0,0,z)$, 
with $z\in\{\SI{-500}{\micro\meter},\,\dots,\,\SI{-9500}{\micro\meter}\}$, allow for evaluating the acoustic pressure signal in both the near and far fields (see Figure \ref{fig:result_64_elements_probe}).
The simulation ran on a high-performance cluster consisting of 40 CPU compute nodes. Each CPU node is configured with two AMD EPYC 7413 processors, each featuring 24 cores, yielding a total of 96 threads per node. These nodes are provisioned with 512 GB of RAM each, providing substantial memory capacity to support large-scale simulations and data-intensive computations within a parallel computing framework. By utilizing 21 nodes, and 1935 processes, the simulation required approximately 12 hours to complete.

\begin{figure}[H]
    \centering
    
    \begin{minipage}{0.62\linewidth}
        \centering
        \begin{tikzpicture}[remember picture]
            \node[anchor=south west, inner sep=0] (img) 
                {\includegraphics[width=\linewidth]{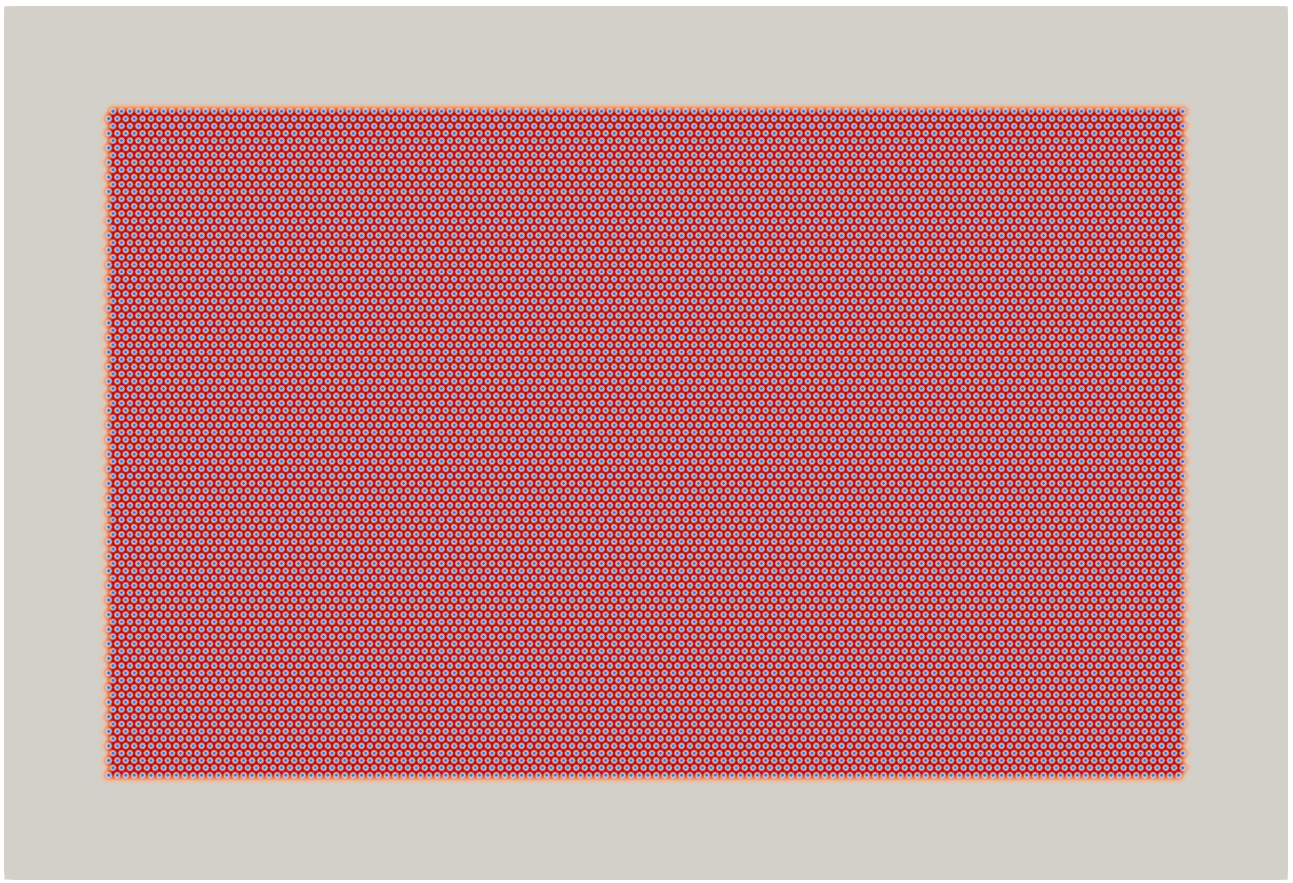}};
            
            \draw[black, thick]
                (4.91,3.5) rectangle (5.91,4.485);

            \coordinate (LTzoom) at (5.91,4.485);
            \coordinate (LBzoom) at (5.91,3.5);
        \end{tikzpicture}
    \end{minipage}%
    \hfill
    \begin{minipage}{0.36\linewidth}
        \centering
        \raisebox{0.5\height - 0.5\totalheight}{%
            \begin{tikzpicture}[remember picture]
                \node[anchor=south west, inner sep=0] (img2)
                    {\includegraphics[width=\linewidth]{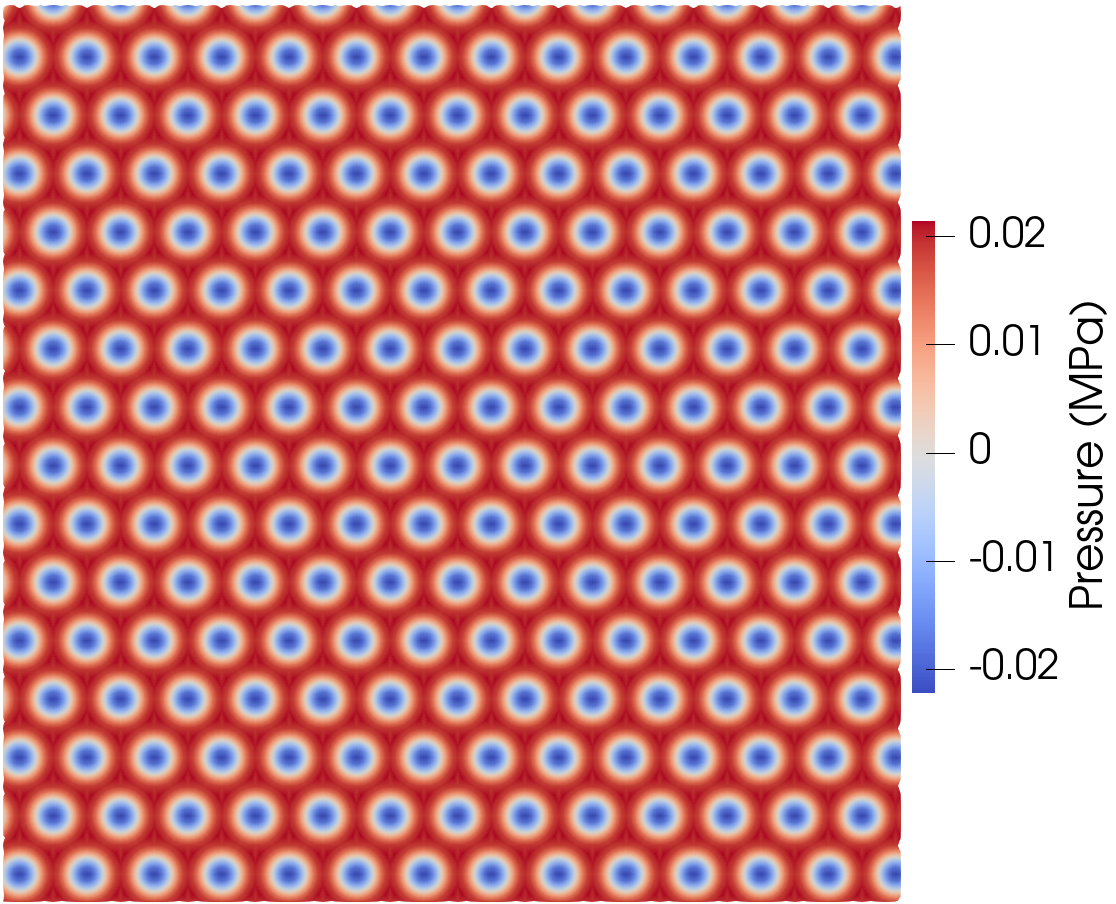}};
                
                \draw[black, thick]
                    (0.03,0.03) rectangle (5.1,5.1);

                \coordinate (RTzoom) at (0.03,5.1);
                \coordinate (RBzoom) at (0.03,0.03);
            \end{tikzpicture}%
        }
    \end{minipage}

    \begin{tikzpicture}[overlay, remember picture]
        \draw[thick]
            (LTzoom) -- (RTzoom);
        \draw[thick]
            (LBzoom) -- (RBzoom);    
    \end{tikzpicture}

    \caption*{}
\end{figure}

\subsection{Performance comparisons and scalability test}\label{sec:scalability}
A comparative performance analysis was carried out to evaluate the effectiveness of the DG interface optimization described in Section \ref{sec:dg_optimization} and to quantify the computational and communication overhead associated with the evaluation of the integrals over $\Gamma_\text{PMUT}$ in equations (\ref{eq:f_i}) and (\ref{eq:g_m}), as discussed in Section \ref{sec:domain_decomposition}. The first test case employed a computational domain consisting of a 20-element array distributed across 315 processes. Two configurations were investigated: one without DG interface optimization and one with the optimization enabled. In the absence of DG optimization, the setup of $318698$ interfaces required approximately $\SI{160000}{\second}$. By contrast, enabling DG interface optimization reduced the setup time to approximately $\SI{630}{\second}$. This represents a reduction of more than two orders of magnitude in interface setup time, underscoring the critical importance of DG optimization for improving scalability and significantly reducing preprocessing overhead. The second test case was performed on a computational domain comprising a single-element array partitioned across 96 processes. In the first configuration, an arbitrary mesh partitioning was employed, with no constraint enforcing the association of an entire PMUT with a single process. Under this configuration, the solution of equation (\ref{eq:semi_PMUT}) incurred an inter-process communication cost of approximately $\SI{0.167}{\second}$ per timestep, arising solely from the exchange of local contribution to the PMUTs solution on each process with all other processes. This communication overhead substantially exceeded the pure computational cost, thereby indicating poor parallel efficiency.

In the second configuration, the tailored partitioning strategy described in Section \ref{sec:domain_decomposition} was adopted. This approach resulted in an average synchronization waiting time of approximately $\SI{0.0008}{\second}$ per PMUT per timestep for each process, corresponding to a reduction of more than two orders of magnitude. It is worth noting that, in the former configuration, increasing either the number of membranes or the number of processes leads to a corresponding increase in communication time, whereas in the latter configuration, the communication cost remains, on average, constant. This behavior is further illustrated by the scalability analysis presented in the following.

Overall, these results demonstrate that DG interface optimization is essential for reducing preprocessing time in large-scale simulations, while PMUTs performance is highly sensitive to communication patterns and domain decomposition strategies. Inappropriate partitioning can result in communication-dominated runtimes, thereby severely limiting parallel efficiency. These findings underscore the necessity of balanced mesh partitioning and efficient data exchange mechanisms to achieve scalable and efficient parallel performance.

\begin{figure}[H]
    \begin{minipage}{0.3\linewidth}
        \centering
        \includegraphics[width=1.65\linewidth]{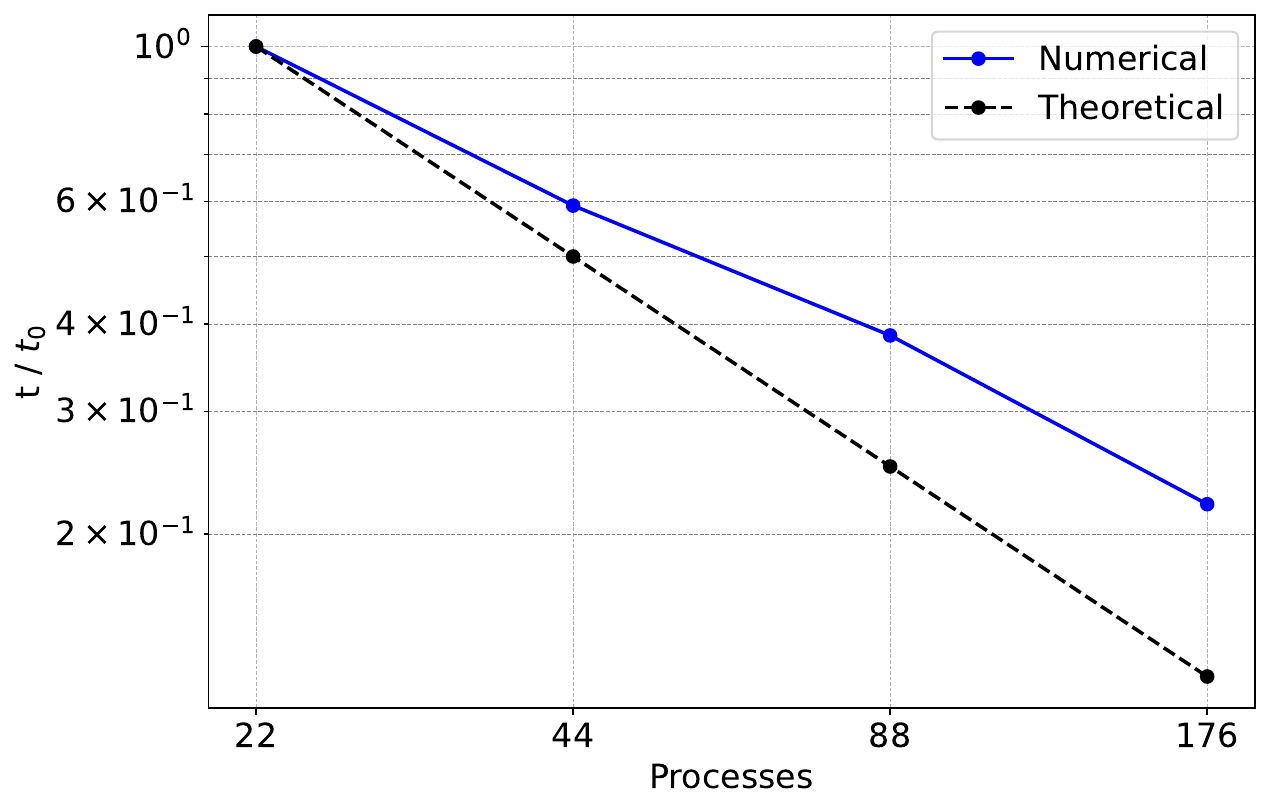}
    \end{minipage}\hspace{11em}
    \begin{minipage}{0.3\linewidth}
        \centering
        \includegraphics[width=1.55\linewidth]{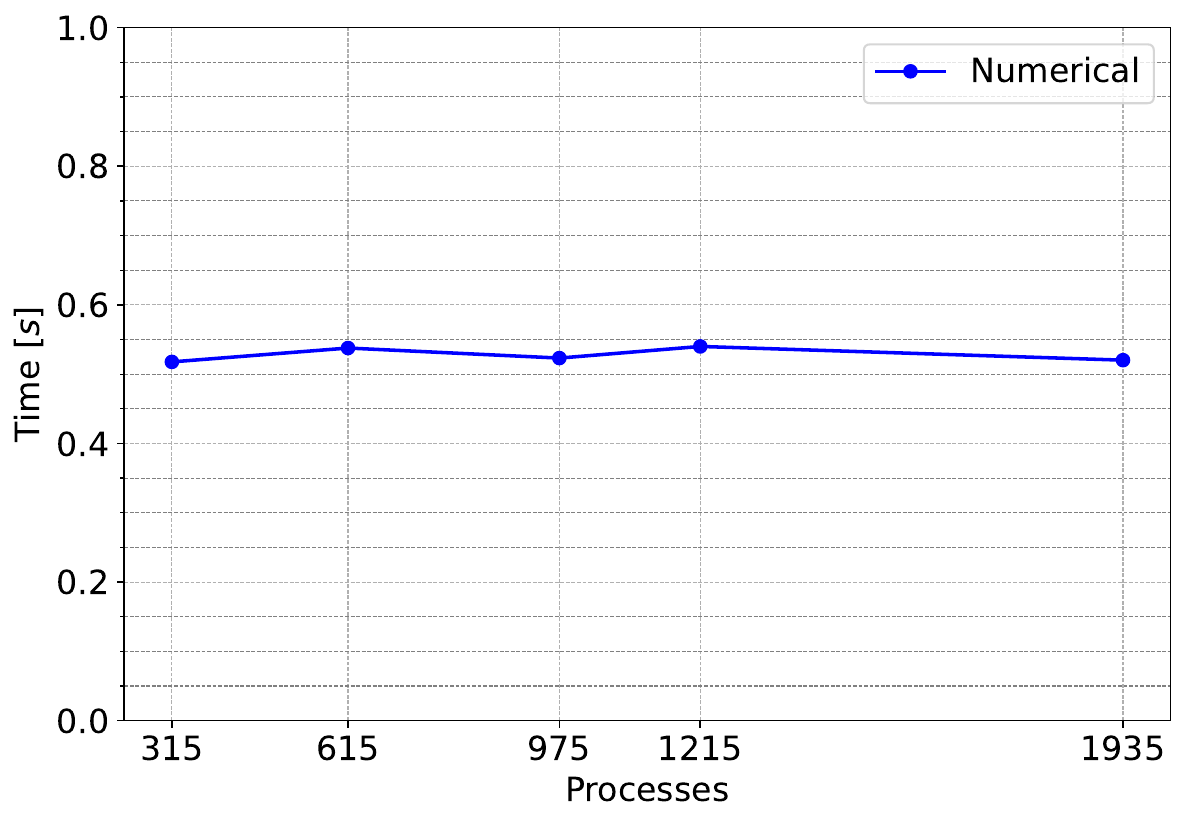}
    \end{minipage}
    \caption{Test case of Section \ref{sec:scalability}. The figure shows the theoretical trend (black) and the real one (blue) of the strong scalability test conducted across four different simulations (left). The value on the $y-$axes is normalized with respect to the time $t_0$ of the first simulation. On the right, the weak scalability test performed on five different simulations.}
    \label{fig:scalability}
\end{figure}

Scalability refers to the ability of a computational system to maintain or improve performance as workload and resources increase. In high-performance and parallel computing, it depends on the efficient use of processors, memory, and communication, requiring well-designed parallel algorithms, effective memory management, and controlled communication overhead to ensure competitive performance.
In the specific context of the problem under investigation, the main constraints that influence scalability are the ones discussed in Section \ref{sec:domain_decomposition} and can be categorized into three main factors:

\begin{enumerate}[label=\roman*)]
\item PMUTs blocks: the scalability of the system is directly affected by the configuration and number of PMUT array blocks. As the array size increases, the computational load associated with the membrane becomes dominant. Efficient data partitioning and optimized communication between array elements are therefore essential to maintain scalability in large configurations.

\item Acoustic propagation domain: the computational complexity of simulating the acoustic propagation domain scales with both the spatial resolution and the physical model. High-resolution domains require finer discretization, leading to increased memory usage and processing demands. Load balancing and domain decomposition strategies are of paramount importance to ensure that the computations remain scalable across multiple processes. Moreover, simulations with the PML enlarge the DOFs and hence require an additional complexity to be taken into account.

\item DG interfaces: they represent another critical scalability constraint. These interfaces govern the communication and numerical flux exchanges between subdomains. As the number of hexahedral elements increases, so does the volume of interfacial data transfer, potentially leading to communication bottlenecks. Maintaining scalability in this context requires efficient inter-process communication schemes and optimized data structures that ease rapid access and exchange of boundary information.
\end{enumerate}
Overcoming these challenges requires a combination of algorithmic optimization, advanced parallelization strategies, and system design to ensure that performance scales effectively with increasing computational resources. A strong scalability test is performed to assess the validity of the approach discussed in the previous sections. Four simulations are conducted within the same computational domain, exploiting the cluster of the Mathematical Department of Politecnico di Milano (DMAT), whose characteristics are reported in the previous section. The number of processes selected is $22, 44, 88$ and $176$, respectively. In this way, communication and load balancing are performed both among processes within a single node and between different nodes. Figure \ref{fig:scalability} (left) shows the comparison between the theoretical and normalized real-time resources required in the four cases. The plot in logarithmic scale shows a very good trend compared to the expected behavior. Moreover, a weak scalability test was assessed over three different simulations proportionally increasing both the number of degrees of freedom and the number of processes, thereby preserving a fixed ratio of approximately 140000 DOFs per process. Figure \ref{fig:scalability} (right) shows an almost ideal weak scalability, by holding the execution timestep constant.  

\section{Conclusions and future developements}\label{sec:conclusion}
This work presented a comprehensive methodology for the efficient numerical simulation of large-scale Piezoelectric Micromachined Ultrasonic Transducer arrays, based on a hybrid approach consisting of model order reduction techniques and a high-performance Discontinuous Galerkin Spectral Element framework implemented within the SPectral Elements in Elastodynamics with Discontinuous Galerkin (SPEED) environment. By integrating the mechanical modal behavior of individual PMUTs into an acoustic domain discretized through high-order elements, the proposed approach achieves a balance between accuracy and computational efficiency. The resulting model successfully captures both transmission and reception phenomena, enabling scalable and high-fidelity simulations of complex array configurations. The implemented domain decomposition and optimization strategies proved effective in enhancing computational performance, ensuring that the framework remains suitable for massively parallel simulations involving thousands of transducers. The numerical results confirm the robustness and scalability of the proposed approach, making it a promising tool for the design and optimization of next-generation PMUT-based ultrasonic systems. Beyond the present developments, future work will focus on extending the current framework to account for thermoacoustic effects, which are increasingly relevant as PMUTs operate at higher power densities and within miniaturized environments. Incorporating thermo-mechanical coupling and heat generation due to dielectric and mechanical losses will enable the study of temperature-induced frequency shifts, performance degradation, and thermal noise in PMUT arrays. Moreover, the integration of nonlinear thermoelastic models could provide deeper insight into self-heating phenomena and their influence on both emission and reception behavior. Such extensions would pave the way toward multiphysics simulations capable of predicting the full electromechanical–acoustic–thermal response of PMUT systems under realistic operational conditions. Ultimately, these developments will support the design of more reliable, energy-efficient, and thermally stable ultrasonic microsystems, broadening their applicability in demanding environments such as biomedical diagnostics, industrial sensing, and harsh-condition monitoring.

\section*{Acknowledgments}
This work is part of the \textit{PROgramma per la realizzazione di dispositivi ad Ultrasuoni e loro promozione Democratica} (PROUD) project funded by the Italian Ministry of Enterprises and Made in Italy (MIMIT): Prog. n.: F/310211/01-05/X56 (Decree n.4187, 19-12-2023). The authors would like to thank STMicroelectronics for their valuable support and contributions within the framework of the PROUD project. This work also benefited from helpful discussions and technical support from all project partners. 
P.F. Antonietti, M.G. Garroni, I. Mazzieri, and N. Prolini are members of INdAM-GNCS group. The present research is part of the activities of ``Dipartimento di Eccellenza 2023-2027''.
The simulations discussed in this work were, in part, performed on the HPC Cluster of the Department of Mathematics of Politecnico di Milano, which was funded by MUR grant ``Dipartimento di Eccellenza 2023-2027''.


\end{document}